\RequirePackage{pdf14}

\PassOptionsToPackage{usenames,dvipsnames}{xcolor}
\PassOptionsToPackage{breaklinks=true}{hyperref}

\documentclass[acmtog]{acmart}

\usepackage{booktabs} 
\usepackage{alltt}
\usepackage{wrapfig}

\usepackage[l2tabu,orthodox]{nag}

\usepackage[T1]{fontenc}
\usepackage[utf8]{inputenc}

\usepackage{amsthm}
\usepackage{amssymb}
\usepackage{mathtools} 
\usepackage{textcomp}
\usepackage[binary-units]{siunitx} 
\usepackage{physics} 
\usepackage{stmaryrd}
\usepackage{mathrsfs}

\usepackage{ifthen}
\usepackage{float} 
\usepackage{hyperref}
\usepackage{cleveref}
\usepackage{csquotes}
\usepackage{calc}
\usepackage[usenames,dvipsnames]{xcolor}
\usepackage{enumitem}

\usepackage{array}
\usepackage{booktabs}
\usepackage{multirow}

\usepackage[percent]{overpic}

\newcommand{\specialcell}[2][c]{%
	\begin{tabular}[#1]{@{}c@{}}#2\end{tabular}}

\newcommand\restr[2]{{
  \left.\kern-\nulldelimiterspace 
  #1 
  \vphantom{\big|} 
  \right|_{#2} 
  }}

\newcommand{\cM}{\mathcal{M}}
\newcommand{\cH}{\mathcal{H}}
\newcommand{\cP}{\mathcal{P}}
\newcommand{\cQ}{\mathcal{Q}}
\newcommand{\cS}{\mathcal{S}}

\newcommand{\param}[1]{{\widehat{#1}}}
\newcommand{\mat}[1]{{\mathbf{#1}}}
\newcommand{\vect}[1]{{\boldsymbol{\mathbf{#1}}}}

\newcommand{\numK}{k}
\newcommand{\kernel}{\psi}

\newcommand{\basis}{\phi}
\newcommand{\reals}{\mathbb{R}}
\newcommand{\threed}{\mathbb{R}^3}

\newcommand{\bx}{\mathbf{x}}
\newcommand{\bv}{\mathbf{v}}

\newcommand{\bu}{\mathbf{u}}
\newcommand{\bg}{\mathbf{g}}
\newcommand{\kc}{\mathbf{z}}
\newcommand{\dD}{\mathrm{D}}
\newcommand{\monom}{q}

\newcommand{\T}{^\mathsf{T}}

\newcommand{\pdiff}[2]{\frac{\partial #2}{\partial #1}}

\usepackage[ruled]{algorithm2e} 

\SetAlFnt{\small}
\SetAlCapFnt{\small}
\SetAlCapNameFnt{\small}
\SetAlCapHSkip{0pt}
\IncMargin{-\parindent}

\setcopyright{none}

\setcitestyle{square}

\begin{document}

\title{Poly-Spline Finite Element Method}

\author{Teseo Schneider}
\affiliation{
  \institution{New York University}
  \country{USA}
}
\email{}
\author{J\'er\'emie Dumas}
\affiliation{
  \institution{New York University, nTopology}
  \country{USA}
}
\author{Xifeng Gao}
\affiliation{
    \institution{New York University, Florida State University}
  	\country{USA}
}
\author{Mario Botsch}
\affiliation{
    \institution{Bielefeld University}
    \country{Germany}
}
\author{Daniele Panozzo}
\affiliation{
  \institution{New York University}
  \country{USA}
}
\author{Denis Zorin}
\affiliation{
  \institution{New York University}
  \country{USA}
}

\renewcommand\shortauthors{Schneider, T. et al.}

\begin{abstract}
We introduce an integrated meshing and finite element method pipeline enabling solution of partial differential equations in the volume enclosed by a boundary representation. We construct a hybrid hexahedral-dominant mesh, which contains a small number of star-shaped polyhedra, and build a set of high-order bases on its elements, combining triquadratic B-splines, triquadratic hexahedra, and harmonic elements. We demonstrate that our approach converges cubically under refinement, while requiring around 50\% of the degrees of freedom than a similarly dense hexahedral mesh composed of triquadratic hexahedra. We validate our approach solving Poisson's equation on a large collection of models, which are automatically processed by our algorithm, only requiring the user to provide boundary conditions on their surface.

\end{abstract}


\graphicspath{{./figs}}

\begin{teaserfigure}
	\centering
	\includegraphics[width=\linewidth]{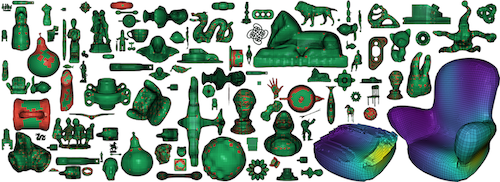}
	\caption{A selection of the automatically generated pure hexahedral and hexahedral-dominant meshes in our test set. The colors denote the type of basis used. In the bottom-right, we show the result of a Poisson problem solved over a hex-dominant, polyhedral mesh.}
	\label{fig:teaser}
\end{teaserfigure}

\maketitle


\section{Introduction}
\label{sec:intro}

The numerical solution of partial differential equations is ubiquitous in computer graphics and engineering applications, ranging from the computation of UV maps and skinning weights, to the simulation of elastic deformations, fluids, and light scattering.

The finite element method (FEM) is the most commonly used discretization of PDEs, especially in the context of structural and thermal analysis, due to its generality and rich selection of off-the-shelf commercial implementations. Ideally, a PDE solver should be a ``black box'': the user provides as input the domain boundary, boundary conditions, and the governing equations, and the code returns an evaluator that can compute the value of the solution at any point of the input domain. This is surprisingly far from being the case for all existing open-source or commercial software, despite the research efforts in this direction and the large academic and industrial interest.

To a large extent, this is due to treating meshing and FEM basis construction as two disjoint problems. The FEM basis construction may make a seemingly innocuous assumption (e.g., on the geometry of elements), that lead to exceedingly difficult requirements for meshing software. For example, commonly used bases for tetrahedra are sensitive to the tetrahedron shape, so tetrahedral mesh generators have to guarantee good element shape everywhere:  a difficult task which, for some surfaces, does not have a fully satisfactory solution. 
Alternatively, if few assumptions are made on mesh generation (e.g., one can use elements that work on arbitrary polyhedral domains), the basis and stiffness matrix constructions can become very expensive. 

This state of matters presents a fundamental problem for applications that require fully automatic, robust processing of large collections of meshes of varying sizes, an increasingly common situation as large collections of geometric data become available. Most importantly, this situation arises in the context of machine learning on geometric and physical data, where a neural network could be trained using large numbers of simulations and used to compute efficiently an approximated solution \cite{Chen:Neural:2018,kostrikov2017surface}. Similarly, shape optimization problems often require solving PDEs in the inner optimization loop on a constantly changing domain \cite{Panetta:2015}.

\paragraph{Overview.}
We propose an integrated pipeline, considering meshing and element design as a single challenge: we make the tradeoff between mesh quality and element complexity/cost  \emph{local}, instead of making an a priori decision for the whole pipeline. We generate high quality, simple, and regularly arranged elements for most of the volume of the shape, with more complex and poor quality  polyhedral shapes filling the remaining gaps \cite{Sokolov2016,Gao:2017}. Our idea is to match each element to a basis construction, with well-shaped elements getting the simplest and most efficient basis functions and with complex polyhedral element formulations used only when necessary to handle the transitions between regular regions, which are the ones that are topologically and geometrically more challenging. 

A spline basis on a regular lattice has major advantages over traditional FEM elements, since it has the potential to be \emph{both} accurate and efficient: it has a single degree of freedom per element, except at the boundary, yet, it has full approximation power corresponding to the degree of the spline. This observation is one of the foundations of \emph{isogeometric analysis} in 3D \cite{hughes:2005:isogeometric,Cottrell:2009}. 
Unfortunately, it is easy to define and implement only for fully regular grids, which is not practical for most input geometries. The next best thing are spline bases on
\emph{pure hexahedral} meshes: while smooth constructions for polar configurations exist \cite{Toshniwal:2017}, a solution applicable to general hexahedral meshes whose interior singular curves meet is still elusive, restricting this construction to simple shapes. Padded hexahedral-meshes \cite{Marechal:2009} are necessary to ensure a good boundary approximation for both regular and polycube~\cite{Tarini:2004} hexahedral meshing methods, but they unfortunately cannot be used by these constructions since their interior curve singularities meet in the padding layer.

We propose a hybrid construction that sidesteps these limitations: we use spline elements only on fully regular regions, and fill the elements that are touching singular edges, or that are not hexahedra, with local constructions (harmonic elements for polyhedra, triquadratic polynomial elements for hexahedra). This construction  further relaxes requirements for meshing, since it works on general hexahedral meshes (without any restriction on their singularity structure) but also directly supports \emph{hex-dominant} meshes, which can be robustly generated with modern field-aligned methods \cite{Gao:2017,Sokolov2016}. These meshes consist mostly of well-shaped hexahedra with locally regular mesh structure,  but also contain other general polyhedra. Our construction takes advantage of this high regularity, adding a negligible overhead over the spline FEM basis only for the sparse set of non-regular elements.

We demonstrate that our proposed \emph{Poly-Spline} FEM retains, to a large extent, both the approximation and performance benefits of splines, at the cost of the increasing basis construction complexity, and at the same time, works for a class of meshes that can be robustly generated for most shapes with existing meshing algorithms.

Our method exhibits cubic convergence on a large data set, for a degree of freedom budget comparable to trilinear hexahedral elements, which have only quadratic convergence. To the best of our knowledge, this paper is the first FEM method exploiting the advantages of spline basis that has been validated on a large collection of complex geometries.

\section{Related Work}
\label{sec:related}

When numerically solving PDEs using the finite element method, one has to discretize the spatial domain into finite elements and define shape functions on these elements. Since shape functions, element types, and mesh generation are closely related, we discuss the relevant approaches in tandem.

For complex spatial domains, the discretization is  frequently based on the Delaunay triangulation \cite{Shewchuk:Tri:1996} or Delaunay tetrahedrization \cite{Si:TetGen:2015}, respectively, since those tessellations can be computed in a robust and automatic manner. Due to their simplicity and efficiency, linear shape functions over triangular or tetrahedral elements are often the default choice for graphics applications~\cite{hughes:fem}, although they are known to suffer from locking for stiff PDEs, such as nearly incompressible elastic materials~\cite{hughes:fem}.

This locking problem can be avoided by using bilinear quadrangular or trilinear hexahedral elements ($Q_1$ elements), which have the additional advantage of yielding a higher accuracy for a given  number of elements \cite{FEAD1992,Benzley95}. Triquadratic hexahedral elements ($Q_2$) provide even higher accuracy and faster convergence under mesh refinement (cubic converge in $L_2$-norm for $Q_2$ vs.\ quadratic converge for $Q_1$), but their larger number of degrees of freedom (8 vs. 27) leads to high memory consumption and computational cost.

The main idea of isogeometric analysis (IGA) \cite{hughes:2005:isogeometric,Cottrell:2009,Engvall:2017:IGA} is to employ the same spline basis for defining the CAD geometry as well as for performing numerical analysis. Using quadratic splines on hexahedral elements results in the same cubic convergence order as $Q_2$ elements, but at the much lower cost of \emph{one} degree of freedom per element (comparable to $Q_1$ elements). This efficiency, however, comes at the price of a very complex implementation for non-regular hexahedral meshes. Moreover, generating IGA-compatible meshes from a given general boundary surface is still an open problem \cite{aigner:2009:swept,martin:2010:volumetric,li:2013:surface}.

Concurrent work~\cite{Wei2018} introduces a construction that can handle irregular pure hex meshes, with tensor-product cubic splines used on regular parts.  However, we focus on handling general polygonal meshes and we use quadratic splines (note that, our approach can be easily extended to cubic polynomials if desired).

A standard method for volumetric mesh generation is through hierarchical subdivision of an initial regular hexahedral mesh, leading to so-called octree meshes \cite{marechal2009,Ito2009,Zhang2013}. The T-junctions resulting from adaptive subdivision can be handled by using T-splines \cite{Sederberg:TSL:2004,Veiga:TIGA:2011} as shape functions. While this meshing approach is very robust, it has problems representing geometric features that are not aligned with the principal axes.

Even when giving up splines or T-splines for standard $Q_1$/$Q_2$ elements, the automatic generation of the required hexahedral meshes is problematic. Despite the progress made in this field over the last decade, \emph{automatically} generating pure hexahedral meshes that (i) have sufficient element quality, (ii) are not too dense, and (iii) align to geometric features is still unsolved. Early methods based on paving or sweeping \cite{Owen:2000,Yamakawa:2003,Staten:UPP:2005,Shepherd2008} require complicated handling of special cases and generate too many singularities. Polycube methods~\cite{HexmeshSGP2011,li:2013:surface,livesu2013polycut,huang2014,fang2016all,Fu:2016}, field-aligned methods \cite{NieserSGP11,Huang2011,Li2012,Jiang2013}, and the surface foliation method \cite{LEI2017758} are interesting research venues, but they are currently not robust enough and often fail to produce a valid mesh.

However, if the strict requirement of producing hexahedral elements only is relaxed, field-aligned methods \cite{Sokolov2016,Gao:2017} can robustly and automatically create hex-dominant polyhedral meshes, that is, meshes consisting of mostly, but not exclusively, of hexahedral elements. The idea is to build local volumetric parameterizations aligned with a specified directional field,  and constructing the mesh from traced isolines of that parameterization, inserting general polyhedra if necessary.  Their drawback is that the resulting \emph{hex-dominant} meshes, are not directly supported by most FEM codes.

One option is to split these general polyhedra into standard elements, leading to a mixed FEM formulation. For instance, the field-aligned meshing of Sokolov et al.~\shortcite{Sokolov2016} extract meshes that are composed of hexahedra, tetrahedra, prisms, and pyramids. However, the quality of those split elements is hard to control in general. An interesting alternative is to avoid the splitting of polyhedra and instead incorporate them into the simulation, for instance though mimetic finite differences~\cite{Lipnikov:MFD:2014}, the virtual element method~\cite{Veiga:VEM:2013}, or polyhedral finite elements~\cite{Manzini:NPPFEM:2014}. The latter employ generalized barycentric coordinates as shape functions, such as mean value coordinates \cite{Floater:2005:MVC,Ju:2005:MVC}, harmonic coordinates \cite{Joshi:2007:HCF}, or minimum entropy coordinates \cite{Hormann:2008:MEC}. From those options, harmonic coordinates seem most suitable since they generalize both linear tetrahedra and trilinear hexahedra to general non-convex polyhedra \cite{Martin:PFE:2008,Bishop:2014}. While avoiding splitting or remeshing hex-dominant meshes, the major drawback of polyhedral elements is the high cost for computing and integrating their shape functions.

In the above methods the meshing stage either severely restricts admissible shape functions, or the element type puts (too) strong requirements on the meshing. In contrast, we use the most efficient elements where possible and the most flexible elements where required, which enables the use of robust and automatic hex-dominant mesh generation.

\section{Algorithm Overview}
\label{sec:overview}

In this section, we introduce the main definitions we use in our algorithm description, and outline the structure of the algorithm. We refer to \Cref{app:fem} for a brief introduction to the finite element method and the setup of our mathematical notation.

\begin{figure}
  \centering
  \includegraphics[width=0.4\linewidth]{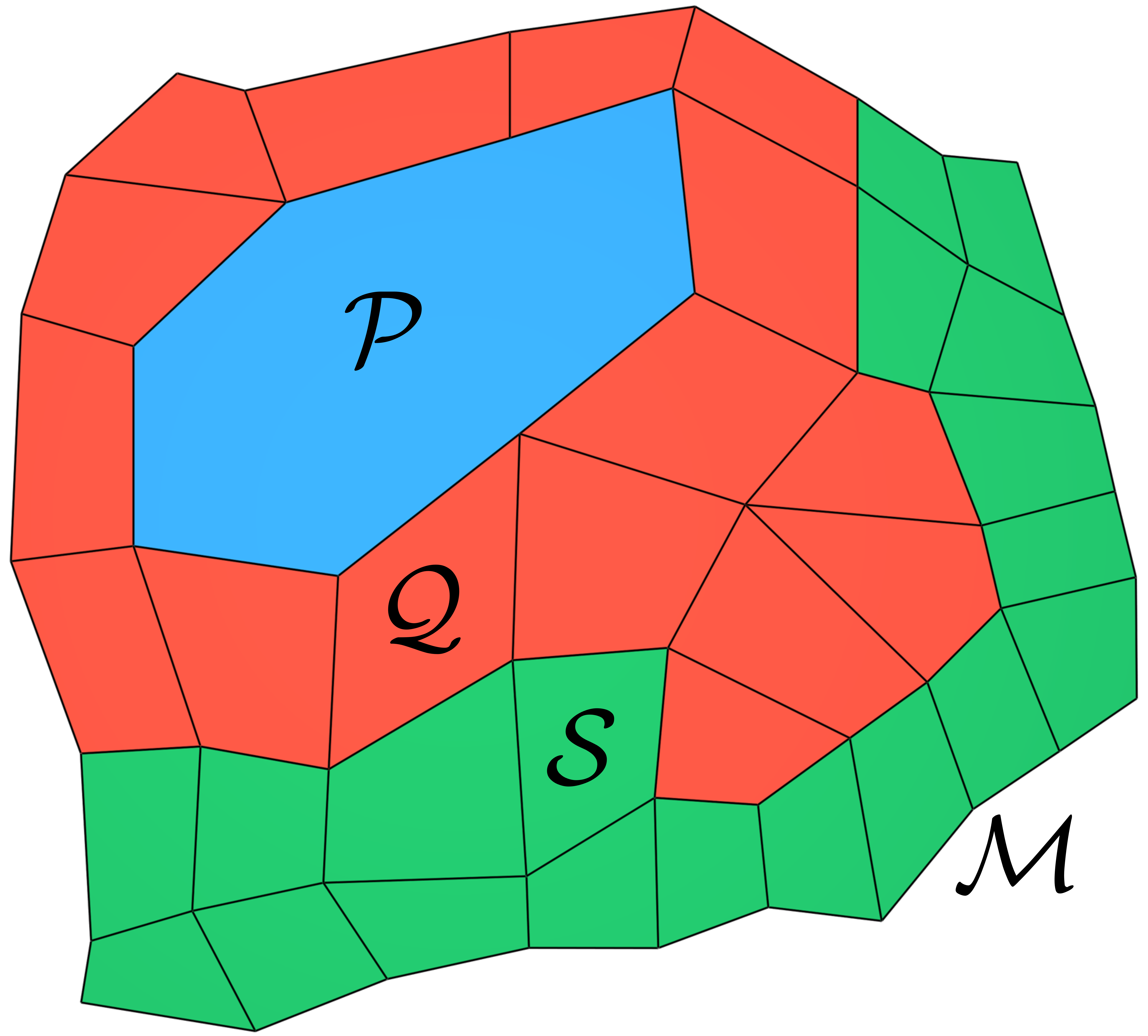}
  \caption{Complexes involved in our construction. In green we show $\cS$, in red $\cQ$, and in blue $\cP$.}
\label{fig:complex}
\end{figure}

\paragraph{Input complex and subcomplexes.}
The input to our algorithm is a 3D polyhedral complex $\cM$, with vertices $\bv_i\in \threed$, $i=1,\ldots,N_V$, consisting of polyhedral cells $C_i$, $i=1,\ldots,N_C$, most of which are hexahedra. \Cref{fig:complex} shows a two-dimensional example of such complex.
The edges, faces, and cells of the mesh are defined combinatorially, that is, edges are defined by pairs of vertices, faces by sequences of edges, and cells by closed surface meshes formed by faces.
We assume that 3D positions of vertices are also provided as input and that $\cM$ is three-manifold, i.e., that there is a way to identify vertices, edges, faces,
and cells with points, curves, surface patches and simple volumes, such that their union is a three-manifold subset of $\threed$.


We assume that for any hexahedron there is at most one non-hexahedral cell sharing one of its faces or edges, which can be achieved by refinement. We also assume that no two polyhedral cells are adjacent, and that no polyhedron touches the boundary, which can also be achieved by merging polyhedral cells and/or refinement.
This preprocessing step (i.e., one step of uniform refinement) is discussed in Section~\ref{sec:refinement}. As a consequence of our refinement, \emph{all faces of $\cM$ are quadrilateral}.

One of the difficulties of using general polyhedral meshes for basis constructions is that, unlike the case of, 
for example, pure tetrahedral meshes, there is no natural way to realize all elements of the mesh in 3D just from vertex positions (e.g., for a tetrahedral mesh, linear interpolation for faces and cells is natural). This requires constructing bases on an explicitly defined parametric domain associated with the input complex.

For this purpose, we define a certain number of complexes related to the original complex $\cM$ (Figure~\ref{fig:complex}).
There are two goals for introducing these: defining the parametric domain for the basis, and defining the \emph{geometric map}, which specifies how the complex is realized in three-dimensional \emph{physical space}.

\begin{itemize}[leftmargin=*]
\item $\cH \subseteq \cM$ is the hexahedral part of $\cM$, consisting of hexahedra $H$. 
\item $\cP = \cM\setminus\cH$ is the non-hexahedral part of $\cM$, consisting of polyhedra $P$.
\item $\cS \subseteq \cH$ is the complex consisting of spline-compatible hexahedra
  $S$ defined in \Cref{sec:spline-compatible}. 
 \item $\cQ = \cH \setminus \cS$ is the spline-incompatible pure-hexahedral part of $\cM$.
\end{itemize}
\noindent
Note that the sub-complexes of $\cM$ are nested:  $\cS \subseteq \cH \subseteq \cM$.

In the context of finite elements, the distinction between \emph{parametric space} and \emph{physical space} is critical:
the bases on the hexahedral part of the mesh are defined in terms of parametric space coordinates, where all hexahedra are unit cubes; this makes it possible to define simple, accurate, and efficient bases.
However, the derivatives in the PDE are taken with respect to physical space variables, and the unknown functions are naturally defined on the physical space. Remapping these functions to the parametric space is necessary to discretize the PDE using our basis. 
We define parametric domains $\param{\cM}$, $\param{\cH}$, $\param{\cS}$, and $\param{\cQ}$ corresponding to ${\cM}$, ${\cH}$, ${\cS}$, and ${\cQ}$, respectively.
 $\param{\cH}$ consists of unit cubes $\param{H}$, one per hexahedron $H$ with corresponding faces identified, and $\param{\cS}$ and $\param{\cQ}$ are its subcomplexes.  The complete parametric space $\param{\cM}$ is obtained by 
 adding a set of polyhedra for $\cP$, defined using the geometric map 
 as described below. For polyhedra, physical and parametric space coincide.
 

\paragraph{Geometric map and complex embedding.}

\begin{figure}
\centering
\includegraphics[width=0.8\linewidth]{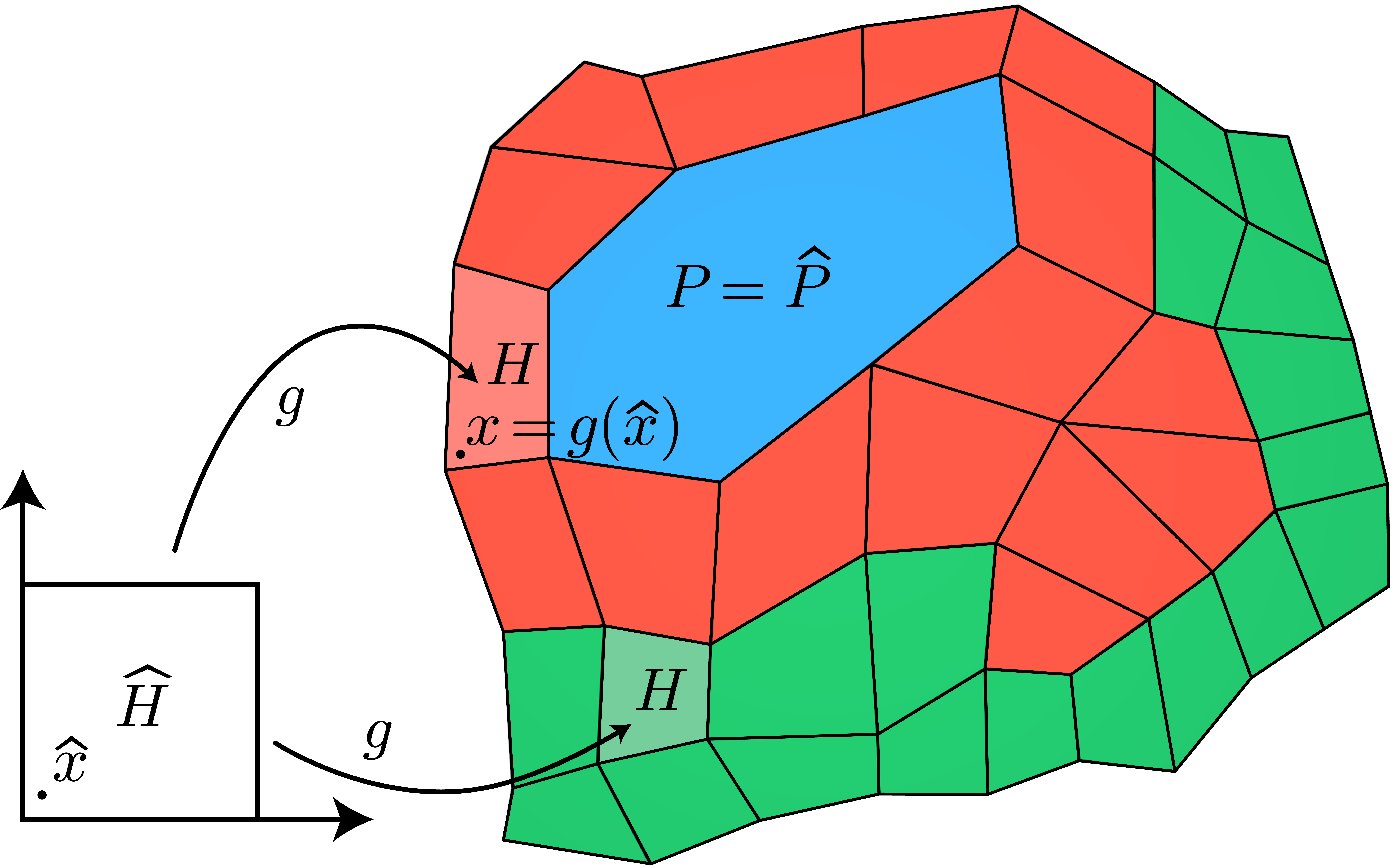}
\caption{Illustration of the geometric mapping.}
\label{fig:geom_mapping}
\end{figure}

The input complex, as it is typical for mesh representations, does not define a complete geometric realization of the complex: rather it only includes vertex positions and element connectivity.
We define a complete geometric realization as the \emph{geometric map} $\bg\colon \param{\cM} \rightarrow \threed$, from the parametric domain $\param{\cM}$ to the physical space. 
We use $\param{\bx}$ for points in the parametric domain, and $\bx$ for points in the physical space, and denote the image of the geometric map by $\Omega = \bg(\param{\cM})$ (\Cref{fig:geom_mapping}).

The definition requires bootstrapping: $\bg$ is first defined on $\param{\cH}$.
 For example, the simplest geometric map $\bg$ on $\param{\cM}$ can be obtained by trilinear interpolation:
$\bg$ restricted to a unit cube $\param{H} \subset \param{\cM}$ is a trilinear interpolation of the positions of the vertices of its associated hexahedron $H$.
We make the following assumption about $\bg(\param{\cH})$: the map is bijective on the faces  of $\cH$, corresponding to the
boundary of any polyhedral cell $P$, and the union of the images of these 
faces does not self-intersect and encloses a volume $P'$. Section~\ref{sec:refinement} explains how this is ensured.
Then we complete $\param{\cM}$ by adding the volume $P'$ as the 
parametric domain for $P$.  We add this volume to the parametric domain $\param{\cM}$, identifying corresponding faces with faces in $\param{\cH}$, and defining 
the geometric map to be the identity on these domains. 

The simplest trilinear map is adequate for elements of the mesh outside the regular part $\cS$, but is insufficient for accuracy on the regular part, as discussed below.
We consider a more complex definition of $\bg$ ensuring $C^1$ smoothness across interior edges and faces of $\cS$, described in Section~\ref{sec:geommap}, after we describe our basis construction.
Our construction is \emph{isoparametric}, that is, it uses the same basis for the geometric map as for the solution.

In other words, on $\cQ$ we use the standard tri-qua\-dra\-tic geometric map that maps each reference cube $[0,1]^3$ to the actual hex-element in the mesh. On $\cS$ we use a $C^1$ spline mapping, explained in \Cref{sec:geommap}. On the polyhedral part, the geometric map is the identity, thus all quantities are defined directly on the physical domain.

\paragraph{Overview of the basis and discretization construction.}

Given an input complex $\cM$, we construct a set of bases $\param{\basis}_i \colon \param{\cM} \to \reals$, $i=1,\ldots,N$, such that:

\begin{itemize}[leftmargin=*]
\item the restriction of basis function $\param{\basis}_i$ to spline compatible hexahedral domains $\param{S} \in \param{\cS}$ is a spline basis function;
\item the restriction to hexahedra $\param{Q} \in \param{\cQ}$ is a standard
  triquadratic ($Q_2$) element function;
\item the restriction to polyhedra $\param{P} \in \param{\cP}$ (or $P \in \cP$) is a harmonic-based
  nonconformal, third-order accurate basis function.
\end{itemize}
\noindent
The degrees of freedom (dofs) corresponding to basis functions $\param{\basis}_i$ are associated with:
\begin{itemize}[leftmargin=*]
\item each hexahedron either in $\cS$ or adjacent to a spline-compatible one (\emph{spline cell dofs});
\item each boundary vertex, edge, or face of $\cS$ (\emph{spline boundary dofs}); these are needed to have correct approximation on the boundary;
\item each vertex, edge, face, and cell of $\cQ$ (\emph{triquadratic element  degrees of freedom}).
\end{itemize}
The total number of degrees of freedom is denoted by $N$.
While most of the construction is independent of the choice of PDE
(we assume it to be second-order), with a notable exception of the 
consistency condition for polyhedral elements, we use the Poisson equation 
to be more specific. 

Note that hexahedra adjacent to $\cS$, but not in $\cS$ (i.e., hexahedra in $\cQ$)
get both spline dofs and triquadratic element dofs: such a cell may have $\geq 28$ dofs instead of 27.

Polyhedral cells are not assigned separate degrees of freedom: the basis functions
with support overlapping polyhedra are those associated with dofs at incident hexahedra.

We assemble the standard stiffness matrix for an elliptic PDE, element-by-element, performing integration on the hexahedra $\param{H}$ of $\param{\cM}$ and polyhedra $P$.
The entry $K_{ij}$ of the stiffness matrix $\mat{K}$ for the Poisson equation is computed as follows:

\begin{equation}
 K_{ij} = \sum_{\param{C} \in \param{\cM}} \int_{\bg(\param{C})}
 \nabla \basis_i(\bx) \cdot \nabla\basis_j(\bx) \,\dd\bx,
\label{eq:stiffness}
\end{equation}
where $\basis_i  = \param{\basis}_i \circ \bg^{-1}$.
The actual integration is performed on the elements in the parametric domain
$\param{\cM}$, using a change of variables $\bx = \bg(\param{\bx})$ for every element:
\begin{equation}
K_{ij} = \sum_{\param{C} \in \param{\cM}} \int_{\param{C}} \nabla \param{\basis}_i(\param{\bx})\T \, \mat{A}(\param{\bx}) \, \nabla \param{\basis}_j(\param{\bx}) \, \abs{\dD\bg}\, \dd\param{\bx}
\label{eq:stiffness-param}
\end{equation}
where $\mat{A}(\param{\bx})$ is the metric tensor of the geometric map $\bg$ at $\param{\bx}$, given by $\dD\bg^{-1} \dD\bg^{\mathsf{-T}}$, with $\dD\bg$ being the Jacobian of $\bg$.

In the next sections, we describe the construction of the basis on each element type, the geometric map, and the stiffness matrix construction.


\newcommand{\uh}{u_h}
\newcommand{\vh}{v_h}

\section{Basis construction}
\label{sec:basis}

We seek to construct a basis on $\Omega = \bg(\param{\cM})$ that has the following properties:
\begin{enumerate}
\item it is $C^0$ everywhere on $\Omega$, $C^1$ at regular edges and
  vertices, and $C^\infty$ within each $H$ and $P$ (polynomials on hexahedra).
\item it has approximation order 3 on each $H$ and $P$.
\end{enumerate}

 \begin{figure}
 \centering
 \includegraphics[width=0.4\linewidth]{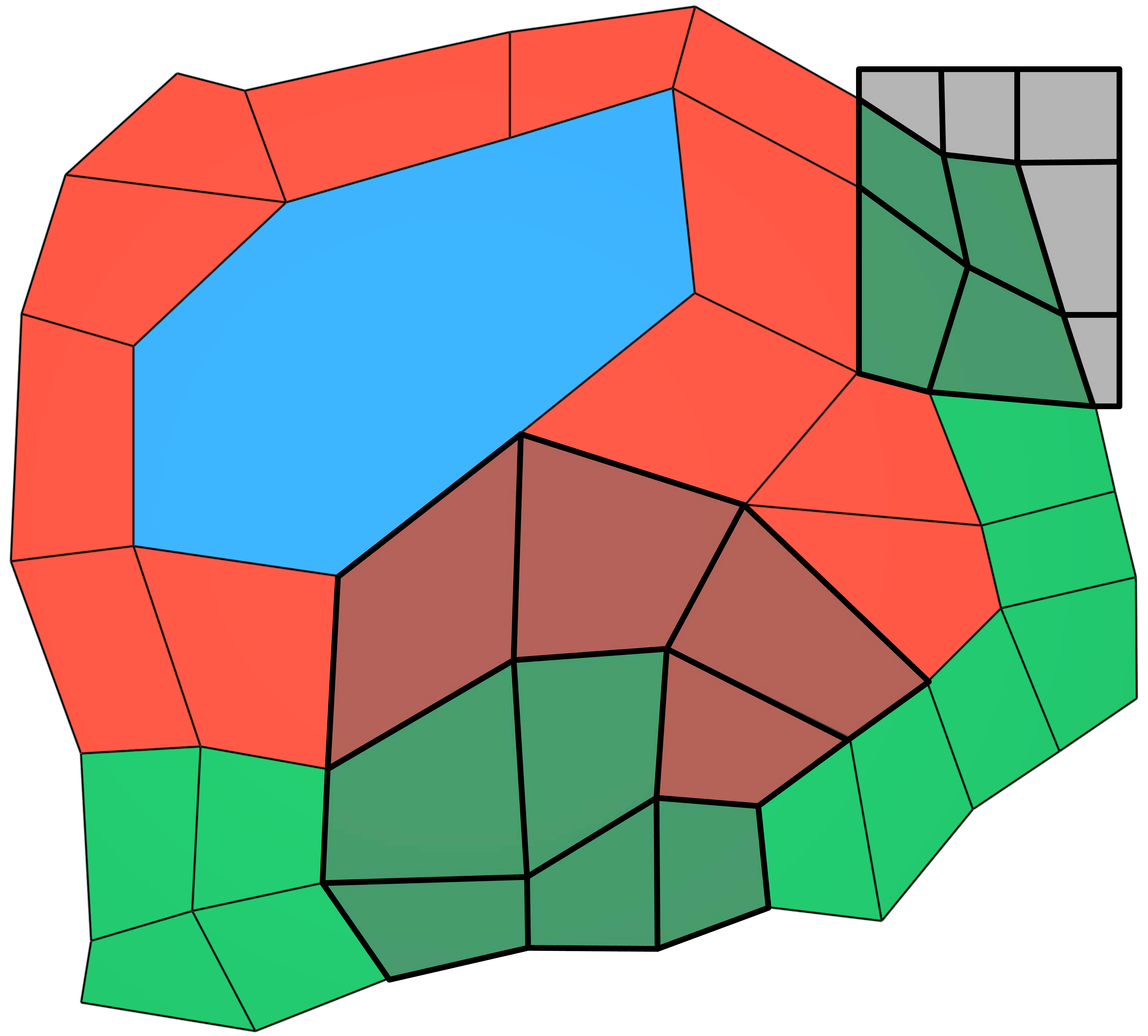}
 \caption{Spline local grid (shown in dark), for an internal and a boundary quadrilateral. The color codes are as defined in \Cref{fig:complex}.}
\label{fig:hex-local-grid}
\end{figure}

 \begin{figure}
 \centering
 \hfill
  \includegraphics[width=.4\linewidth]{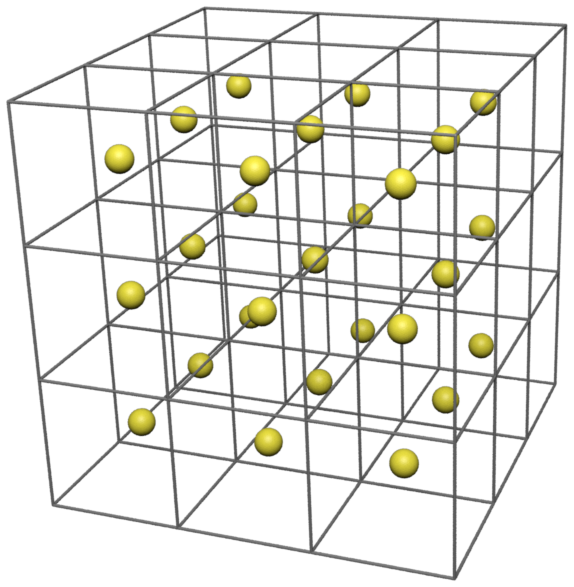}\hfill
  \includegraphics[width=.4\linewidth]{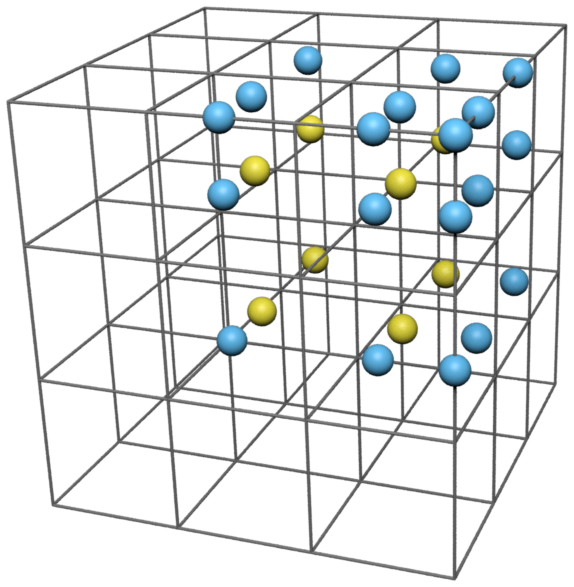}\hfill{}
 \caption{Spline hex degrees of freedom for a central element and a corner one.}
\label{fig:spline-dofs}
\end{figure}

The unknown function $u$ on the domain $\Omega$ is approximated by $\uh = \sum_{i=1}^N u_i \basis_i$, where $\basis_i$ are the basis functions.
The support of each basis function is a union of a set of the images under $\bg$ of cells in $\param{\cM}$.

The actual representation of the basis, which allows us to perform per-element construction of the stiffness matrix, consists of three parts.
The first two parts are local: we define a \emph{local} set of dofs and a \emph{local} basis.
For  hexahedral elements, there are several types of local polynomial bases, each coming with its set of local dofs, associated with a local \emph{control mesh} for the element. These basis functions are encoded as sets of polynomial coefficients.
For polyhedral elements, all local basis functions are weighted combinations of harmonic kernel functions and a triquadratic polynomial, so these are encoded as kernel centers, weights and polynomial coefficients.

The third part is the \emph{local-to-global} linear map that represents local dofs in terms of the global ones. Importantly, unlike most standard FEM
formulations, our local-to-global maps are not necessarily simply identifying local dofs the global ones: some local dofs are linear combinations of global ones.
These maps are formally represented by $m \times N$ matrices, where $m$ is a small number of local dofs, and $N$ is the total number of global dofs.
However, as the elements local dofs depend only on nearby global dofs, these
matrices have a small number of nonzeros and can be encoded in a compact form.

In the following, we consider the construction of these three elements (set of local basis
functions, set of local dofs, local-to-global map) for each of our three element
types. But before we can construct the basis for each element, hexahedral elements need to be classified into $\cS$ (spline-compatible) and $\cQ$.

\subsection{Spline-compatible hexahedral elements}
\label{sec:spline-compatible}

We define a hexahedron $H$ to be spline-compatible, if its one-ring cell neighborhood is a $3\times 3\times 3$ regular grid, possibly cut on one or more sides if $H$ is on the boundary, see \Cref{fig:hex-local-grid}.

\emph{The local dofs} of this element type form a $3\times 3 \times 3$
grid (for interior elements), with the element in the center (Figure~\ref{fig:spline-dofs} left);
for boundary elements, there are still 27 dofs, ensuring a full triquadratic polynomial reproduction.
If a single layer with 9 dofs is missing, we add an extra
degree of freedom for each face of the local $3 \times 3 \times 2$ grid
corresponding to the boundary. Other cases are handled in a similar manner; e.g. the configuration for a regular corner is shown in Figure~\ref{fig:spline-dofs}, right.

\begin{figure}
 \centering
 \hfill
  \includegraphics[width=.4\linewidth]{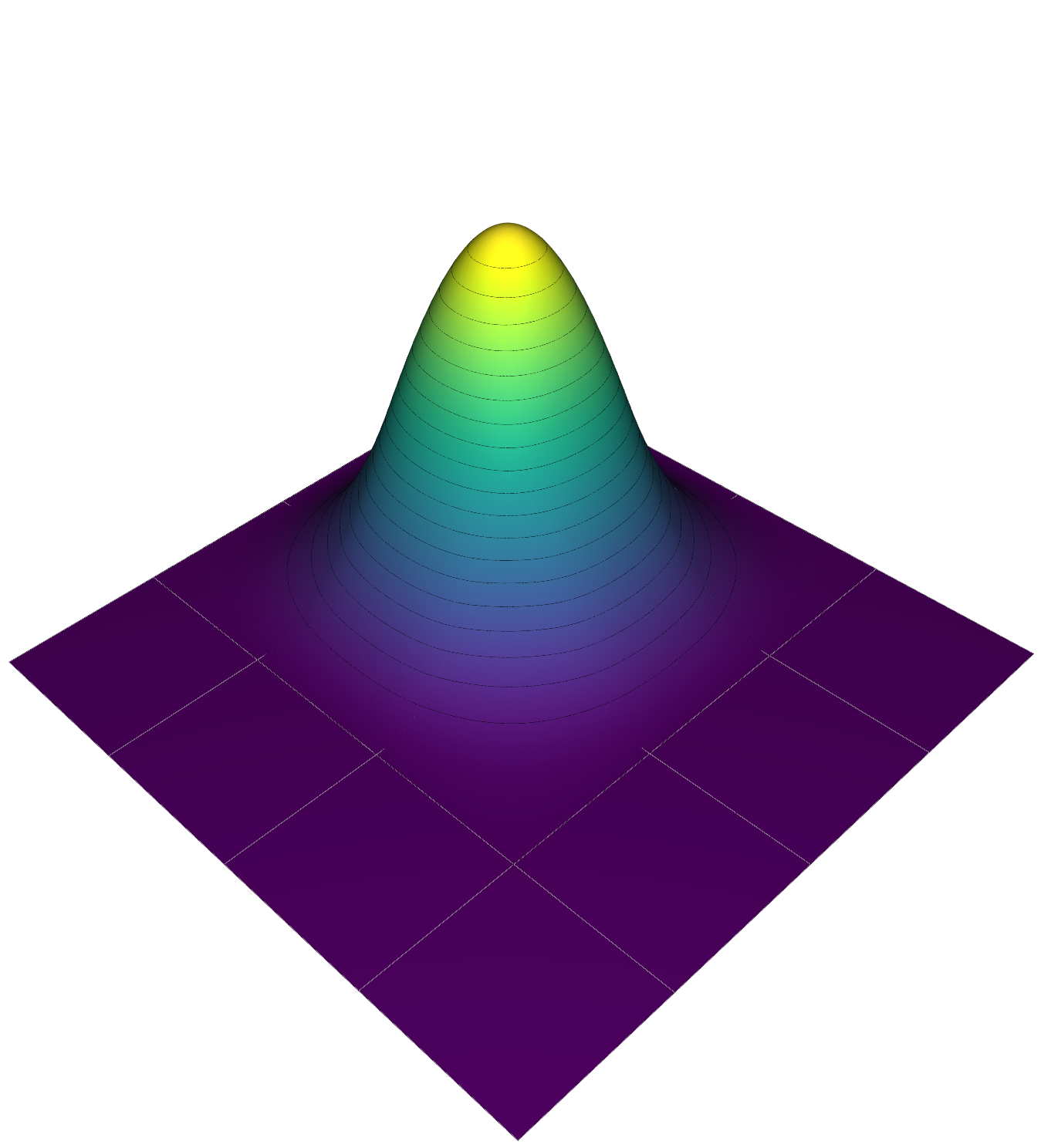}\hfill
  \includegraphics[width=.4\linewidth]{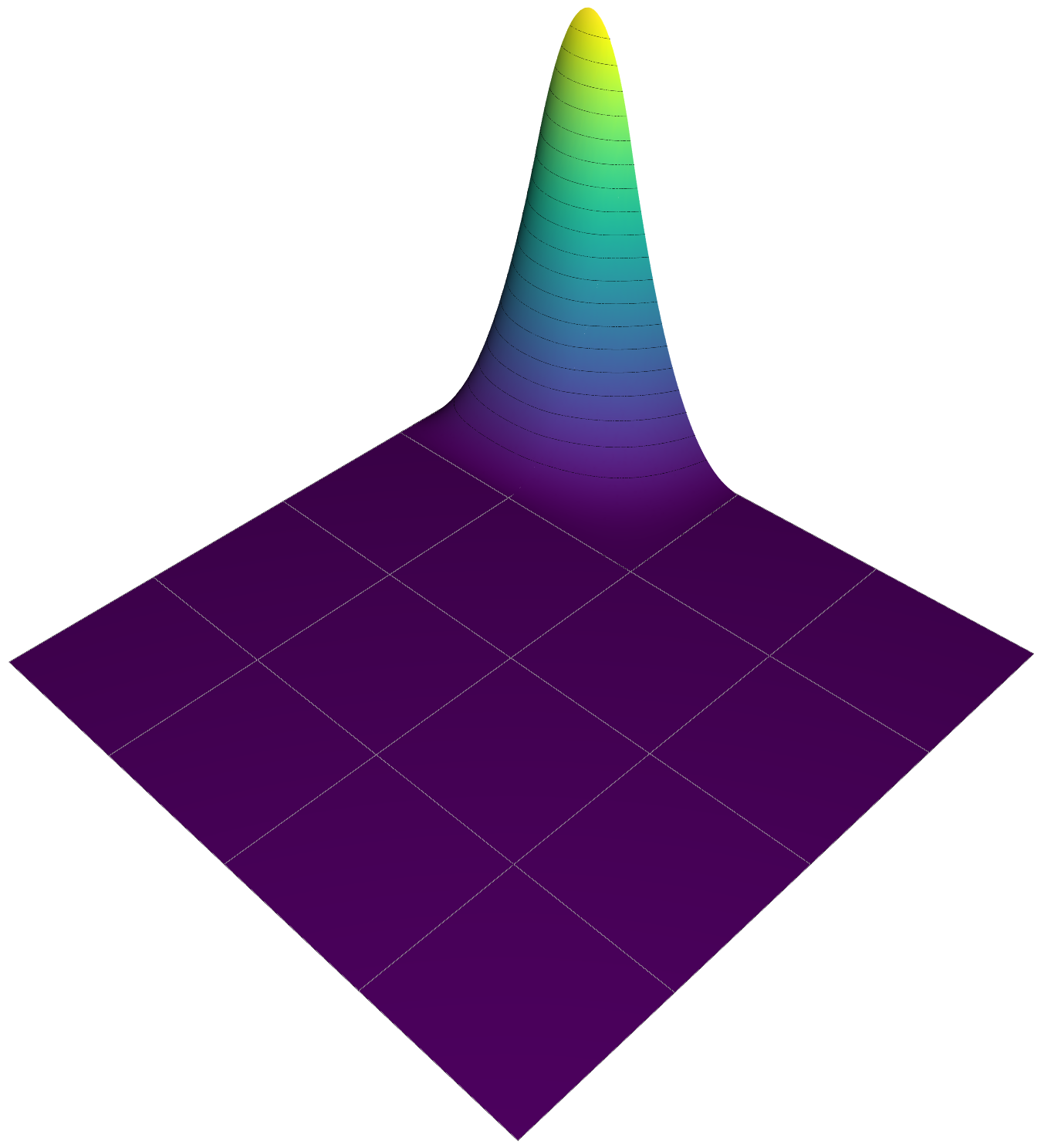}\hfill{}
 \caption{Plot of the spline bases for a regular 2D grid.}
\label{fig:spline-bases}
\end{figure}

\emph{The basis functions} in this case are just the standard triquadratic
uniform spline basis functions for interior hexahedra.  For the boundary case,
we use the knot vector $[0,0,0,1,2,3]$ in the direction perpendicular to the
boundary. \Cref{fig:spline-bases} shows an example of the bases in 2D, for an internal node on the left and for a boundary node on the right.
Finally, the \emph{local-to-global} map simply identifies local basis dofs with corresponding global ones.

Compared to a standard $Q_2$ element, the ratio of degrees of freedom to the number of elements is much lower (a single degree of freedom per element for splines), although the approximation order is the same.

\subsection{$Q_2$ hexahedral elements}
\label{seq:q2}

This element is used for all remaining hexahedra. It is a standard element, widely used in finite element codes. \emph{Local dofs} for this element are associated with the element vertices, edge midpoints, face centers, and cell centers (Figure~\ref{fig:basis-nodes}).

\begin{figure}
   \centering
     \includegraphics[width=0.35\linewidth]{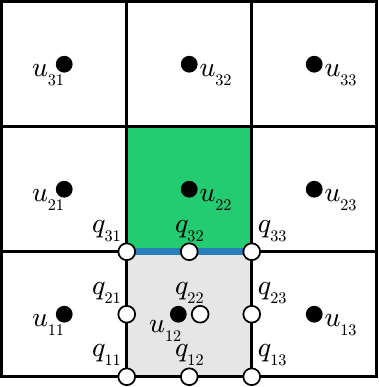}
  \caption{Local-to-global map for a $Q_2$ element (gray) adjacent to a single spline element (green).}
\label{fig:Q2-dofs}
\end{figure}

The \emph{local basis functions} for the element are obtained as the tensor product of the interpolating quadratic bases on the interval $[0,1]$, consisting of $(t-\frac{1}{2})(t-1)$, $(t-\frac{1}{2})t$, and $(t-1)t$ (\Cref{app:pk}).
The only complicated part in the case of $Q_2$ elements is the definition of the
local-to-global map. For the two-dimensional setting, it is illustrated in Figure~\ref{fig:Q2-dofs}.

The difficulty in the construction of this map is due to the interface between spline elements and $Q_2$ elements, and the need to ensure continuity between the two.
In the two-dimensional case, suppose that a $Q_2$ element $Q\in\cQ$ shares an edge with exactly one quad spline element $S\in\cS$.
Let $u_{ij}$, $i,j= 1,\ldots, 3$, be the global dofs of the spline element,
and let $q_{ij}$, $i,j=1,\ldots, 3$, be the degrees of freedom of the $Q_2$ element, as shown in the picture.

In this case, we ensure $C^0$ continuity of the basis by expressing the
values of the polynomials on $Q \in \cQ$ at the shared boundary points in terms of global degrees of freedom.
Since both the $Q_2$ and the spline basis restricted to an edge are quadratic polynomials, they only need to be equal on three distinct points of the edge to ensure continuity.
By noticing that the $Q_2$ basis is interpolatory at the nodes, it is enough to evaluate the spline basis at these edge nodes.

For the two-dimensional example in \Cref{fig:Q2-dofs}, the local-to-global map for the local dofs $q_{31}$, $q_{32}$ and $q_{33}$ along the edge (in blue) that the $Q_2$ element shares with the spline is obtained as follows:
\begin{equation}
  \begin{split}
    q_{31} & = \frac{1}{4} \left(u_{11} + u_{12} + u_{21} + u_{22} \right),\\
    q_{32} &=  \frac{3}{8} \left(u_{12} + u_{22}\right) + \frac{1}{16}\left(u_{11}+u_{21}+u_{13}+u_{23} \right),\\
    q_{33} & = \frac{1}{4} \left(u_{12} + u_{13} + u_{22} + u_{23} \right).\\
  \end{split}
  \label{eq:locglobq2}
\end{equation}

In 3D, the construction is similar. We first identify all spline bases overlapping with a local dof $q_{ij}$ on the boundary of a $Q_2$ element (i.e., a vertex, edge, or face dof). To determine the weights of the local-to-global map, we evaluate each spline basis on the local dof $q_{ij}$ and set it as weight.

\begin{figure}
 \centering
 \hfill
 \includegraphics[width=.2\linewidth]{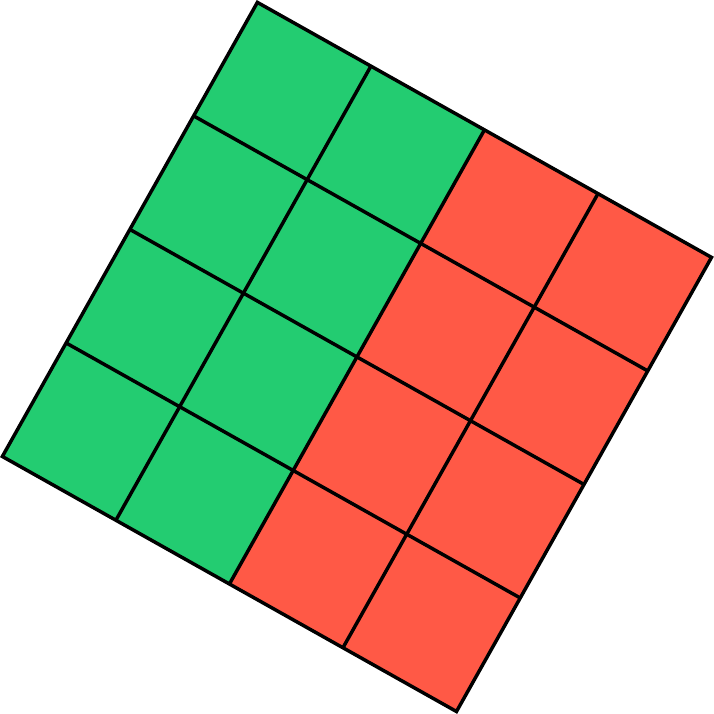}\hfill
  \includegraphics[width=.37\linewidth]{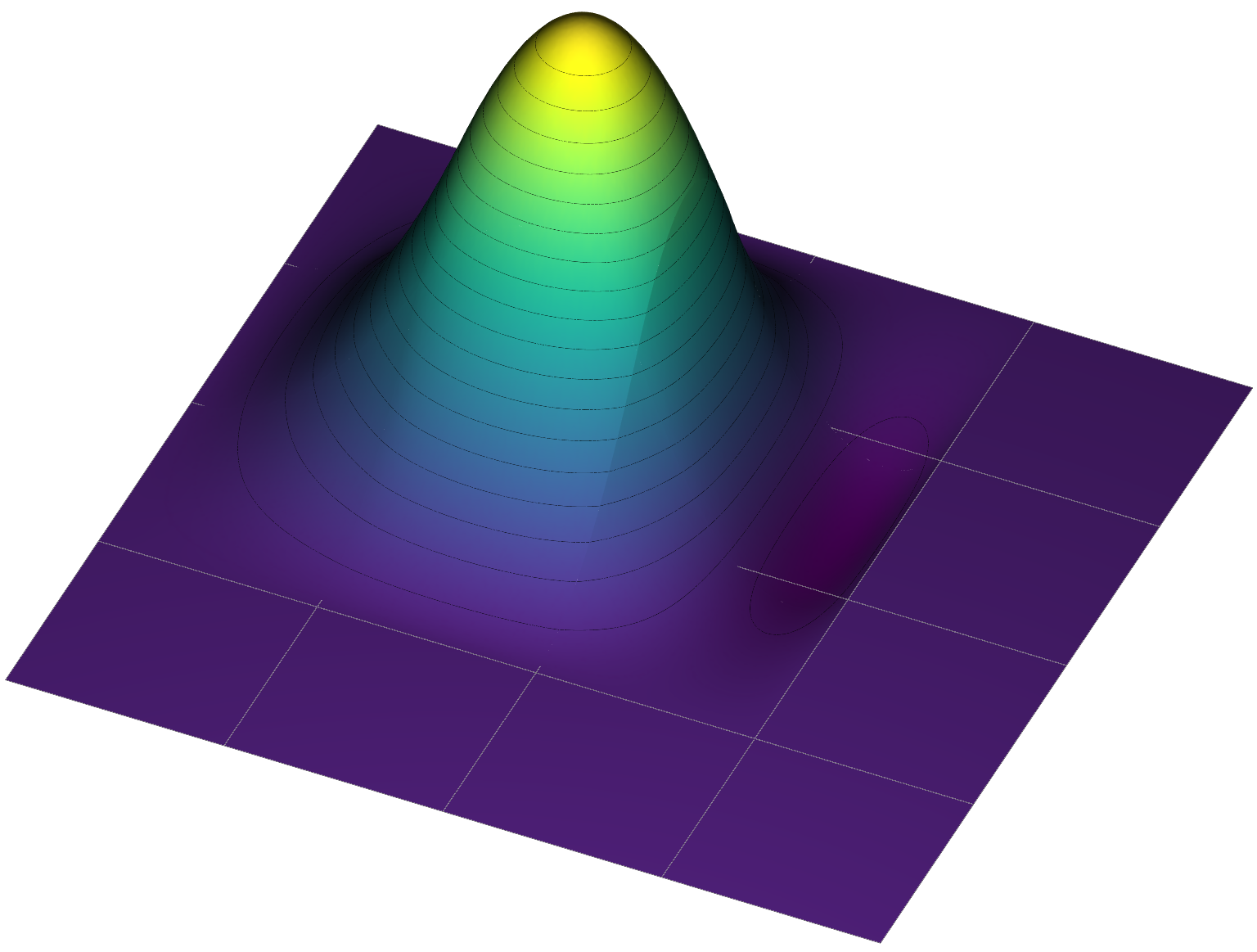}\hfill
  \includegraphics[width=.37\linewidth]{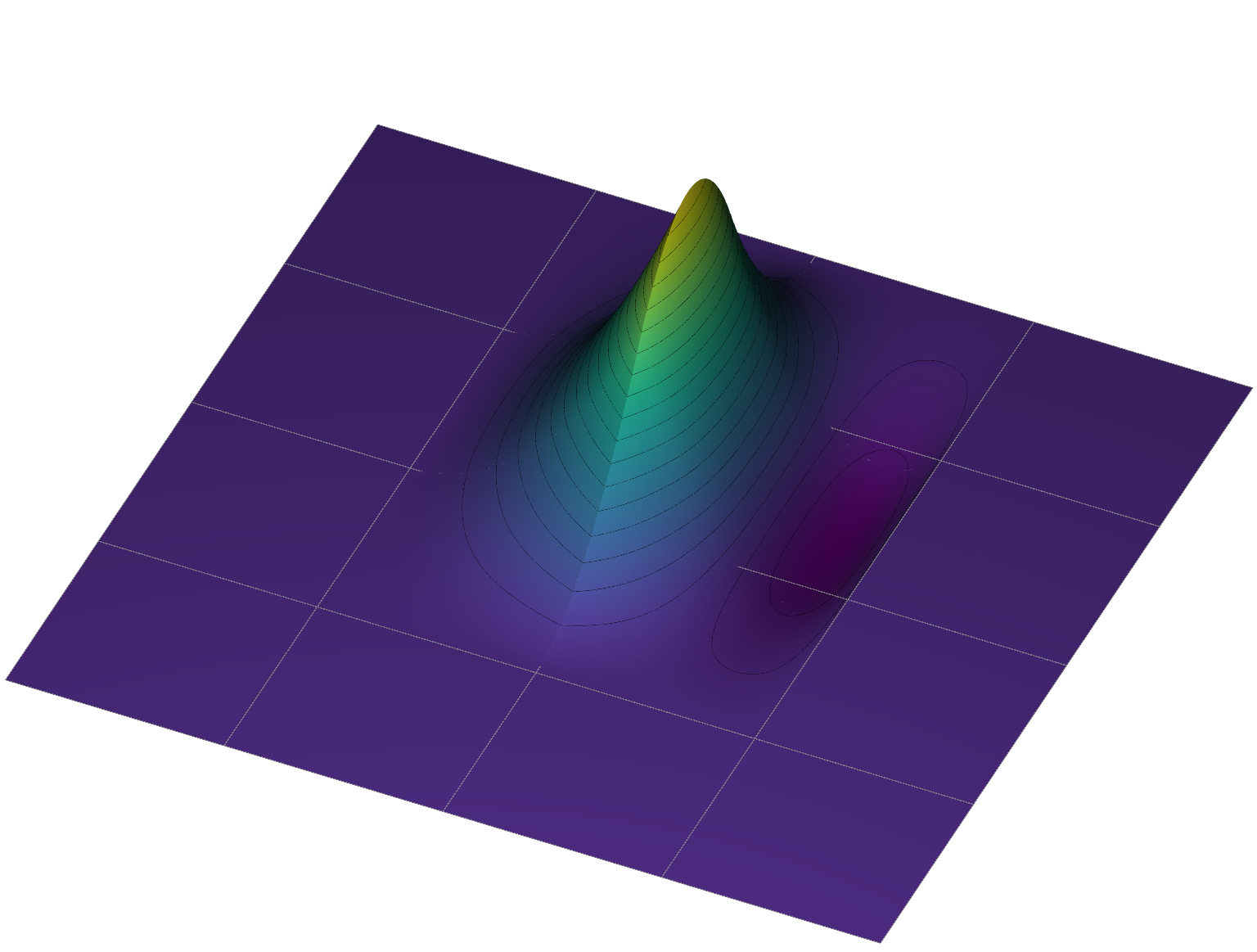}
 \caption{Plot of the bases on a junction between a regular (green) and an irregular (red) part for a regular 2D grid.}
\label{fig:spline-q2-bases}
\end{figure}

The remaining degrees of freedom of the $Q_2$ element are identified with global $Q_2$ degrees of freedom at the same locations. We note once again that at the center of cells in $\cQ$ with neighboring cells in $\cS$, there are \emph{two} dofs, one spline dof and one $Q_2$ dof. \Cref{fig:spline-q2-bases} shows an example of two basis functions on the transition from the regular part on the left to the ``irregular'' part on the right. We clearly see that on the regular part the bases are splines and on the irregular one are the standard $Q_2$ basis function: on the interface the functions are only $C_0$.

\subsection{Basis construction on polyhedral cells}
\label{sec:polyhedral}

The construction of the basis on the polyhedral cells is quite different from the construction of the basis on hexahedra. For hexahedra, the basis functions are defined on the parametric domain $\param{\cM}$, and are remapped to $\Omega\subset \threed$ via the geometric map.
For polyhedra, we construct the basis directly in physical space.

On possible option to construct the basis on polyhedral cells is to split each polyhedral cells into tetrahedra. This approach has two main disadvantages: (i) it requires the use of pyramids to ensure conformity to the neighboring hexahedra, (ii) it is difficult to guarantee a sufficient element quality after subdivision. Instead, we follow the general approach of \cite{Martin:PFE:2008} with two important alterations designed to ensure third-order convergence.

Recall that all polyhedron faces are quadrilateral, and all
polyhedra are surrounded by hexahedra, specifically $Q_2$ hexahedra as their neighborhood is not regular. Moreover, since we always perform an initial refinement step, there are no polyhedral cells touching each other.
We use the degrees of freedom on the faces of these elements as degrees of freedom for the polyhedra, therefore the \emph{local-to-global} map in this case is trivial.

Each dof is already associated with a basis function $\basis_j$ defined on the hexahedra adjacent to the polyhedron.
We construct the extension of $\basis_j$ to the polyhedron $P$ from $\numK$  harmonic kernels $\kernel_i(\bx) = \norm{\bx-\kc_i}^{-1}$ centered at positions $\kc_i$ outside the polyhedron and quadratic monomials $\monom_d(\bx)$, $d=1, \ldots, 10$, as
\begin{align}
    \restr{\basis_j}{P}(\bx) 
    &= \sum_{i=1}^{\numK} w^j_i \kernel_i(\bx) +
  	   \sum_{d=1}^{10} a^j_d \monom_d(\bx) \label{eq:basisdef} \\
    &= \vect{w}^j \cdot \vect{\kernel}(\bx) + 
       \vect{a}^j \cdot \vect{\monom}(\bx), \nonumber 
\end{align}
where $\vect{w}^j = (w^j_1, \dots, w^j_{\numK})\T$, $\vect{\kernel} = (\kernel_1, \dots, \kernel_{\numK})\T$, $\vect{a}^j = (a^j_1, \dots, a^j_{10})\T$, and $\vect{q} = (\monom_1, \dots, \monom_{10})\T$. The coefficients $\vect{w}^j$ and $\vect{a}^j$ are $r \times k$ and $r \times 10$ matrices, respectively, with $r=1$ (scalar PDEs) or $r=2,3$ (vector PDEs).
Following \shortcite{Martin:PFE:2008}, the weights $w^j_i, a^j_d \in \reals^k$ are determined
using a least squares fit to match the values of the basis $\basis_j$ evaluated on a set of points sampled on the boundary of the polyhedron $P$.

In \cite{martin2011flexible}, it is shown that this construction automatically
guarantees reproduction of linear polynomials if $\monom_d$ are linear; the quadratic case is fully analogous.
However, this condition is insufficient for high-order convergence, because our basis is \emph{non-conforming}, that is non $C^0$. In the context of the second-order PDEs we are considering, it means that it lacks $C^0$ continuity on the boundary of the polyhedron.
For this type of elements, additional \emph{consistency conditions} are required to ensure high-order convergence. These conditions depend on the PDE that we need to solve.

\begin{figure}
  \centering
  \hfill
  \includegraphics[width=0.35\linewidth]{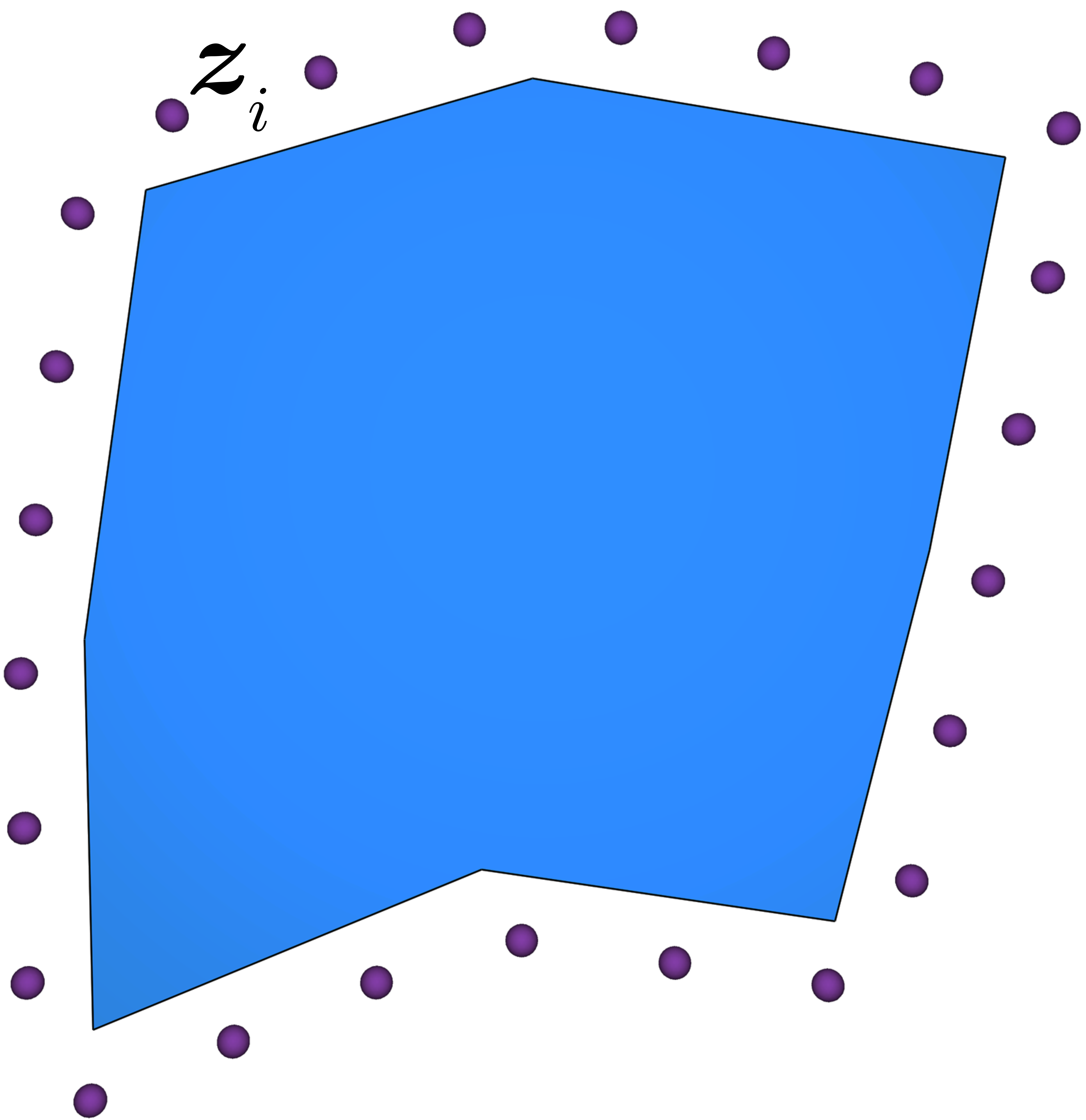}\hfill
  \includegraphics[width=0.45\linewidth]{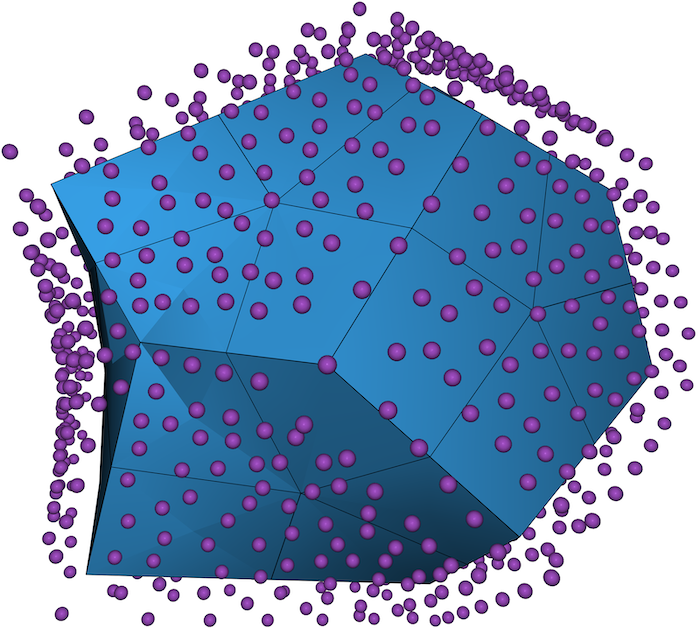}\hfill
  \caption{The local basis for a polygon consists of the set of triquadratic polynomials $\monom_d$
   and harmonic kernels $\kernel_i$ centered at shown locations $\kc_i$.}
\label{fig:poly-dofs}
\end{figure}

\paragraph{FEM theory detour.} \label{sec:poly-contraints} To achieve higher order convergence \emph{three} conditions need to be satisfied: (1) polynomial reproduction; (2) consistency, which we discuss in more detail below and (3) quadrature accuracy.
We refer to standard FEM texts such as \cite{braess2007finite}
for details, as well as to virtual element method literature (e.g., \cite{de2016nonconforming} is closely related).

To satisfy the third condition, we use high-order quadrature on the polyhedron: we decompose it into tetrahedra and use Gaussian quadrature points in each tetrahedron (the decomposition is detailed in \Cref{sec:starshaped}).
The first condition, polynomial reproduction, is ensured by construction of the basis above.

The second constraint, consistency, requires further elaboration. We first derive it for the Poisson 
equation, and then summarize the general form. 
We leave as future work the complete proof of the convergence properties of our method (cf. \cite{de2016nonconforming}),
which requires, in particular, a proper stability analysis. 
Nevertheless, in Section~\ref{sec:evaluation} we provide numerical evidence that our method does converge at the expected rate, and that its conditioning is not affected in a significant way by the presence of nonconforming polyhedral elements.

The standard way to find the solution of a PDE for a finite element system
is to consider its weak form. For the Poisson equation,  find $u$ such that
\begin{equation}
	\int_\Omega  \Delta u \, v  = -
	\int_\Omega \nabla u \cdot \nabla v = \int_\Omega f v,\qquad
	\forall v
\end{equation}
 \emph{Remark.} We omit, for readability, the integration 
variable $\dd\bx$. In the remaining formulas we use integration over 
the physical space exclusively, in practice carried over to the parametric 
space by adding the Jacobian of the geometric map. 

Then, $u$ is approximated by $\uh=\sum_i u_i \basis_i$, and $v$ is taken to be in the space spanned by the basis functions $\basis_j$.
The stiffness matrix entries are obtained as $K_{ij} = \sum_{C} \int_{C} \nabla \basis_i \cdot\nabla\basis_j$, where the integral is computed per element $C$, leading to the discrete system $\mat{K} \bu = \mathbf{f}$ (\Cref{app:fem}).

For general non-conforming elements, however, we cannot rely on this standard approach.
For example, if we consider piecewise-constant elements for the Poisson equation, the stiffness matrix would be all zeros.

However, \emph{for a given PDE}, one can construct converging non-conforming elements. One condition that is typically used, is that the discrete matrix, constructed per element as above, gives us \emph{exact} values of the weak-form integral for all polynomials reproduced by the basis (cf. $k$-consistency property in \cite{de2016nonconforming}).

As our basis reproduces triquadratic monomials (i.e., they are in the span of bases $\basis_i$), we have $\monom_d(\bx)=\sum_i \monom_d^i\basis_i(\bx)$. To ensure consistency, we require that any \emph{nonconforming basis function $\basis_j$} satisfies
\begin{equation}
-\int_{\bg(\param{\cM})} \Delta \monom_d \basis_j = \sum_i K_{ij} \monom_d^i
\label{eq:consistency}
\end{equation}
for all triquadratic monomials $\monom_d$.

To convert this equation to an equation for the unknown coefficients $w^j_i$ and
$a^j_d$, we observe that
\begin{equation}
	\sum_i K_{ij} \monom_d^i =
	\int_{\bg(\param{\cM})}\bigg(\sum_i \monom_d^i \nabla \basis_i\bigg) \cdot \nabla \phi_j =
	\int_{\bg(\param{\cM})} \nabla \monom_d \cdot \nabla\basis_j
\end{equation}
due to the polynomial reproduction property. Separating the integral into the part over the hexahedra $g(\param\cM\setminus P)$ and over the polyhedron $P=g(P)$, we write
\begin{align}
  \sum_i K_{ij} \monom_d^i &=
  C_H + \int_P \nabla \monom_d \cdot \nabla \Big(\vect{w}^j \cdot \vect{\kernel} + \vect{a}^j \cdot \vect{\monom}\Big) \\
  & = C_H + \vect{b}\T \vect{w}^j + \vect{c}\T \vect{a}^j \nonumber
\end{align}
where 
\begin{equation*}
C_H = \sum_{\param{C}\in\param{\cM}\setminus P} \int_{\bg(\param{C})} \nabla \monom_d \cdot \nabla \basis_j,\quad
\vect{b} =  \int_P \nabla \monom_d \cdot \nabla \vect{\kernel}, \quad
\vect{c} =  \int_P \nabla \monom_d \cdot \nabla \vect{q}.
\end{equation*}

Similarly, the left-hand side of \Cref{eq:consistency} is reduced to a
linear combination of $\vect{w}^j$ and $\vect{a}^j$.
This forms a set of additional constraints for the coefficients of the basis functions on the polyhedron.
To enforce them on each polyhedron, we solve a constrained least squares system for each nonconforming basis function and store the obtained coefficients.

Importantly, the addition of constraints to the least squares system does not violate 
the polynomial reproduction property on the polyhedron. This can be seen as follows.
Let $\vh$ be the linear combination of basis functions $\phi_i$ overlapping $P$ that yields a triquadratic mononomial $q_d$ when restricted to  $P$. Then $\vh$ is continuous on $\Omega$:  the samples at the points of the boundary are from a quadratic function, therefore, match 
exactly the quadratic continuation to adjacent hexahedra. 

The consistency condition (Equation~\ref{eq:consistency}) applied to $\vh$ simply states that it satisfies the integration by parts formula, which it does as it is $C^0$ at the element boundaries, and smooth on the elements:
\[ 
-\sum_{C} \int_C \Delta \monom_d  \vh = \sum_C \int_C \nabla \monom_d \cdot \nabla \vh .
\]
We conclude that $\vh$ is in the space defined by the consistency constraint, and imposing this constraint preserves polynomial  reproduction. See \Cref{app:constraints} for the complete list of constraints for the Poisson equation.

More generally, for a linear PDE and for any polynomial $\monom$ (for vector  PDEs, e.g.,  elasticity,  this means that all components are polynomial) we require
     \[
     a(\monom_d, \vh) = a_h(\monom_d, \vh),
      \quad\text{where}\quad
     a(u, v) = \int_\Omega \mathcal{F}(x, u, \nabla u, \Delta u) v
     \]
     where $\mathcal{F}$ is a linear function of its arguments depending on $u$, and $a_h$ is defined as a sum of integral over $\Omega$ after formal integration by parts of $\mathcal{F}$, to eliminate the second-order derivatives.  For a conforming $C^0$ basis, this condition automatically follows from the integration by parts formulas, which are applicable. 
     We now split the two bilinear forms as $a=a^H+a^P$ and $a_h=a_h^H + a_h^P$ where $a^H$ and $a_h^H$ contains the integral over the hexahedral known part, and $a^P$ and $a_h^P$ the integral over the polyhedral unknown part. Thus, for a basis $\basis_j$ we obtain the following set of constraints
     \begin{align*}
     a^H(\monom_d,\basis_j) - a^H_h(\monom_d,\basis_j) 
     &= a^H\left(
         \monom_d, \vect{w}^j \cdot \vect{\kernel}(\bx) + \vect{a}^j \cdot \vect{\monom}(\bx)
         \right) \\
     & -a^H_h\left(
         \monom_d, \vect{w}^j \cdot \vect{\kernel}(\bx) + \vect{a}^j \cdot \vect{\monom}(\bx)
     \right).
     \end{align*}

For a scalar-valued PDE, we have the same number of constraints (5 in 2D and 9 in 3D) as monomials $\monom_d$, thus we are guaranteed to have a solution that respects the constraints for any $k > 0$.
For vector PDEs (e.g., elasticity), we impose the additional constraints such that the 
coefficients $\left( w^j_{i} \right)_\alpha$ are the same for all dimensions $\alpha = 1,\ldots,r$, $r=2$ or $r=3$, which simplifies the implementation, but increases the number of required centers $\vect{z}_j$, 
    so that all constraints can be satisfied. More explicitly, for vector PDEs we require that the 
    constraints  
    \[
    a(q^s_d \vect{e}_\alpha, \phi^s_j \vect{e}_\beta) =   a_h(q^s_d \vect{e}_\alpha, \phi^s_j \vect{e}_\beta)
    \]
    for $\alpha,\beta = 1,\ldots,r$ are satisfied, with $q^s$ and $\phi_j^s$ denoting scalar polynomials 
    and scalar basis functions respectively, defined as in \eqref{eq:basisdef} for dimension 1, and $\vect{e}_\alpha$ is the unit vector for axis $\alpha$.
     For dimensions 2 and 3, the number of monomials $q$ is $5$ and $9$ respectively.  The number of constraints is given by $r^2q - q$,  and thus we will need at least 15 $\vect{z}_i$ in 2D and 72 in 3D to ensure that the constraints are respected.

\subsection{Imposing boundary conditions}

We consider two standard types of boundary conditions: Dirichlet (fixed function values on the boundary) and Neumann (fixed normal derivatives at the boundary). Neumann (also known as natural) boundary conditions are handled in the context of the variational formulation of the problem as extra integral terms, in the case of inhomogeneous conditions.  Homogeneous conditions do not require any special treatment and are imposed automatically in the weak formulation.

We assume that the Dirichlet conditions are given as a continuous function defined on the boundary of the domain. For all boundary dofs, we sample the boundary condition on the faces of the domain and perform a least-squares fit to retrieve the nodal values.


\section{Geometric map construction}
\label{sec:geommap}

The geometric map is a map from $\param{\cM}$ to $\Omega\subset\threed$, defined per element.
Its primary purpose is to allow us to construct basis functions $\param{\basis}_i$ on reference domains (i.e., the elements of $\param{\cM}$ that are unit cubes), and then to remap them
to the physical space as $\basis_i = \param{\basis}_i \circ \bg^{-1}$.
As the local basis on the polyhedral elements is constructed directly in the physical space, $\bg$ is the identity on these elements. 

The requirements for the geometric map are distinct for the spline and $Q_2$ elements, and are matched by using spline basis itself for $\cS$ and trilinear interpolation for $Q_2$ elements.

Because of the geometric mapping $\bg$, for the quadratic spline, the basis $\basis_i$ does not reproduce polynomials in the physical space; nevertheless, the  approximation properties of the basis are retained \cite{bazilevs2006isogeometric}.

For $Q_2$ elements, Arnold et al.~\shortcite{arnold2002approximation} shows that \emph{bilinear} maps are sufficient, and in fact allow to retain reproduction of triquadratic polynomials in the physical space.
This is very important for the basis construction on polyhedral elements, as polynomial reproduction on these elements depends on reproduction of polynomials on the polyhedron boundary.

\paragraph{Computing the geometry map.}
If we assume that the input only has vertex positions $\vect{v}_i$ for $\cM$, we solve the equations  $\bg(\param{\bx}_i) = \vect{v}_i$, which is a linear system of equations in terms of coefficients of $\bg$ in the basis we choose.
In the trilinear basis, the system is trivial, as the coefficients coincide with the values at $\bx_i$, and these are simply set to $\vect{v}_i$.
For the triquadratic basis, this is not the case, and a linear system needs to be solved. If the system is under-determined, we find the least-norm solution.

\section{Mesh preprocessing and refinement}
\label{sec:refinement}

Without loss of generality, we restrict the meshing discussion to 2D, as the algorithm introduced in this section extends naturally to 3D.

For the sake of simplicity, in this discussion the term polygon refers to non-quadrilateral elements.
As previously mentioned, our method can be applied to hybrid meshes without two adjacent polygons and without polygons touching the boundary, which we ensure with one step of refinement. While our construction could be extended to support these configurations, we favored refinement due to its simplicity. Refining polygonal meshes is an interesting problem on its own: while there is a canonical way to refine quads, there are multiple ways to refine a polygon. We propose the use of polar refinement~(\Cref{sec:polar}), which has the added benefit of allowing us to resample large polygons to obtain a uniform element size. However, to avoid self-intersections between edges during the refinement, we impose each polygon be star-shaped. This condition is often, but not always, satisfied by existing hybrid meshers: we thus introduce a simple merging and splitting procedure to convert hybrid meshes into star-shaped polyhedral meshes (\Cref{sec:starshaped}), and then detail our refinement strategy (\Cref{sec:polar}).

Another advantage of restricting ourselves to star-shaped polygons is that partitioning it into triangles (respectively tetrahedra in 3D) is trivial, by introducing a point in the kernel and connecting it to all the boundary faces. This step is required to generate quadrature points for the numerical integration (\Cref{sec:polyhedral}): the quality of the partitioning is usually low, but this is irrelevant for this purpose.

\subsection{Mesh preprocessing}

\begin{figure}
  \centering
  \includegraphics[width=0.7\linewidth]{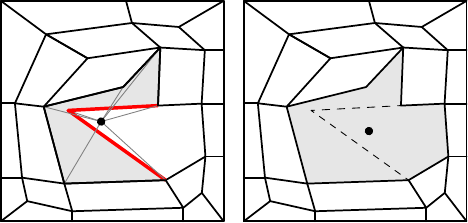}
  \caption{Our algorithm iteratively merges polygons (gray polygon in the first image), until the barycenter of the merged polygon is inside its kernel (gray polygon in the second).}
\label{fig:merging}
\end{figure}

\label{sec:starshaped}
We propose a simple and effective algorithm to convert polygonal meshes into star-shaped polygonal meshes, by combining existing polygons until they are star-shaped (and eventually splitting them if they contain a concave part of the boundary).

For every non-star-shaped polygon, we compute its barycenter and connect it to all its vertices (Figure \ref{fig:merging}, left). This procedure generates a set of intersecting segments (red in Figure~\ref{fig:merging}), which we use to grow the polygon by merging it with the faces incident to each intersecting segment. The procedure is repeated until no more intersections are found, which usually happens in one or two iterations in our experiments. If we reach a concave boundary during the growing procedure, it might be impossible to obtain a star-shaped polyhedron by merging alone: In these cases, we triangulate the polygon, and merge the resulting triangles in star-shaped polygons if possible.

\subsection{Polar refinement}
\label{sec:polar}

\begin{figure}
  \centering
	\includegraphics[width=0.24\linewidth]{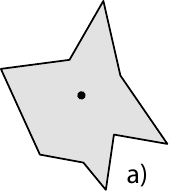}\hfill
	\includegraphics[width=0.24\linewidth]{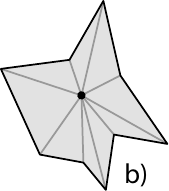}\hfill
	\includegraphics[width=0.24\linewidth]{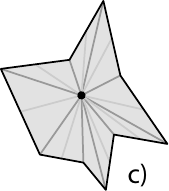}\hfill
	\includegraphics[width=0.24\linewidth]{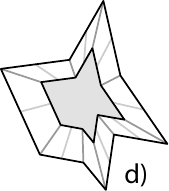}
  \caption{Polar refinement for polygons.}
\label{fig:refinement_procedure}
\end{figure}

\begin{figure}
  \centering
	\includegraphics[width=0.3\linewidth]{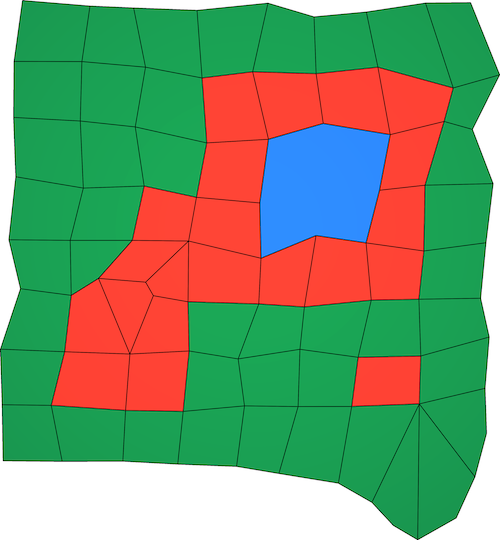}\hfill
	\includegraphics[width=0.3\linewidth]{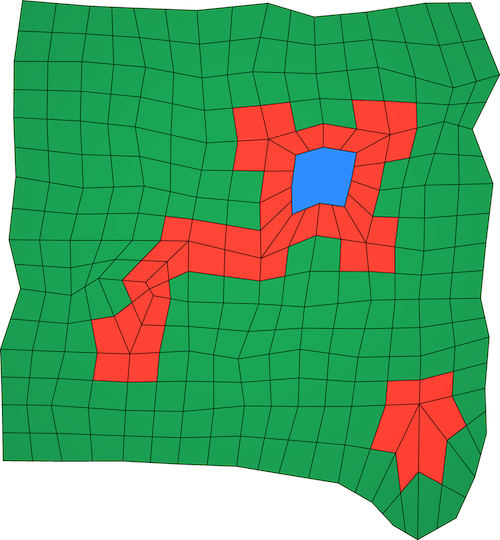}\hfill
	\includegraphics[width=0.3\linewidth]{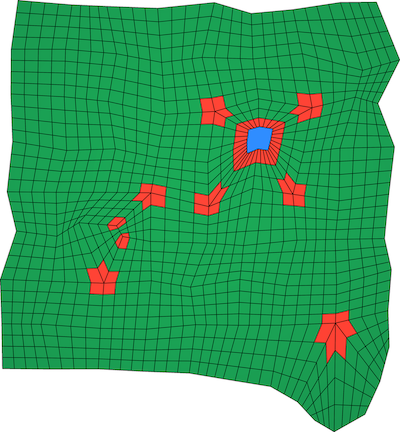}\\[1ex]
	\includegraphics[width=0.3\linewidth]{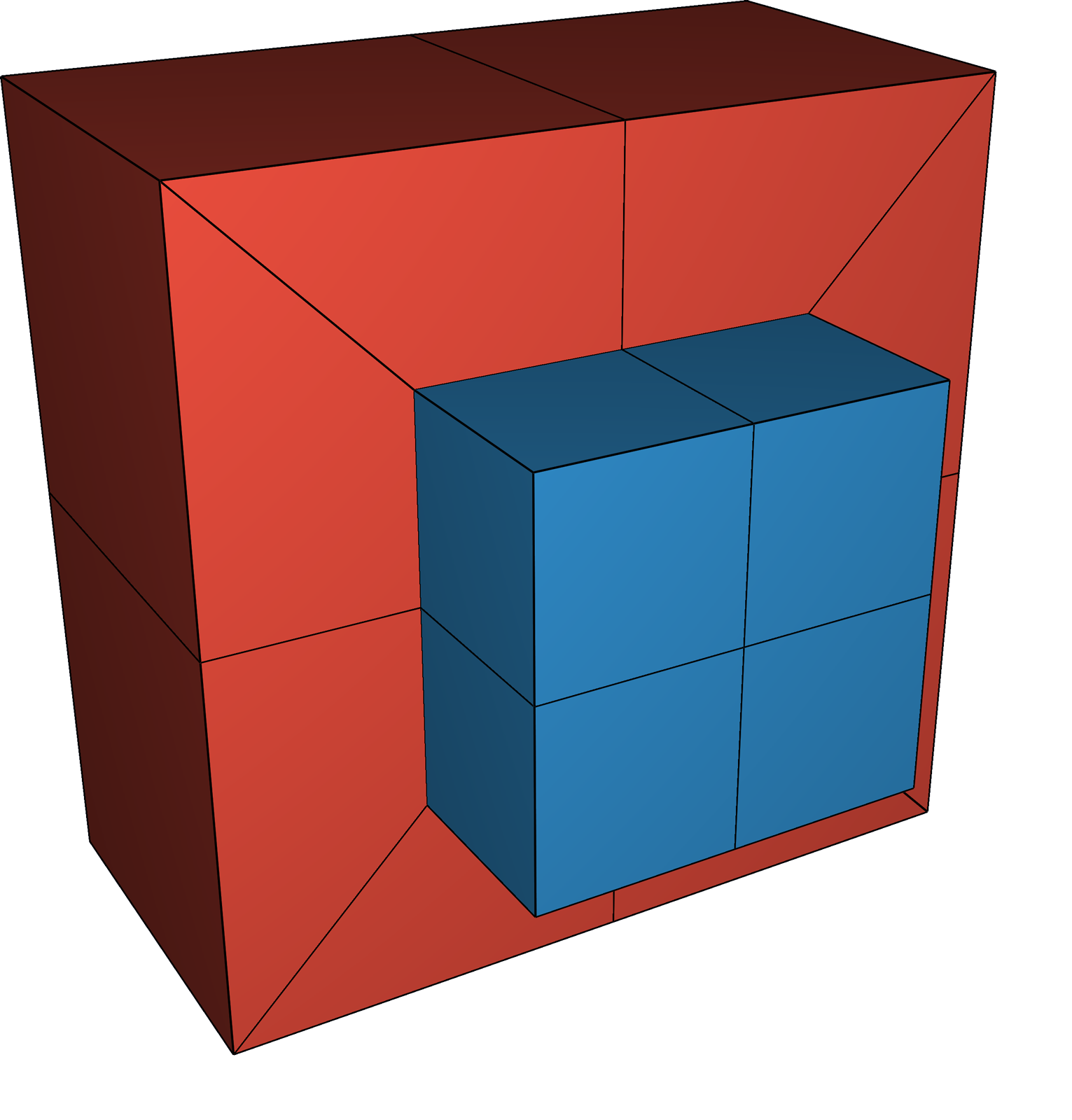}\hfill
	\includegraphics[width=0.3\linewidth]{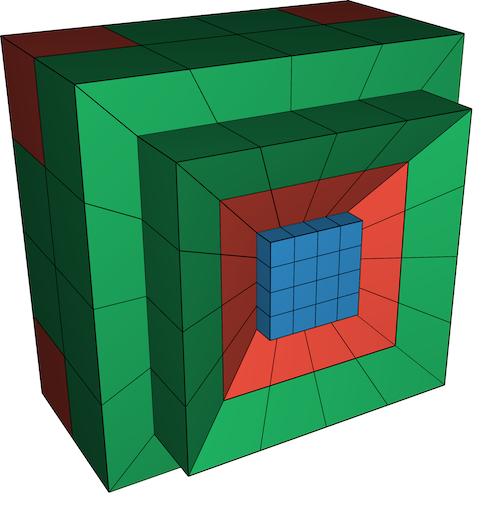}\hfill
	\includegraphics[width=0.3\linewidth]{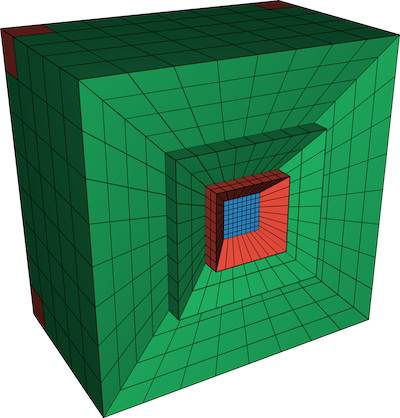}
  \caption{Example of polar refinement for a polygon and a polyhedron. The bottom view is a cut-through of the actual 3D mesh.}
\label{fig:refinement_example}
\end{figure}

Each star-shaped polygon is refined by finding a point in its kernel 
(Figure \ref{fig:refinement_procedure}, a), connecting it to all its vertices (b), splitting each edge with mid-point subdivision and connecting them to the point in the kernel (c), and finally adding rings of quadrilaterals around the boundary (d). Figure~\ref{fig:refinement_example} shows an example of polar refinement in two and three dimensions. The more splits are performed in the edge, the more elements are added. This is a useful feature to homogenize the element size in case the polygons were expanded too much during the mesh preprocessing stage. In our implementation, we split the edges evenly, ensuring that the shortest segment has a length as close as possible to the average edge length of the input mesh.

\section{Evaluation}

\label{sec:evaluation}

\begin{figure*}
	\centering
	\parbox{0.38\linewidth}{
	    \includegraphics[width=0.45\linewidth]{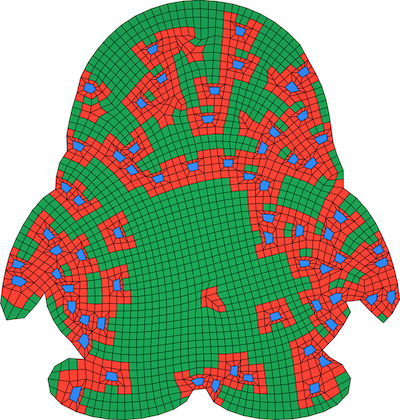}
	    \includegraphics[width=0.45\linewidth]{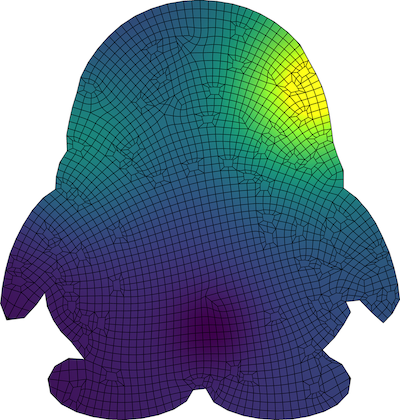}
	    }\hfill{}
	\parbox{0.2\linewidth}{
	    \includegraphics[width=0.9\linewidth]{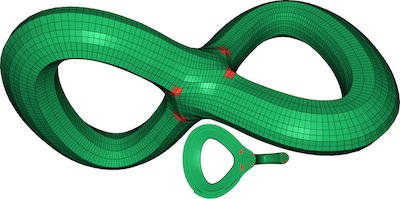}\par
	    \includegraphics[width=0.9\linewidth]{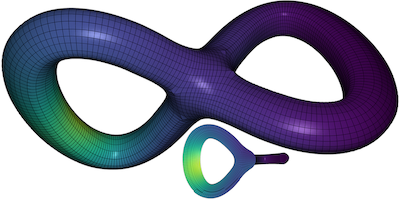}
	}\hfill{}
    \parbox{0.38\linewidth}{
	    \includegraphics[width=0.45\linewidth]{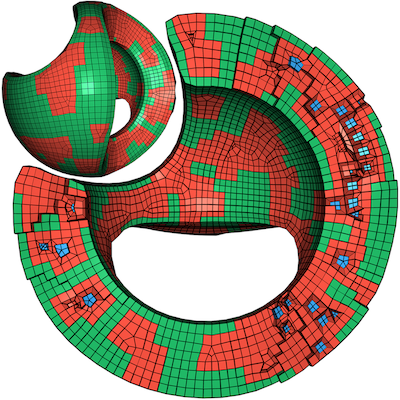}
	    \includegraphics[width=0.45\linewidth]{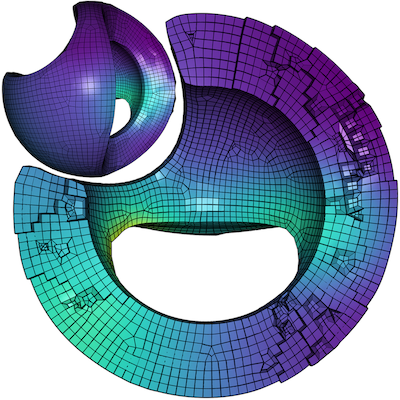}
	 }\\[1.5ex]
	 \parbox{0.38\linewidth}{\centering 2D polygon.}\hfill
	 \parbox{0.2\linewidth}{\centering 3D hexahedral mesh.}\hfill
	 \parbox{0.38\linewidth}{\centering 3D hybrid mesh.}
	\caption{Solution of the Poisson problem different meshes.}
	\label{fig:result_all}
\end{figure*}

\begin{figure*}
	\centering
	\makebox[\columnwidth][c]{
	\rotatebox{90}{\parbox{0.15\linewidth}{\centering\footnotesize Error}}\hfill
	\begin{overpic}[width=0.26\linewidth]
	    {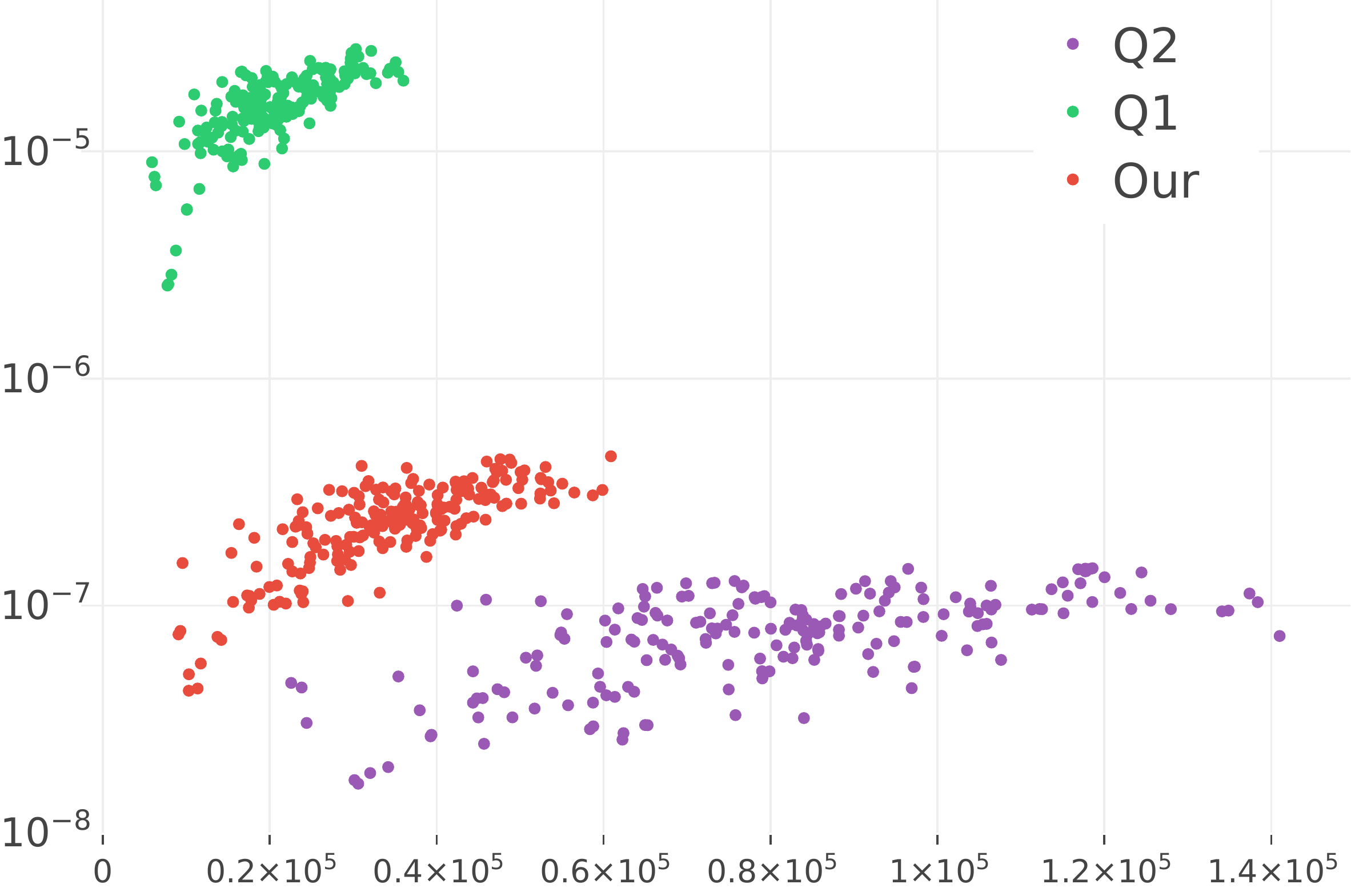}
		\put (50,60) {$L_2$}
	\end{overpic}\hfill{}
	\begin{overpic}[width=0.26\linewidth]
	    {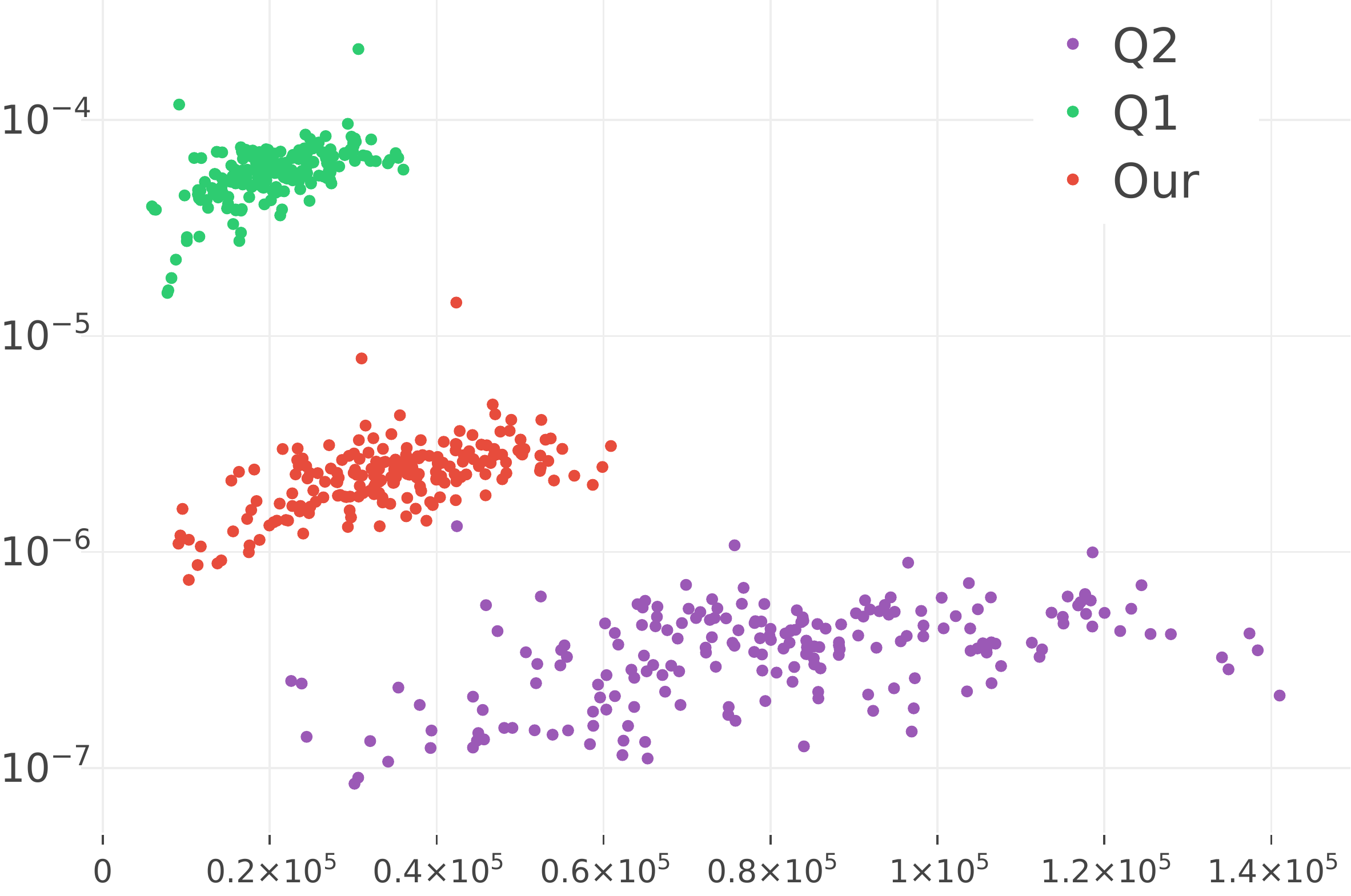}
		\put (50,60) {$L_\infty$}
	\end{overpic}\hfill{}
	\begin{overpic}[width=0.26\linewidth]{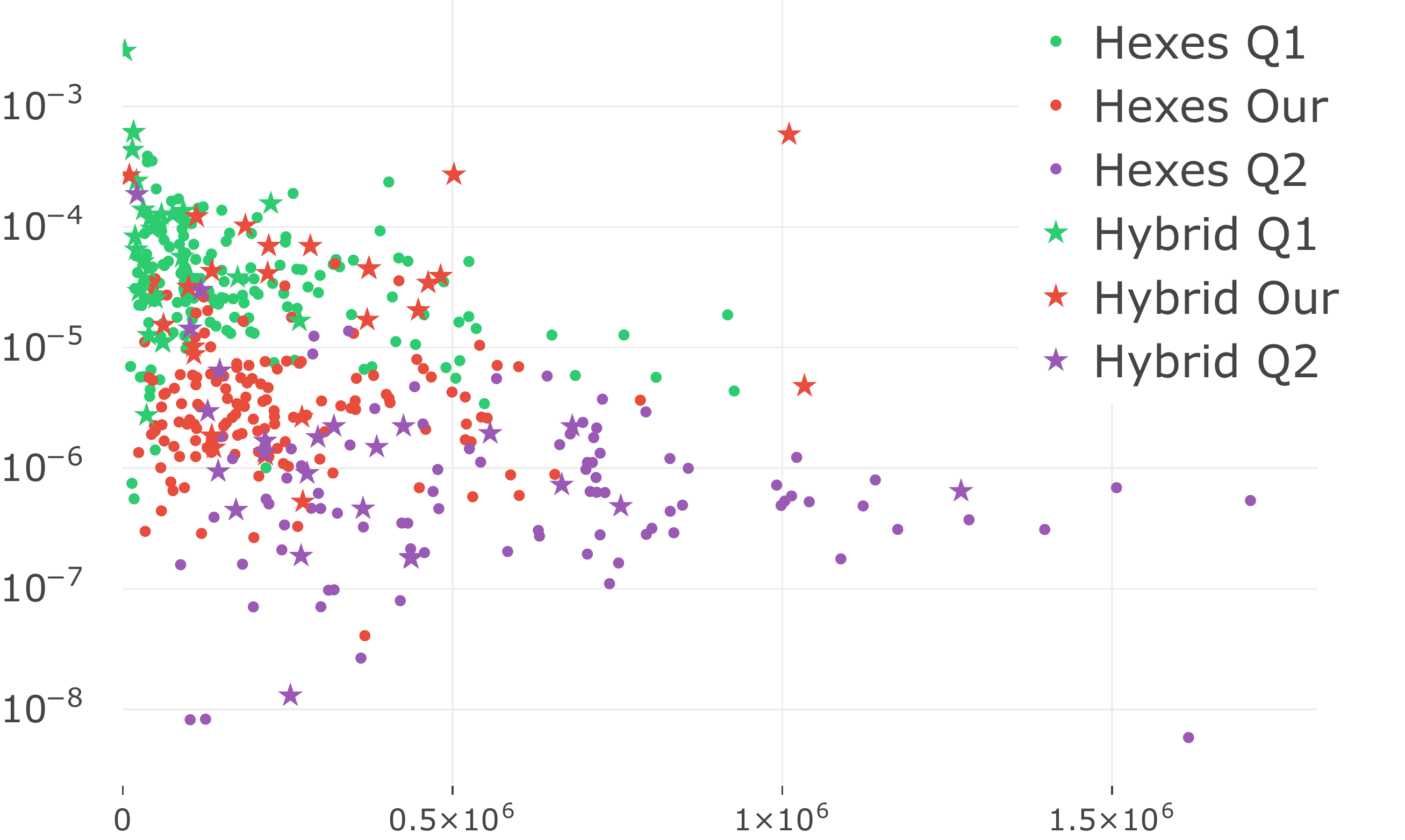}
		\put (50,60) {$L_2$}
	\end{overpic}\hfill{}
	\begin{overpic}[width=0.26\linewidth]{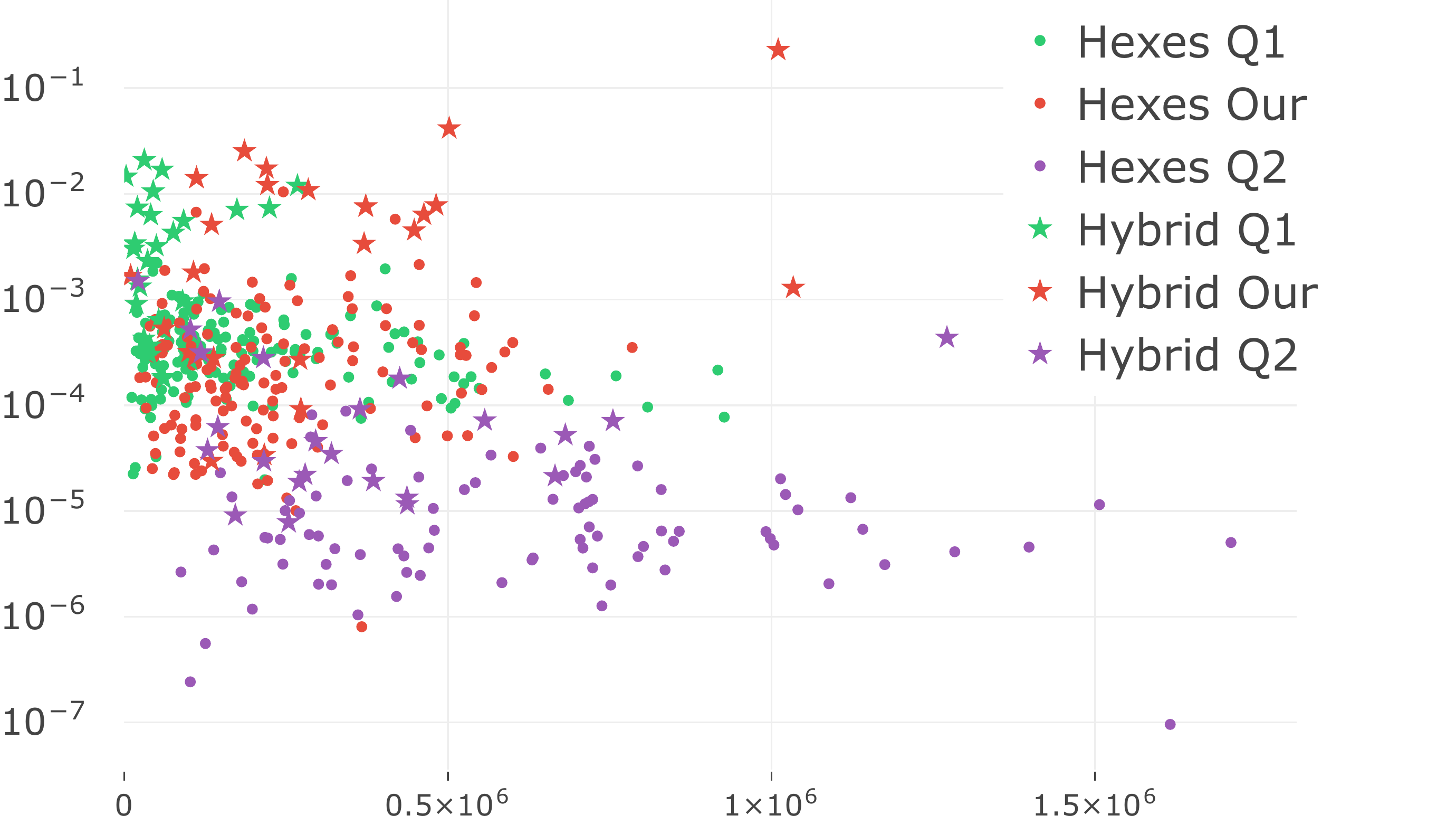}
		\put (50,60) {$L_\infty$}
	\end{overpic}\hfill{}
	}\\
	{\footnotesize Number of DOFs}
	\caption{Scatter plot of the $L_2$ and $L_\infty$ error versus the number of dofs on the 2D (first two) and 3D (last two) dataset.}
	\label{fig:scatter_plot}
\end{figure*}

We demonstrate the robustness of our method by solving the Poisson equation on a dataset of pure hex and hybrid meshes, consisting of 205 star-shaped polygonal meshes in 2D, 165 pure hexehedral meshes in 3D, and 29 star-shaped polyhedral meshes in 3D.
The dataset can be found at \url{https://cims.nyu.edu/gcl/papers/2019-Polyspline-Dataset.zip}.
All those meshes were automatically generated using \cite{Gao:2017,Gao:2017:RSS}. We show a selection of meshes from our dataset in Figures~\ref{fig:teaser} and \ref{fig:result_all}.

We evaluated the performance, memory consumption, and running time of our proposed spline construction compared with standard $Q_1$ and $Q_2$ elements. For our experiments, we compute the approximation error on a standard Franke's test function~\cite{Franke79} in 2D and 3D (\Cref{app:franke}).
Note that in all these experiments, we enforced the consistency constraints on the bases spanning the polyhedral elements, to ensure the proper convergence order.

The 2D experiments were run on a PC with an Intel\textregistered{} Core\texttrademark{}\textbf{} i7-5930K CPU @ 3.50GHz with 64 GB, while the 3D dataset was run on a HPC cluster with a memory limit of 64 GB.

\paragraph{Absolute Errors.} \Cref{fig:scatter_plot} shows a scatter plot of the $L_2$ and $L_\infty$ errors on both 2D and 3D datasets, with respect to the number of bases created by each type of elements ($Q_1$, $Q_2$, Splines), after one step of polar refinement.
The plot shows that in 2D both the $L_2$ and $L_\infty$ errors are about $1.5$ orders of magnitude lower for our splines compared to $Q_1$, while keeping a similar number of dofs. In comparison, $Q_2$ has lower error, but requires a much larger number of dofs. In 3D the spread of both errors is much larger, and the gain in $L_\infty$ is less visible, but still present, compared to $Q_1$.

\begin{figure}
	\centering
	\makebox[\columnwidth][c]{
	\includegraphics[width=0.5\columnwidth]{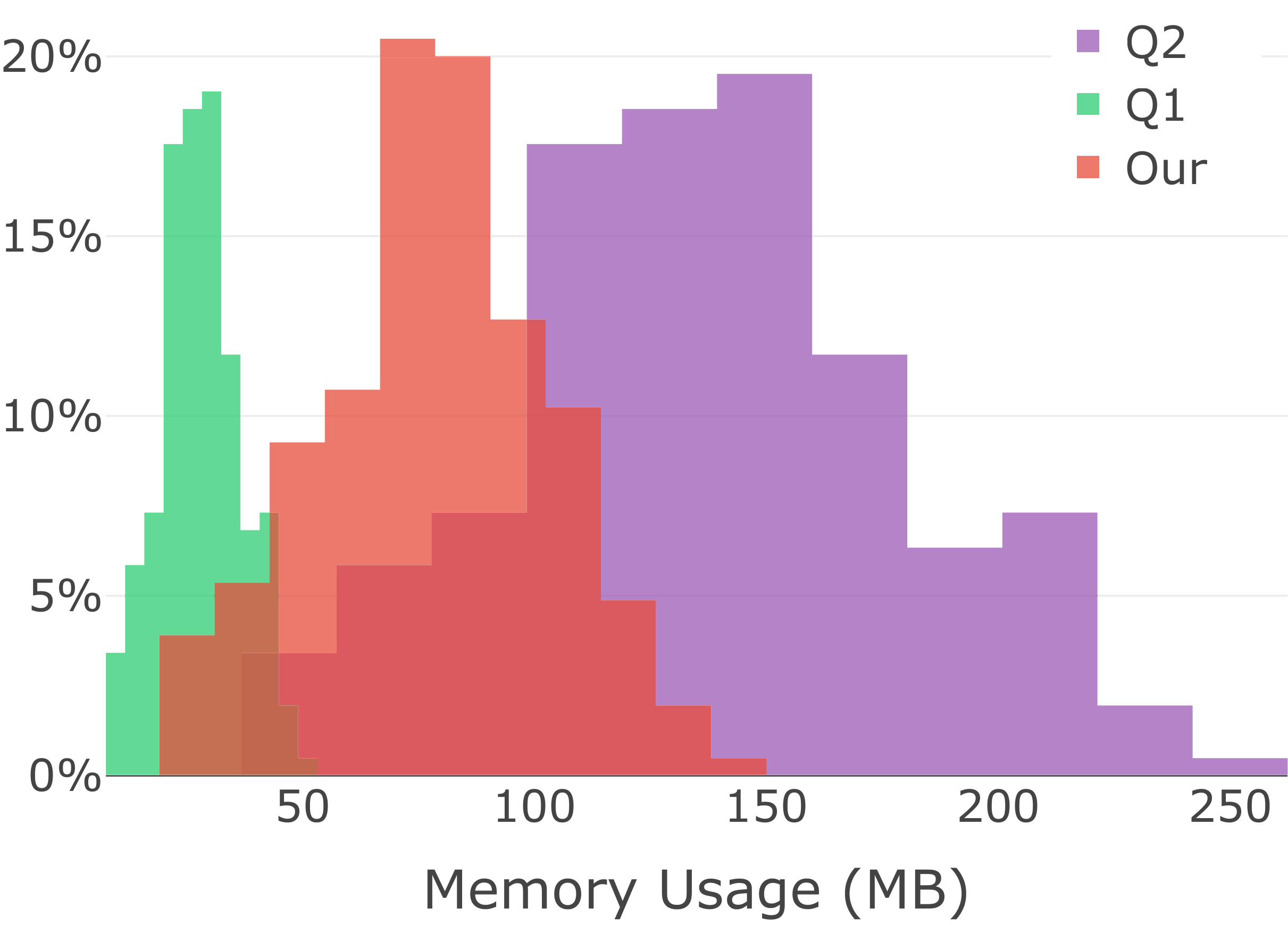}
	\includegraphics[width=0.5\columnwidth]{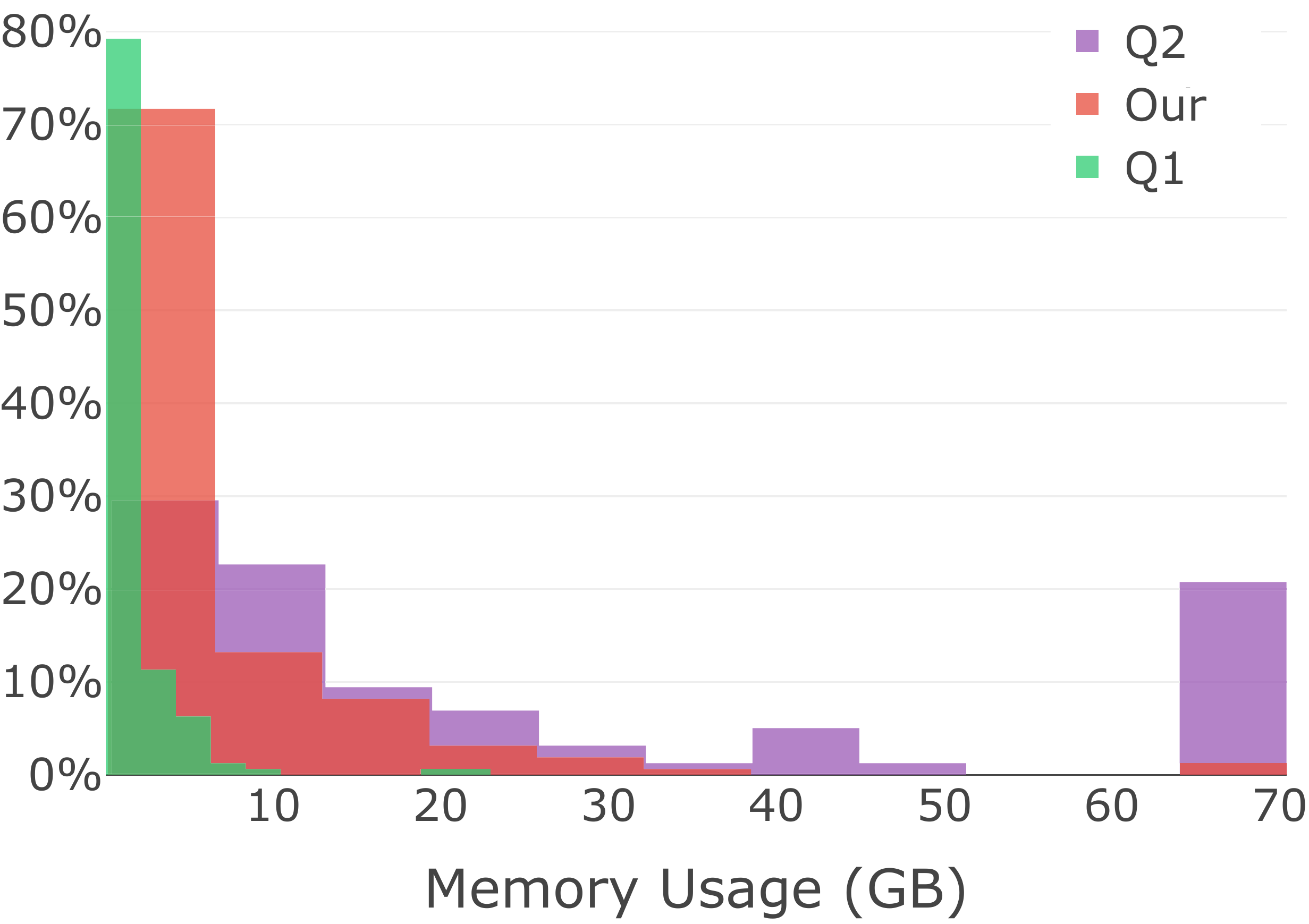}
	}
	\caption{Peak memory for the direct solver as reported by Pardiso. Left: 2D results. right: 3D results.}
	\label{fig:memory}
\end{figure}

\paragraph{Memory.} A histogram of the memory consumption of the solver is presented in \Cref{fig:memory}. The figure shows the peak memory usage as reported by Pardiso \cite{pardiso-6.0a, pardiso-6.0b} when solving the linear system arising from the FEM.
Out of the 159 pure hexahedral models we tested, 33 went out of memory when solving using $Q_2$ elements, while only 2 are too big to solve with our spline bases. On the star-shaped hybrid meshes, one model is too big to solve for both $Q_2$ and our spline construction.
More detailed statistics are reported in \Cref{tab:statistics}. We remark that the error for our method is higher than $Q_2$ because our method has less dofs (50\% less in average) since both meshes have the same number of vertices.

\begin{figure}
	\centering
	\makebox[\columnwidth][c]{
	\rotatebox{90}{\parbox{0.3\linewidth}{\centering\footnotesize Time (s)}}\hfill
	\begin{overpic}[width=0.52\linewidth]{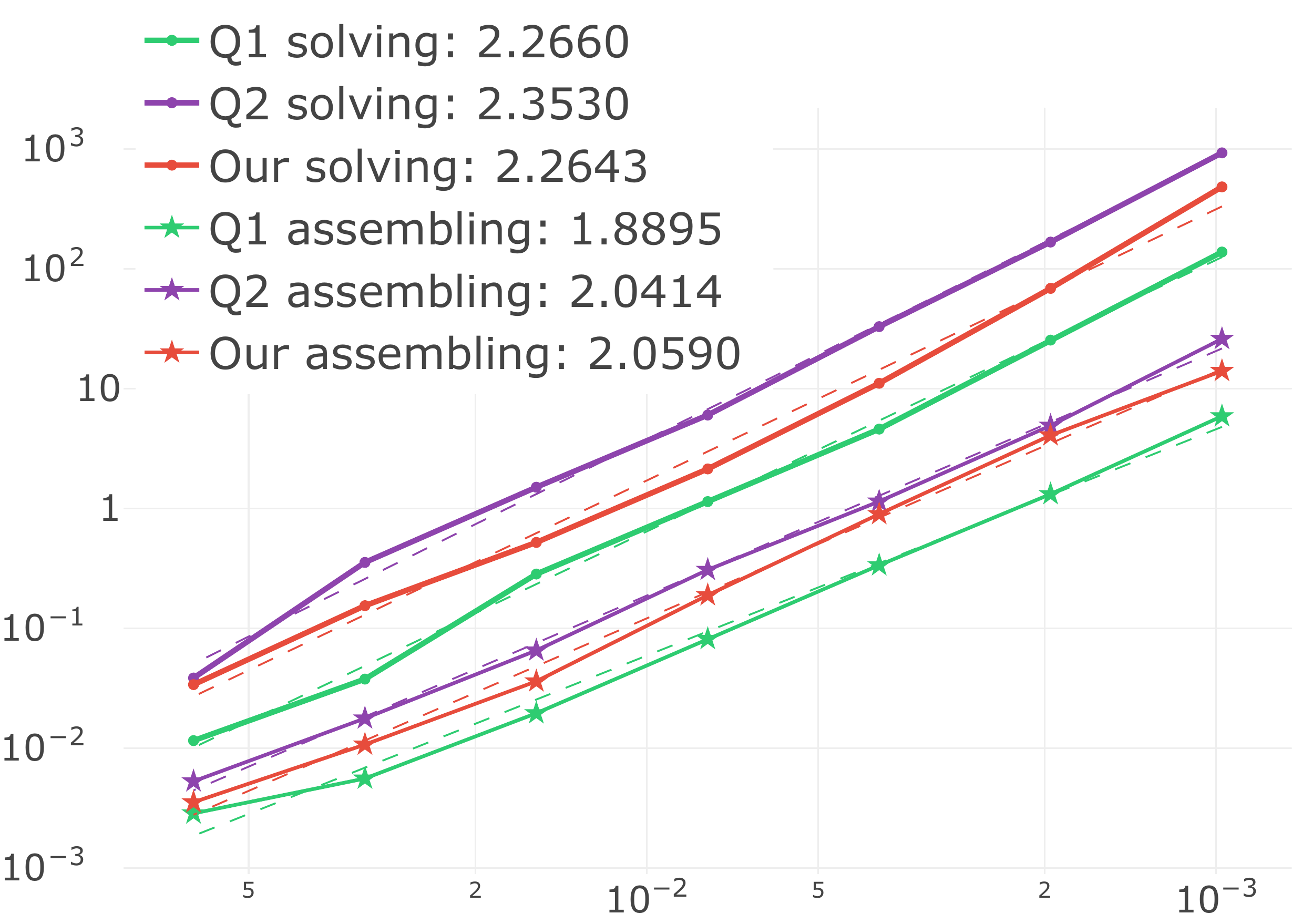}
	\end{overpic}\hfill{}
	\begin{overpic}[width=0.52\linewidth]{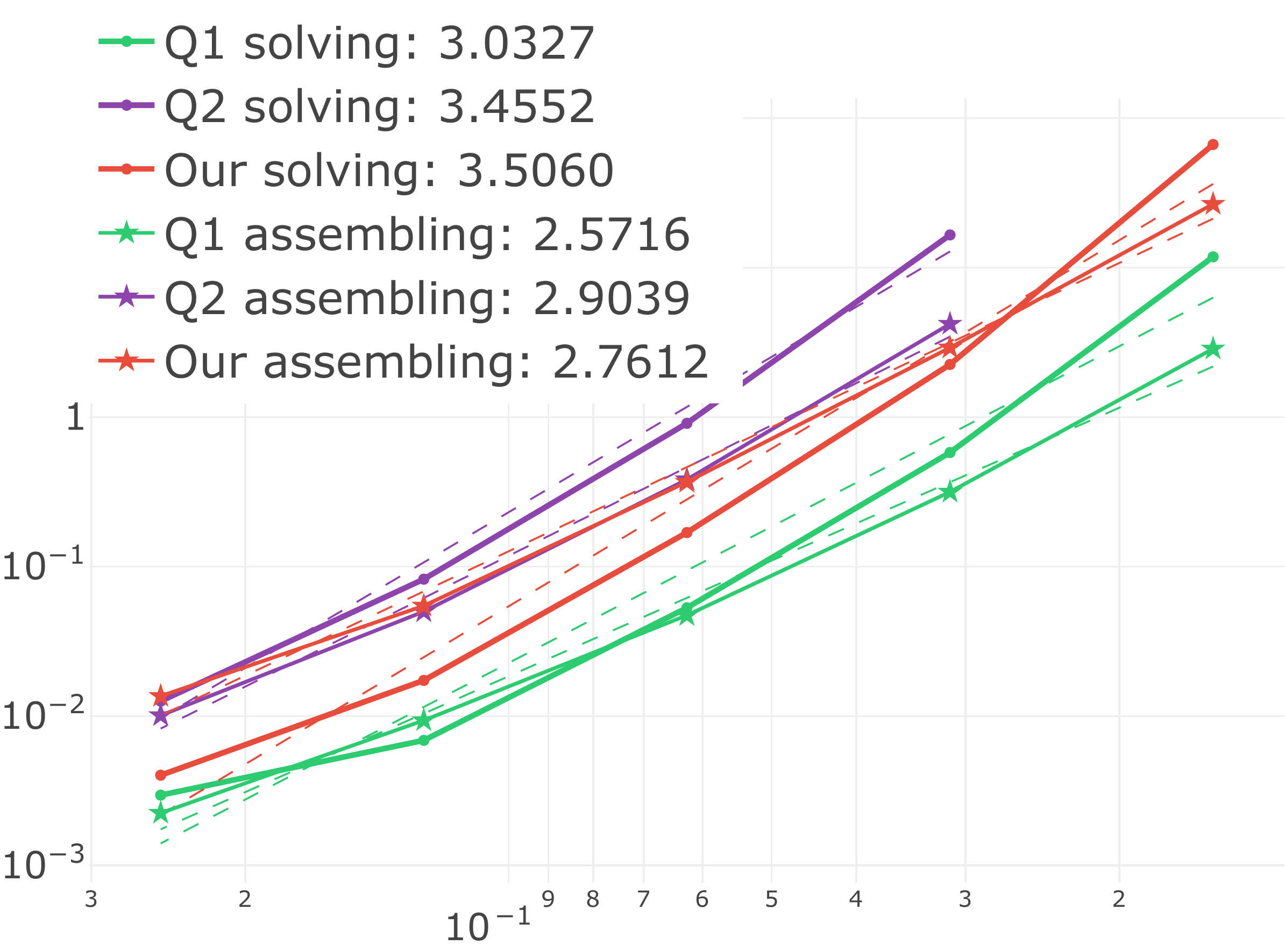}
	\end{overpic}
	}\\
	{\footnotesize Max edge length}
	\caption{
	    Time required to assemble the stiffness matrix and solve the linear system on a regular grid in 2D (left) and 3D (right).
    }
	\label{fig:timings_regular}
\end{figure}

\paragraph{Time.}  \Cref{fig:timings_regular} shows the  assembly time and solve time for solving a Poisson problem on an unit square (cube) under refinement in two (three) dimensions. Note that both steps (assembly and solve) are performed in parallel. For the 2D experiment we used a 3.1~GHz Intel Core i7-7700HQ with 8 threads, while in 3D we used a 3.5~GHz Intel Core i7-5930K with 12 threads (both machines use hyper-threading). In \Cref{tab:statistics} we summarize the timings for the large dataset using a 2.6~GHz Intel Xeon E5-2690v4 with 8 threads. In all cases, the total time is dominated by the solving time.

\begin{figure*}
	\centering
	\makebox[\columnwidth][c]{
	\rotatebox{90}{\parbox{0.15\linewidth}{\centering\footnotesize Error}}\hfill
	\begin{overpic}[width=0.26\linewidth]{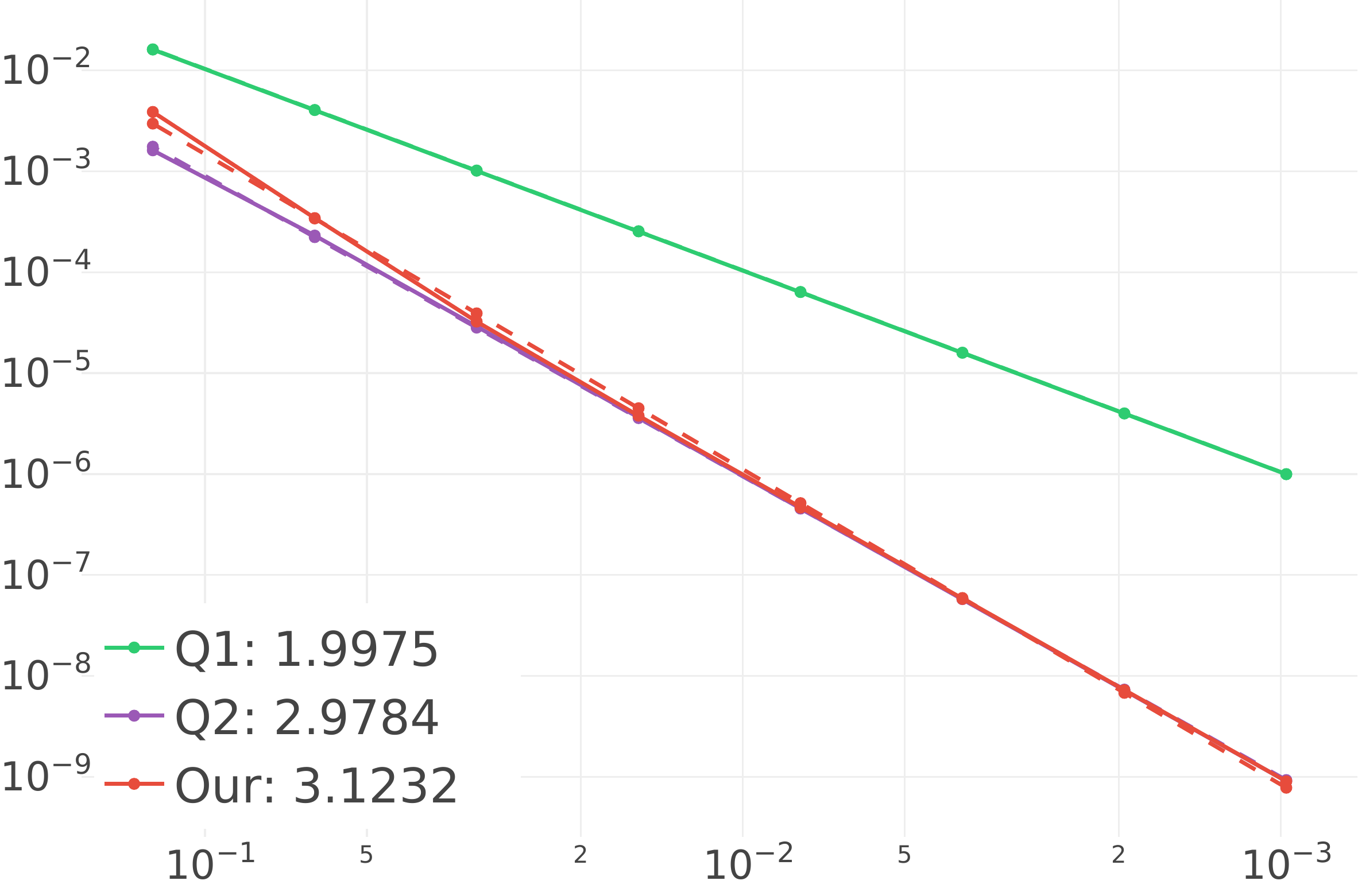}
		\put (50,60) {$L_2$}
	\end{overpic}\hfill{}
	\begin{overpic}[width=0.26\linewidth]{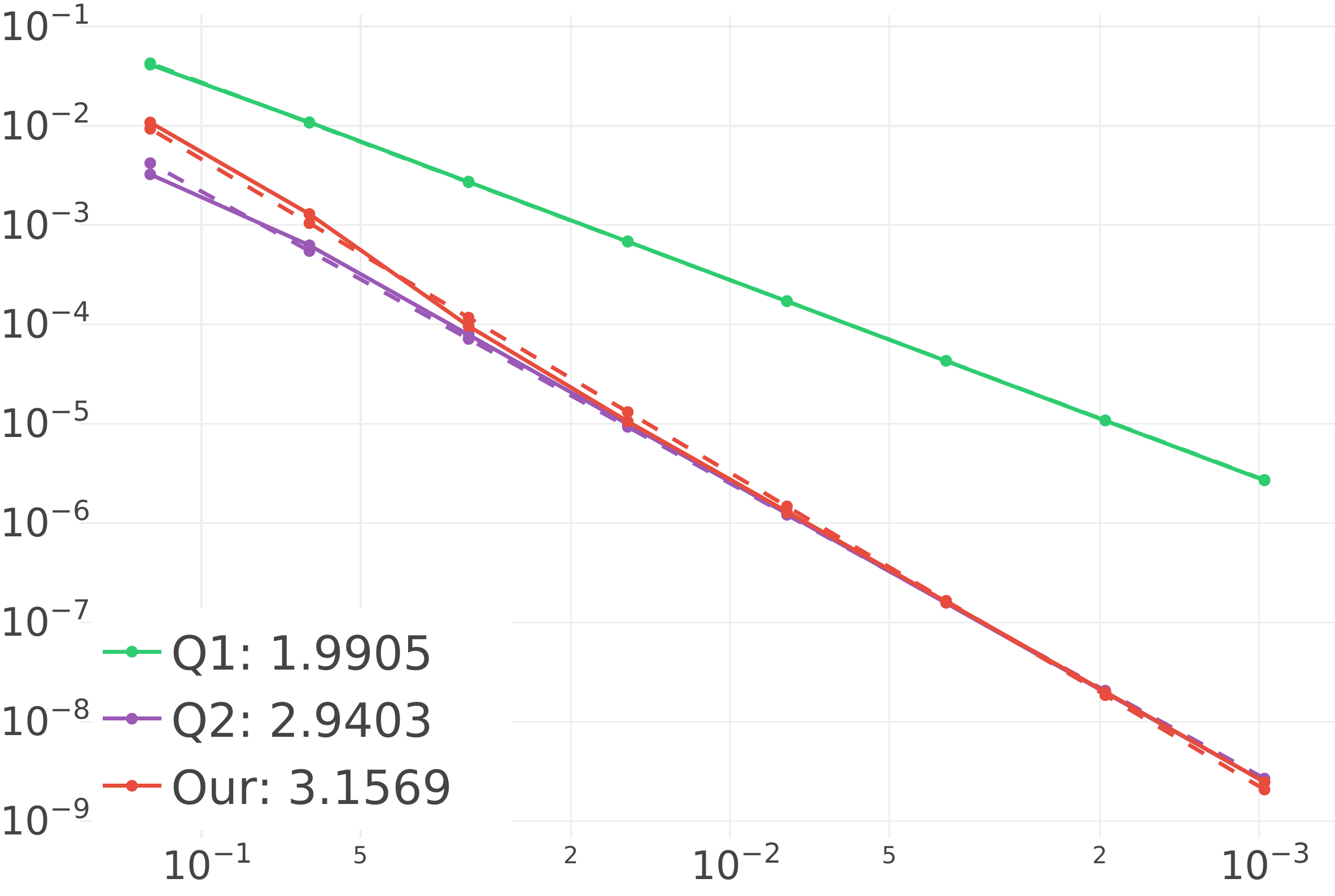}
		\put (50,60) {$L_\infty$}
	\end{overpic}\hfill{}
	\begin{overpic}[width=0.26\linewidth]{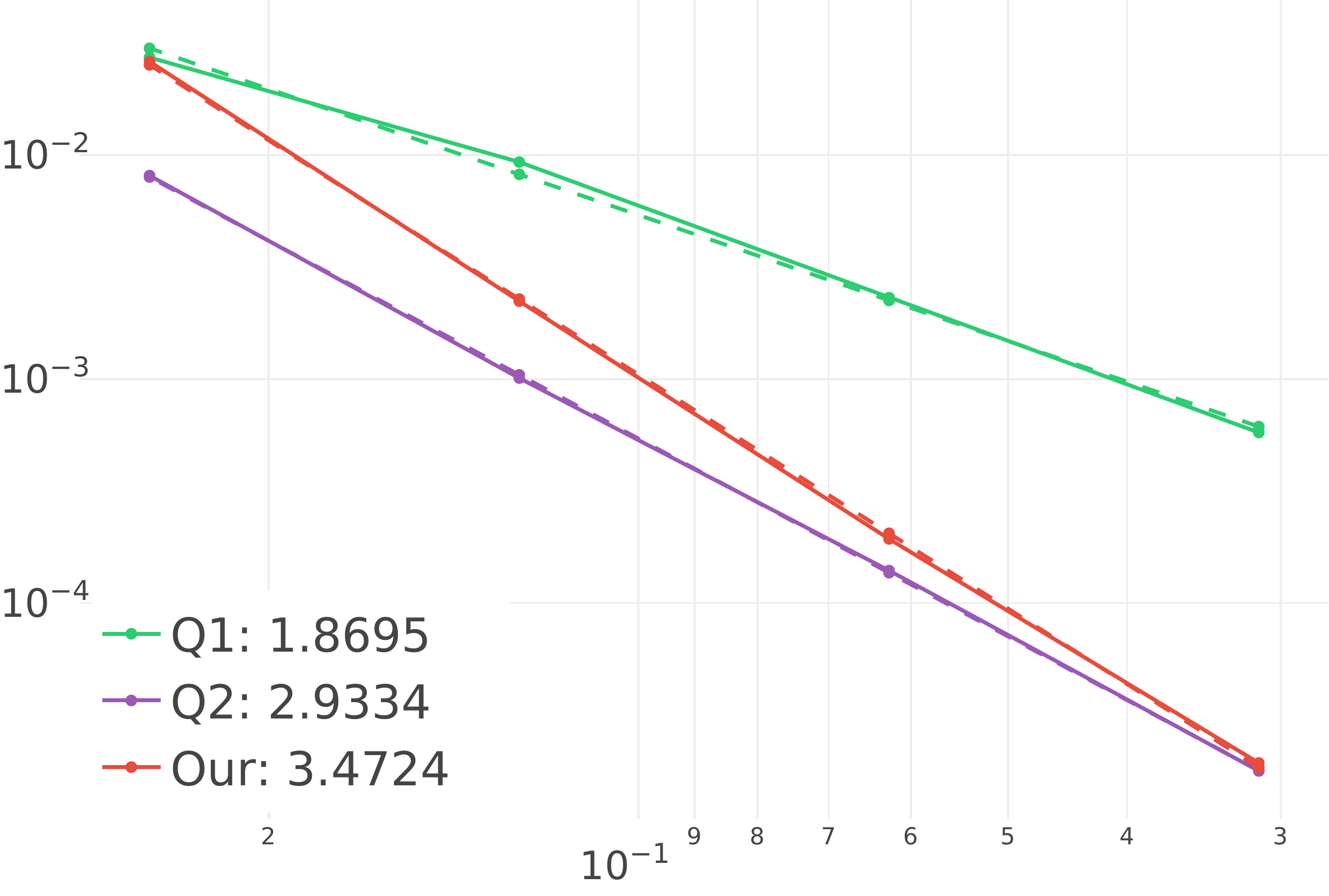}
		\put (50,60) {$L_2$}
	\end{overpic}\hfill{}
	\begin{overpic}[width=0.26\linewidth]{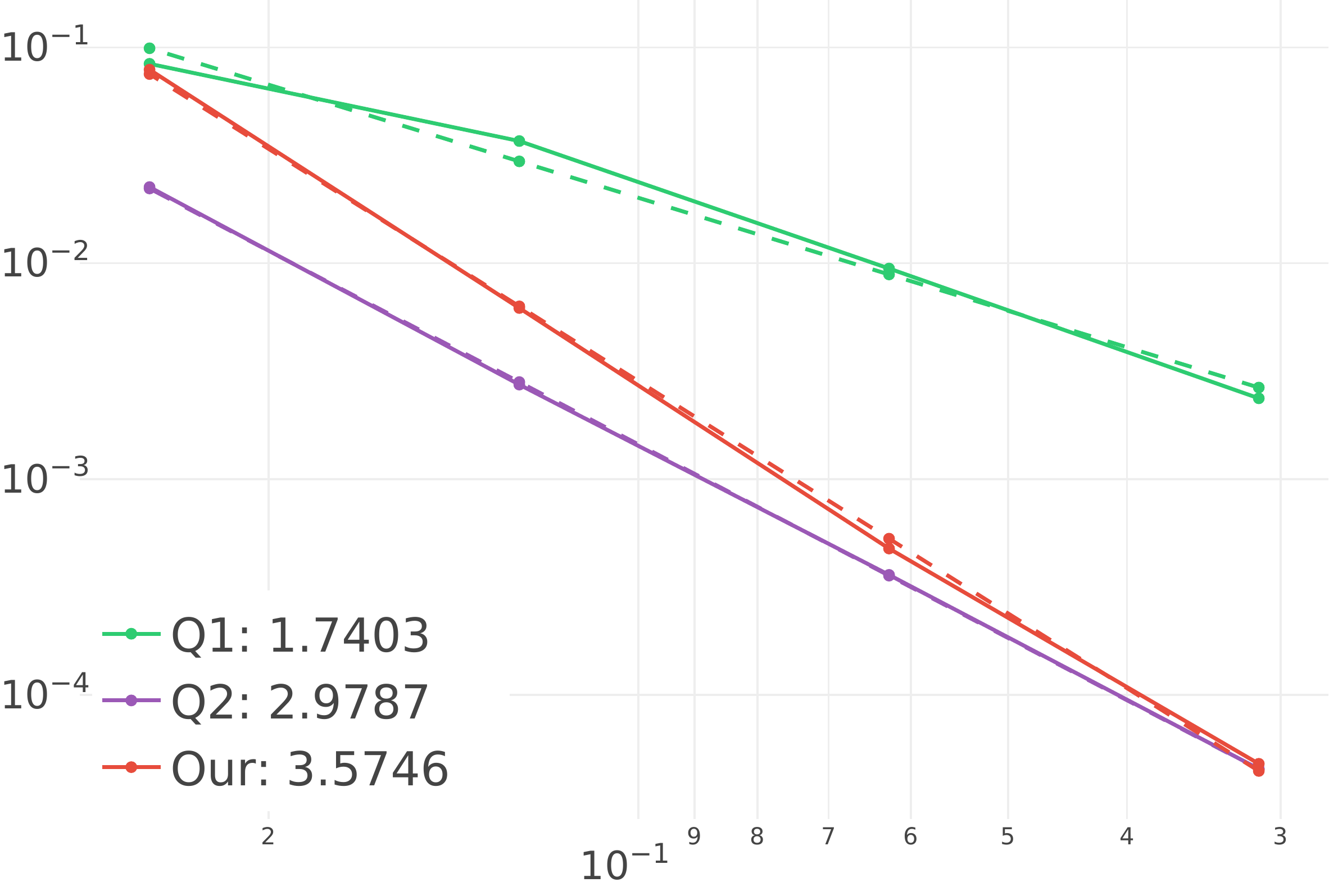}
		\put (50,60) {$L_\infty$}
	\end{overpic}\hfill{}
	}\\
	{\footnotesize Max edge length}

	\caption{
	    Poisosn equation convergence plot in $L_2$ and $L_\infty$ norm on a regular grid in 2D (first two) and 3D (last two). 
    }
	\label{fig:convergence_regular}
\end{figure*}
\begin{figure*}
	\centering
	\makebox[\linewidth][c]{
	\rotatebox{90}{\parbox{0.15\linewidth}{\centering\footnotesize Error}}\hfill
	\begin{overpic}[width=0.26\linewidth]{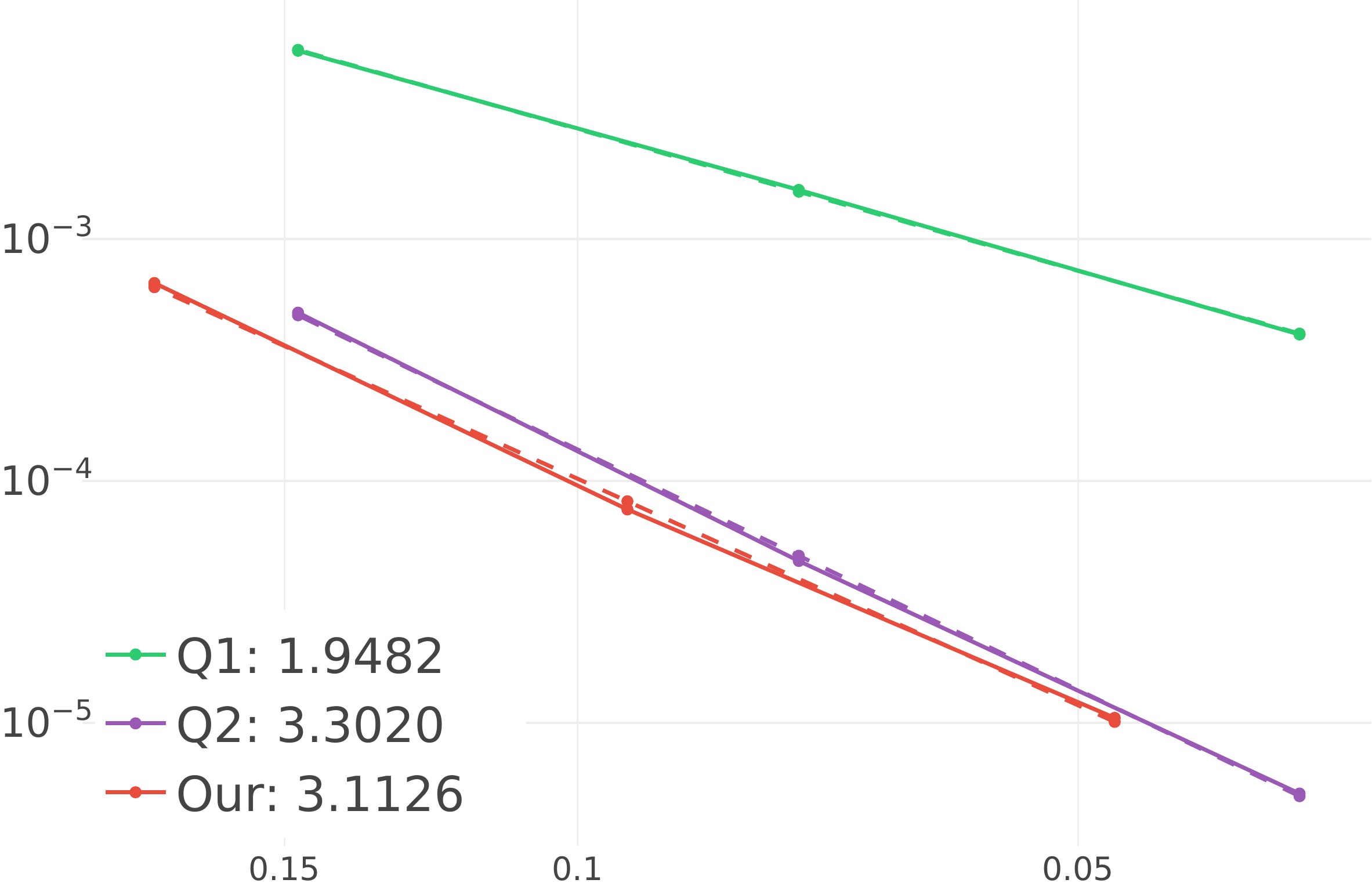}
		\put (50,60) {$L_2$}
	\end{overpic}\hfill{}
	\begin{overpic}[width=0.26\linewidth]{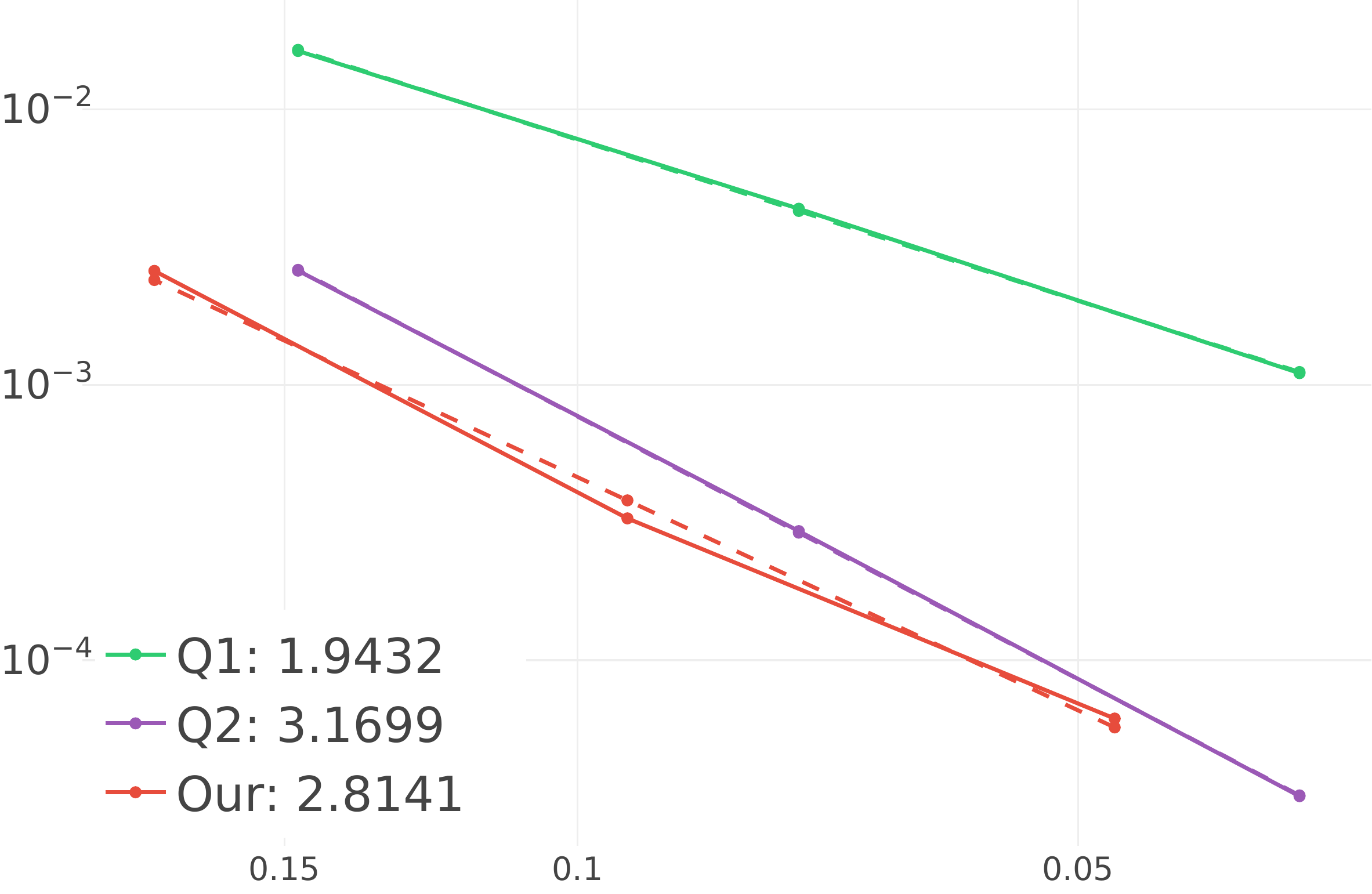}
		\put (50,60) {$L_\infty$}
	\end{overpic}\hfill{}
	\begin{overpic}[width=0.26\linewidth]{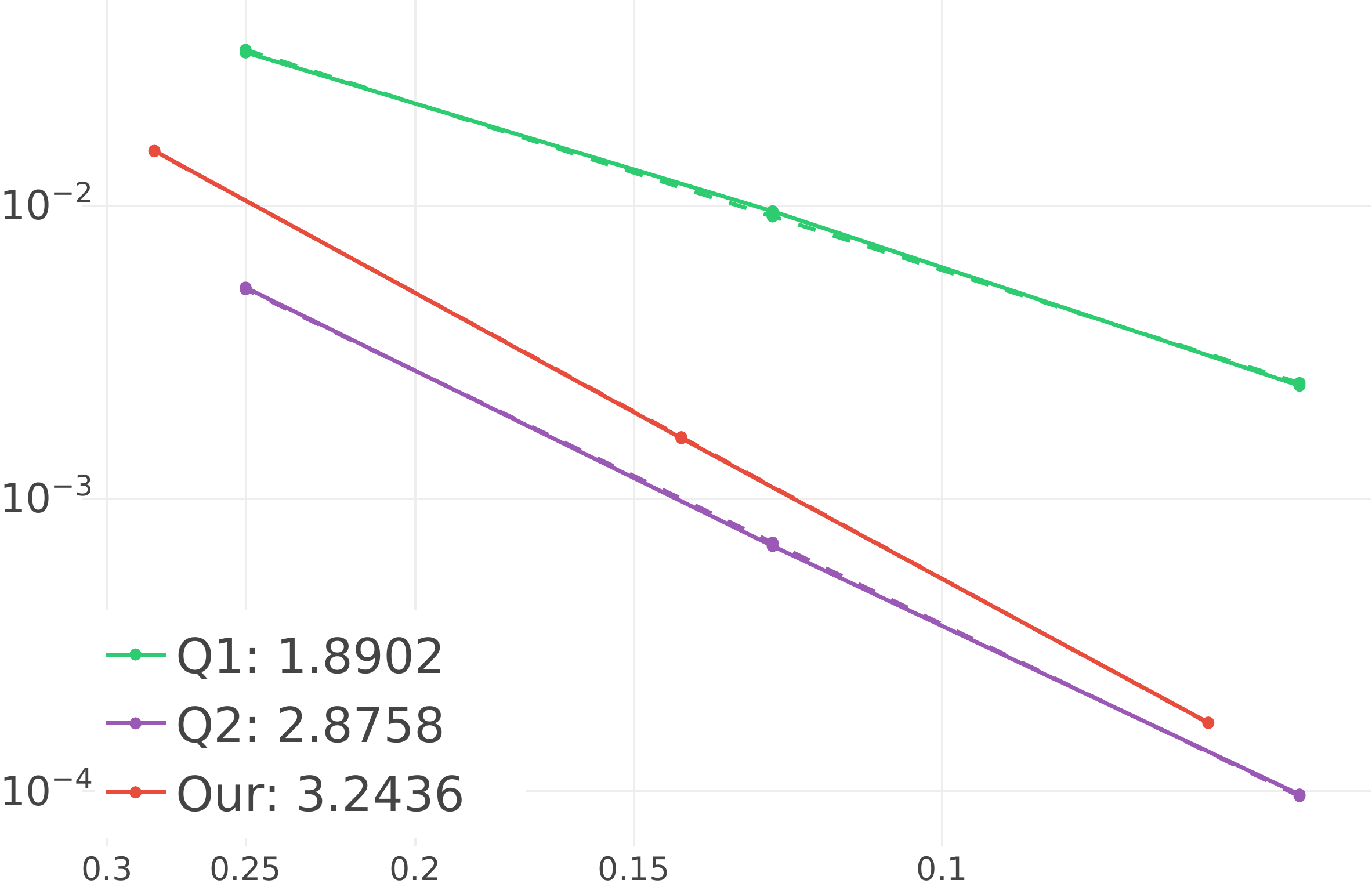}
		\put (50,60) {$L_2$}
	\end{overpic}\hfill{}
	\begin{overpic}[width=0.26\linewidth]{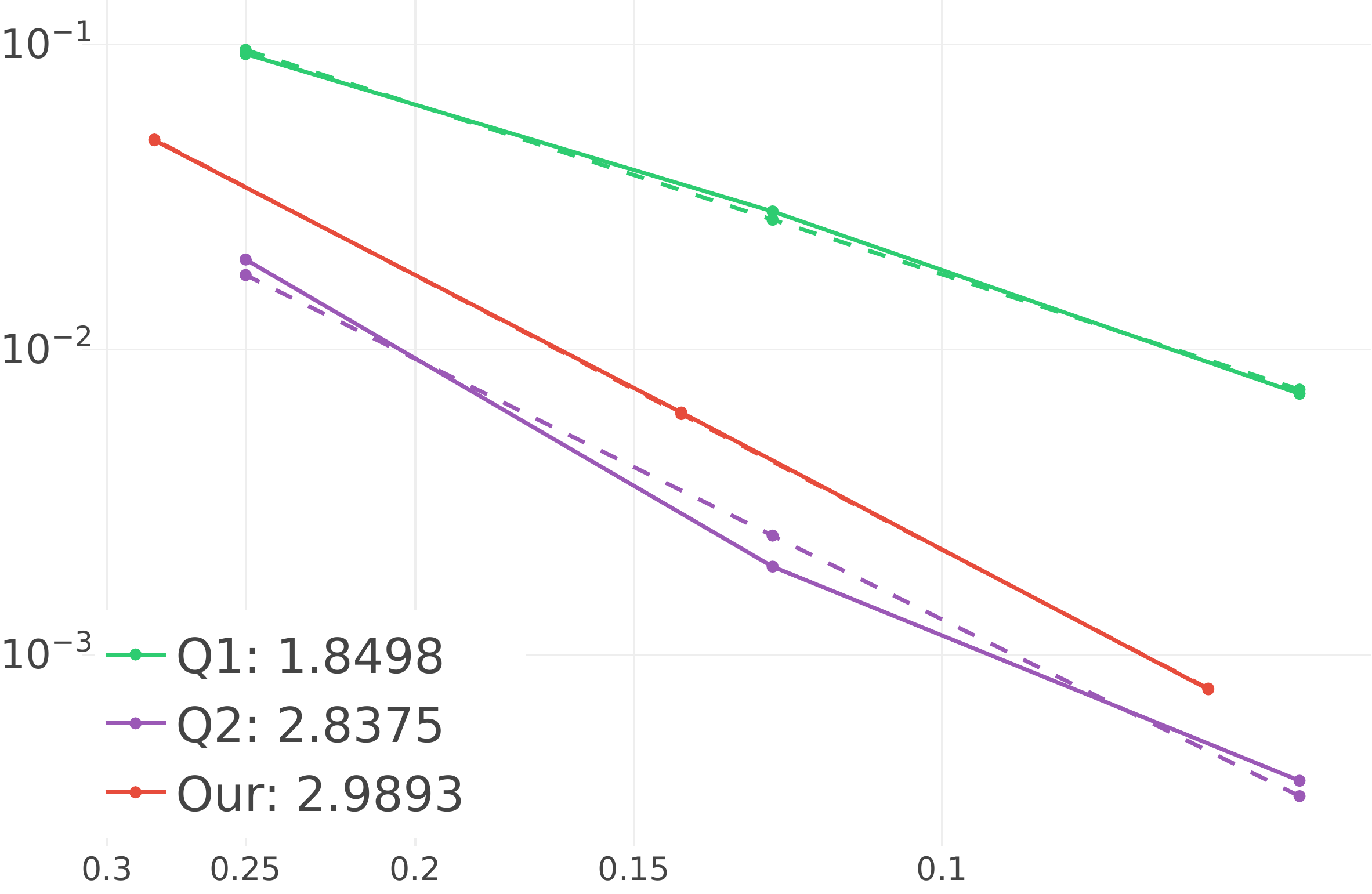}
		\put (50,60) {$L_\infty$}
	\end{overpic}\hfill{}
	}\\
	{\footnotesize Max edge length}

	\caption{
	    Poisson equation convergence plot in $L_2$ and $L_\infty$ norm for a hybrid mesh in 2D (first two) and 3D (last two). Meshes are show in \Cref{fig:refinement_example}. 
	}
	\label{fig:convergence_hybrid}
\end{figure*}

\begin{figure}
	\centering
	\makebox[\columnwidth][c]{
	\rotatebox{90}{\parbox{0.3\linewidth}{\centering\footnotesize Error}}\hfill
	\begin{overpic}[width=0.52\linewidth]{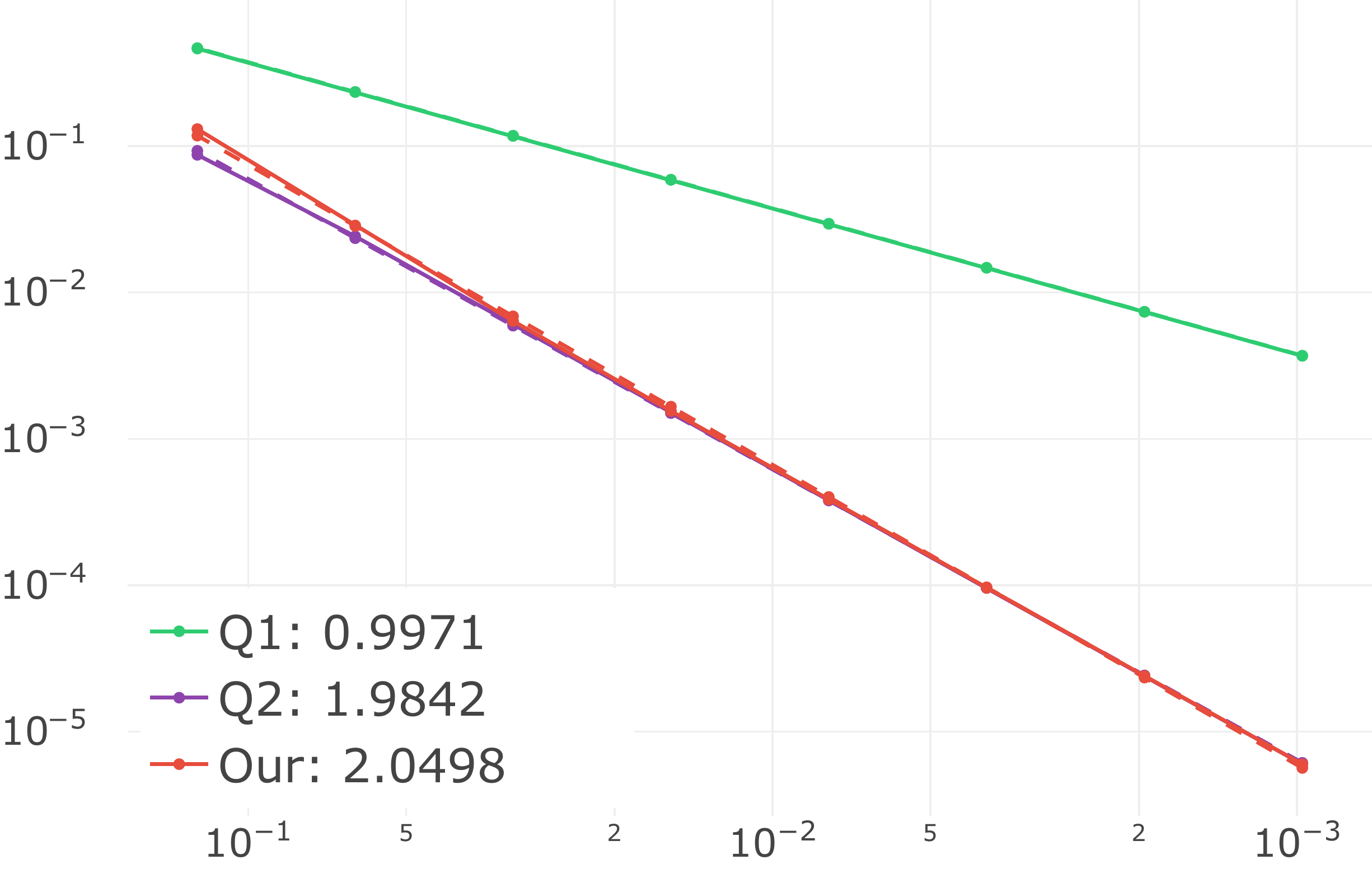}
	\end{overpic}\hfill{}
	\begin{overpic}[width=0.52\linewidth]{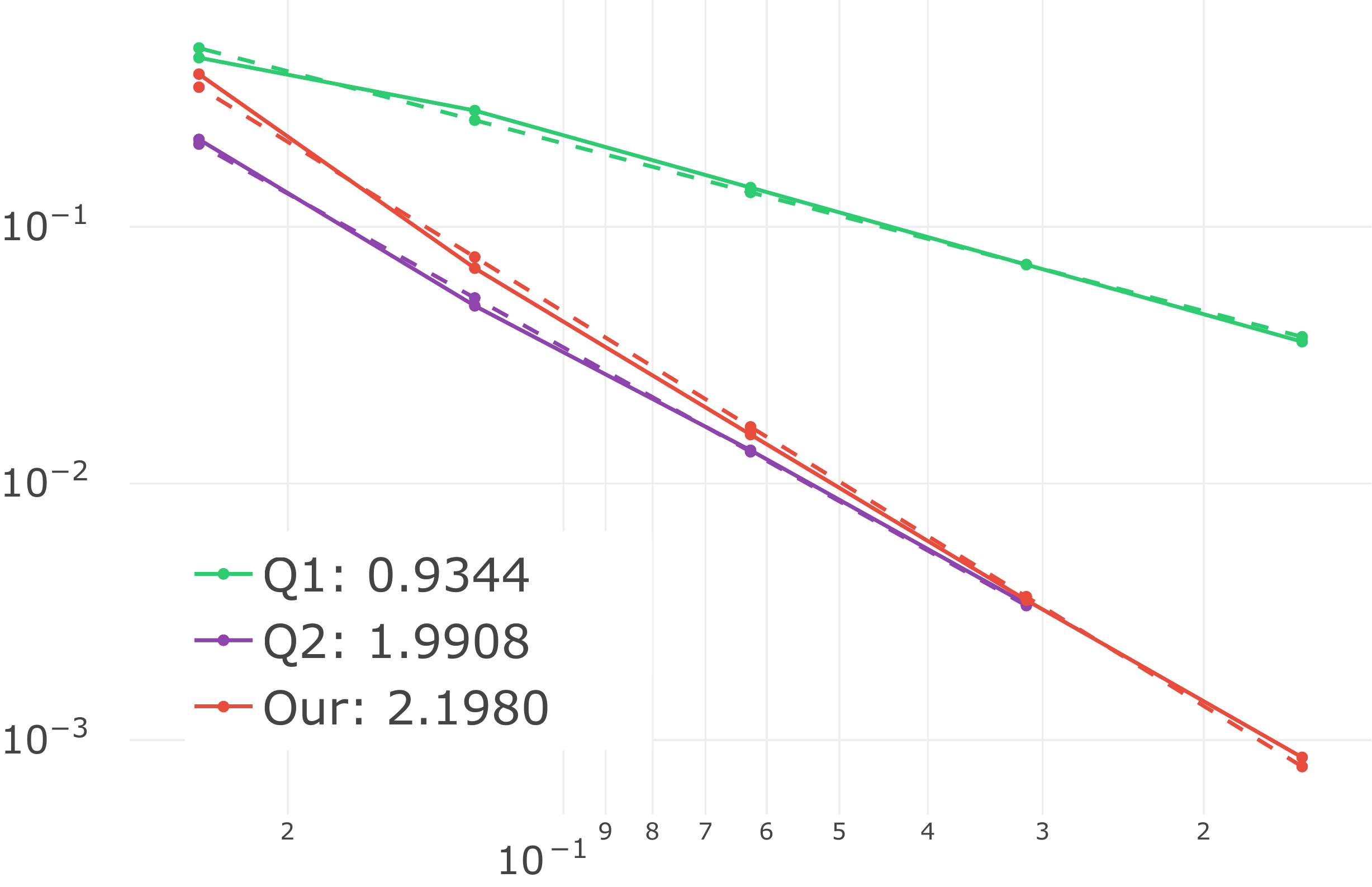}
	\end{overpic}
	}\\
	{\footnotesize Max edge length}
	\caption{
	    Poisson equation convergence plot in $H_1$ norm on a regular grid in 2D (first plot) and 3D (second plot). 
    }
	\label{fig:convergence_h1_regular}
\end{figure}

\paragraph{Convergence.} Figures~\ref{fig:convergence_regular} and~\ref{fig:convergence_h1_regular} show the convergence of spline elements vs $Q_1$ and $Q_2$ for the $L_2$, $L_\infty$, and $H_1$ norms, in the ideal case of a uniform grid, both in 2D and 3D.
This is in a sense the best-case scenario that can be expected for our spline construction: every element is regular and has a $3^2$ or $3^3$ neighborhood. In this situation, splines exhibit a superior convergence $> 3.0$ under both $L_2$, $L_\infty$, and $H_1$ norms.

On a 2D test mesh with mixing polygons and splines (model shown in \Cref{fig:refinement_example} top), we achieved a convergence rate of 2.8 in $L_\infty$, and 3.1 in $L_2$ (\Cref{fig:convergence_hybrid}, left). \Cref{fig:convergence_hybrid} also displays the convergence we obtained on a very simple hybrid 3D mesh, starting from a cube marked as a polyhedron, to which we applied our polar refinement described in Section~\ref{sec:refinement}. On this particular mesh, the splines exhibited a $L_\infty$ convergence similar to $Q_2$, albeit producing an error that is somewhat larger.

\begin{figure}
	\centering
	\rotatebox{90}{\parbox{0.3\linewidth}{\centering\footnotesize Error}}
	\begin{overpic}[width=0.52\linewidth]{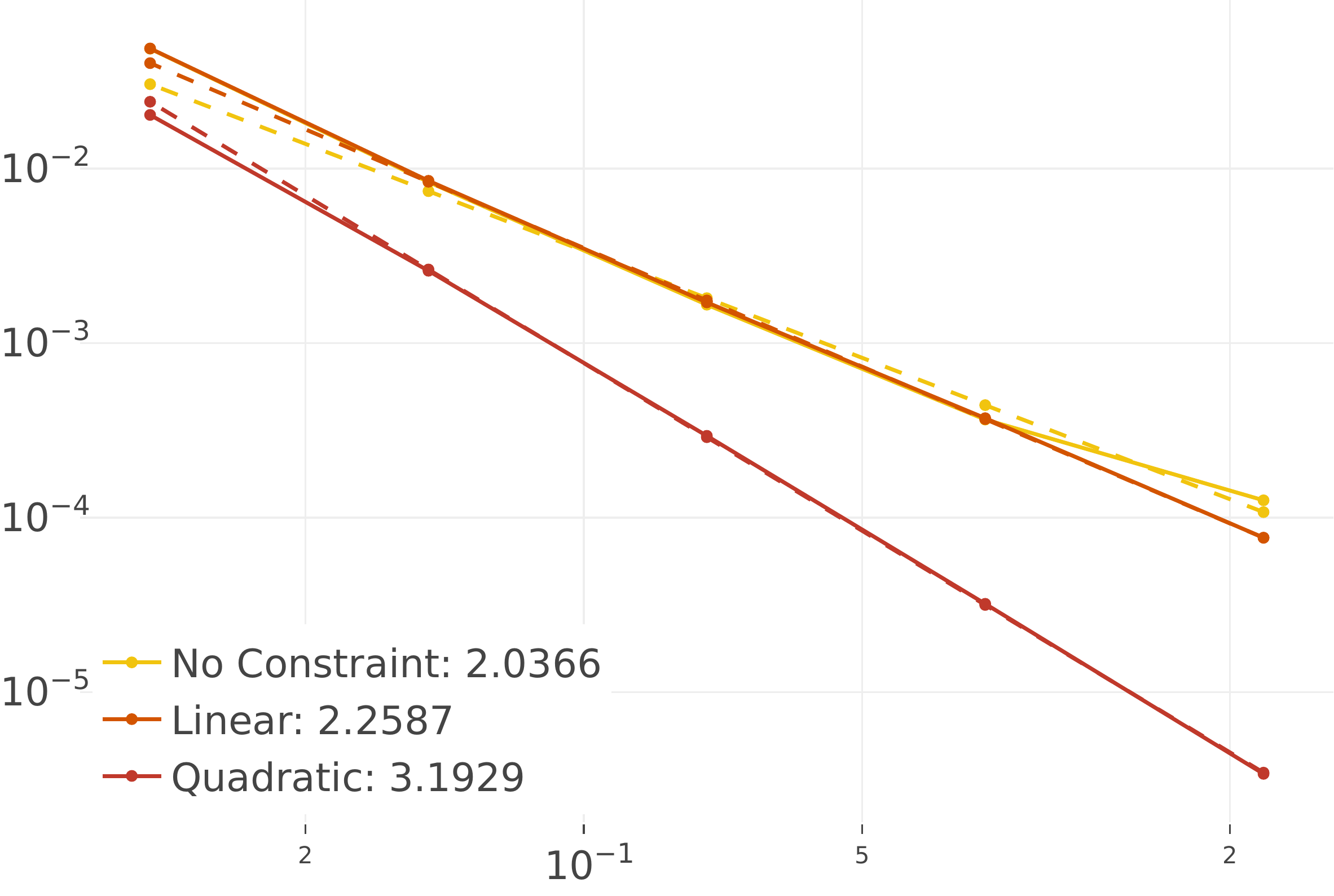}
		\put (50,60) {$L_\infty$}
	\end{overpic}
	\\
	{\footnotesize Max edge length}
	\caption{$L^\infty$ convergence for the different consistency constraints on the polyhedron of Figure~\ref{fig:refinement_example}.}
	\label{fig:integral_constraints}
\end{figure}

\paragraph{Consistency Constraints.} \Cref{fig:integral_constraints} shows the effect of our consistency constraint on the convergence of a polygonal mesh under refinement (the one shown in \Cref{fig:refinement_example}, top), with $Q_2$ elements used on the quadrilateral part. Without imposing any constraint on the bases overlapping the polygon, one can hope at best for a convergence of $\sim{}2.0$, whereas pure $Q_2$ elements should have a convergence rate of $3.0$. With a constraint ensuring linear reproduction for the bases defined on polyhedra, the convergence rate is still only $\sim{}2.5$.
Finally, with the constraints we describe in Section~\ref{sec:polyhedral} to ensure the bases reproduce triquadratic polynomials, we reach the expected convergence rate of $\sim{3.0}$.

\begin{figure}
	\centering
	\includegraphics[width=0.24\linewidth]{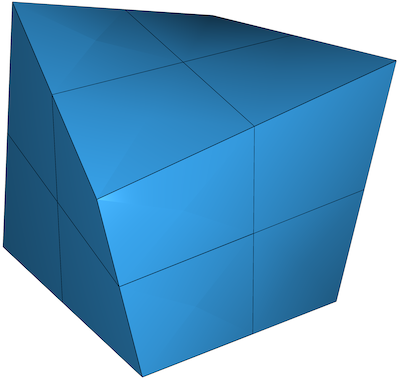}\hfill{}
	\includegraphics[width=0.24\linewidth]{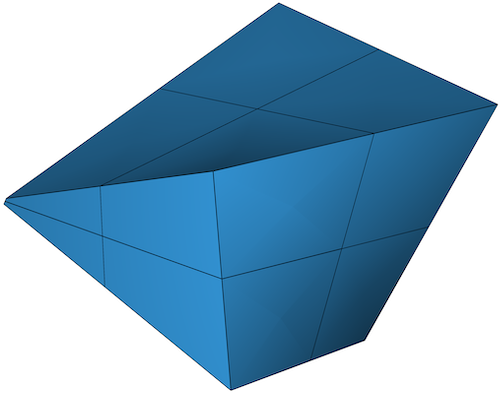}\hfill{}
	\includegraphics[width=0.24\linewidth]{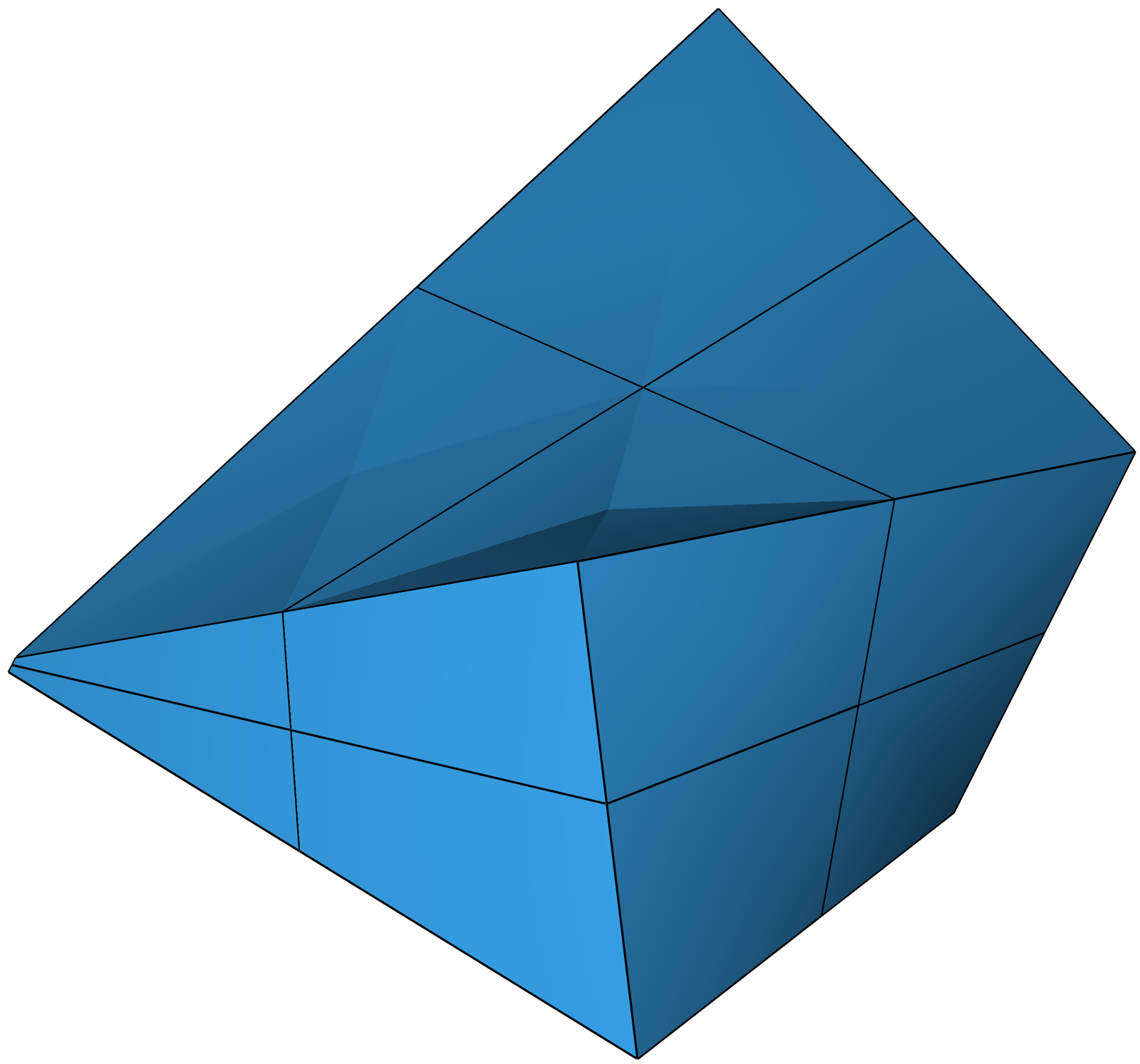}\hfill{}
	\includegraphics[width=0.24\linewidth]{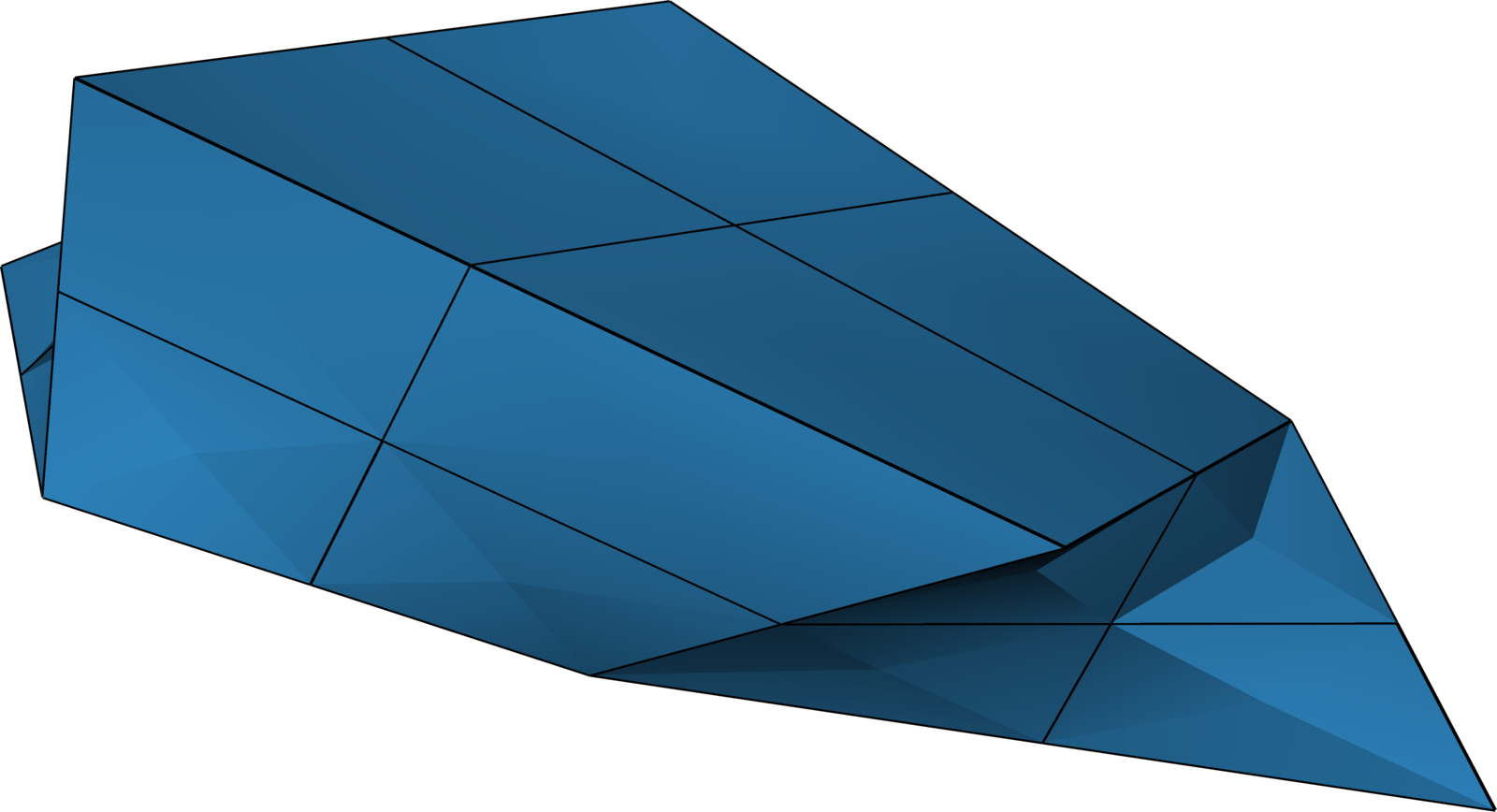}\par
	\caption{Low-quality polyhedra used to evaluate the interpolation errors.}
	\label{fig:poly_error}
\end{figure}

\paragraph{Polyhedral Basis Resilience.} Our polyhedral bases are less susceptible to badly shaped elements than $Q_2$. We computed the $L_2$ and $L_\infty$ interpolation errors for the gradients of the Franke function for 14 badly shaped hexahedra, \Cref{fig:poly_error} shows some of them. The $L_2$ and $L_\infty$ maximum and average errors are 3 times smaller with our polygonal basis.

\begin{figure}
	\centering
	\rotatebox{90}{\parbox{0.3\linewidth}{\centering\footnotesize Condition number}}\hfill
	\includegraphics[width=0.45\linewidth]{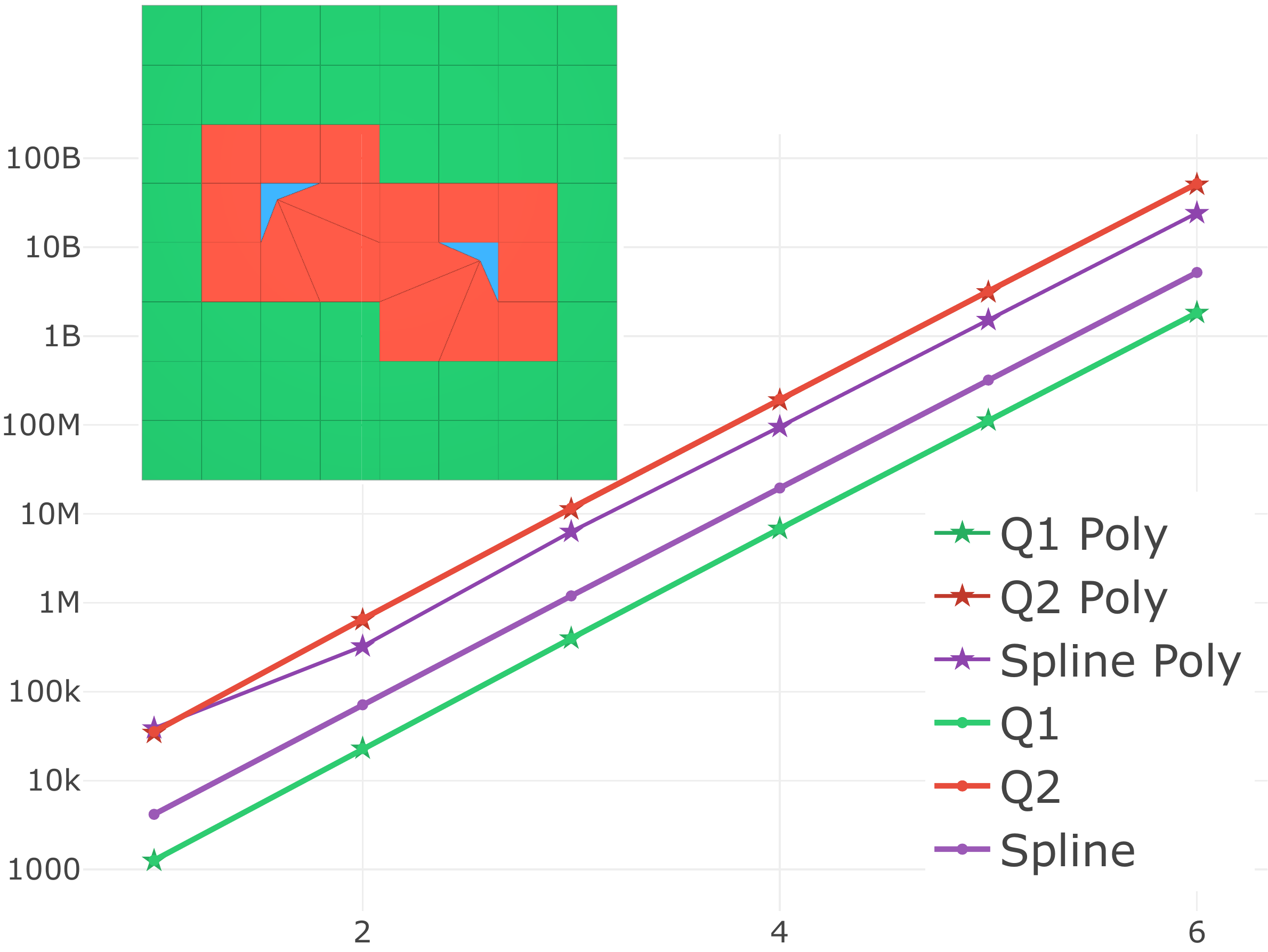}\hfill{}
	\includegraphics[width=0.45\linewidth]{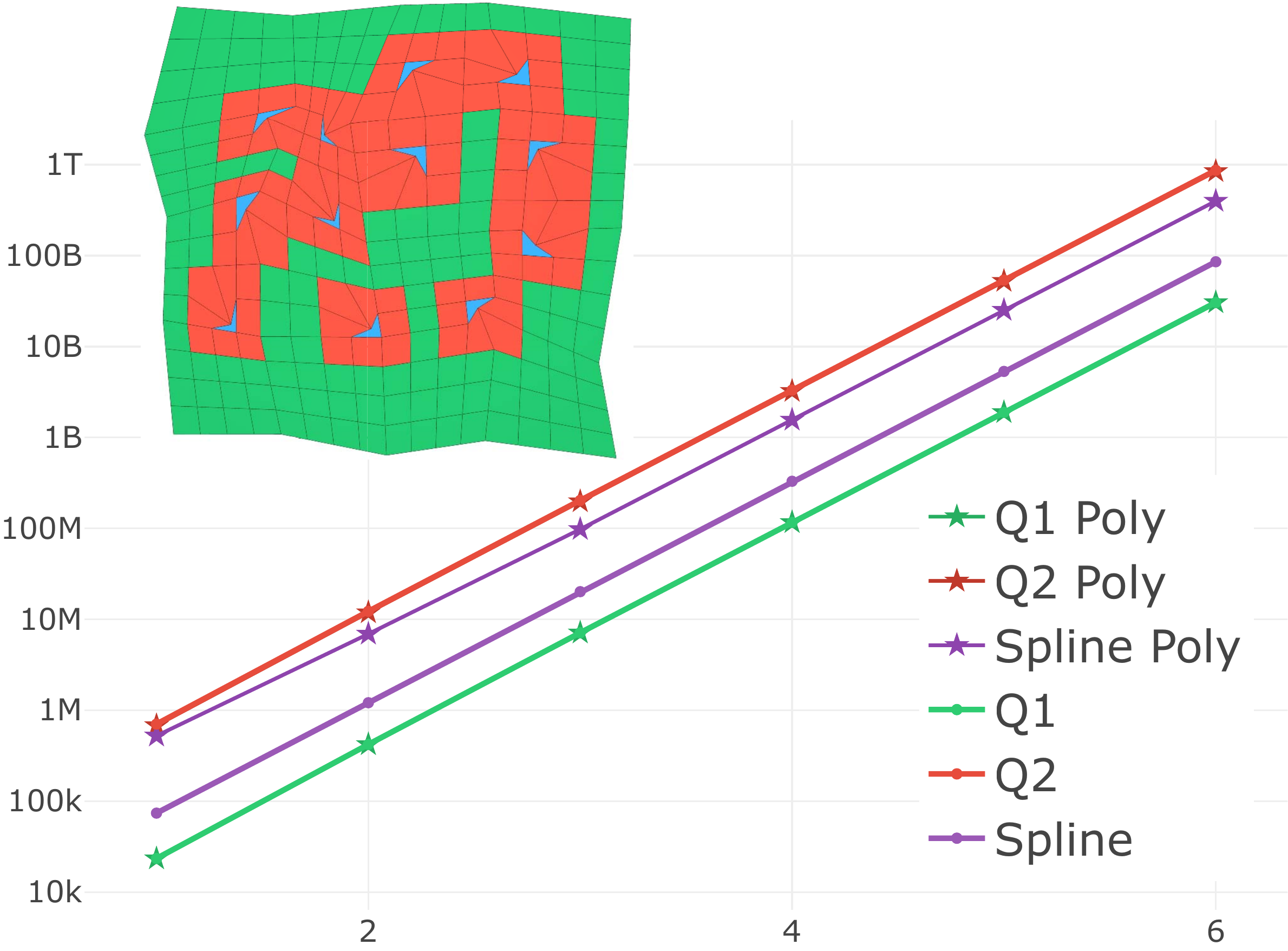}\\
	{\footnotesize Number of refinements}
	\caption{
	Evolution of the condition number of the stiffness matrix for the Poisson problem under refinement.
	For each level of refinement we artificially marked 5\% of the quads as polyhedra and move one random vertex on the diagonal between 20\% 40\%, as shown in the insets figures in blue.
 Note that some of the curves coincide, that is $Q_1$ with $Q_1$ poly and $Q_2$ with $Q_2$ poly.
	}
	\label{fig:spectrum}
\end{figure}

\paragraph{Conditioning and Stability.} An important aspect of our new FE method is the conditioning of the resulting stiffness matrix: this quantity relates to both the stability of the method, and to its performances when an iterative linear solver is used (important only for large problems where direct solvers cannot be used due to their memory requirements). We compute the condition number of the Poisson stiffness matrix on a regular and perturbed grid (Figure \ref{fig:spectrum}). 
In both cases, our discretization has a good conditioning number, slightly higher than pure linear elements, but lower than pure quadratic elements (while sharing the same cubic convergence property).
To evaluate the conditioning of the polyhedral bases we started from a base mesh of good quality, marked $5\%$ of the quads as polygons, and pushed one of the vertices inwards. Even for this extreme distortion of polyhedral elements, the conditioning remained similar to the case when no polyhedral elements are used on the same mesh.

\begin{figure}
	\centering
	\hfill{}
	\includegraphics[width=0.45\linewidth]{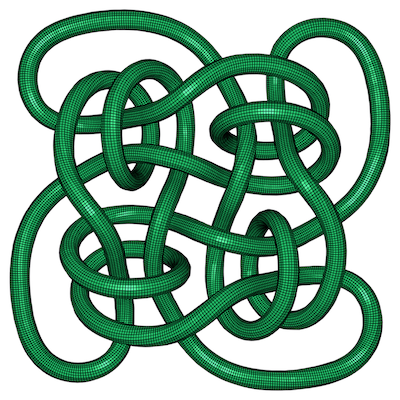}\hfill{}
	\includegraphics[width=0.45\linewidth]{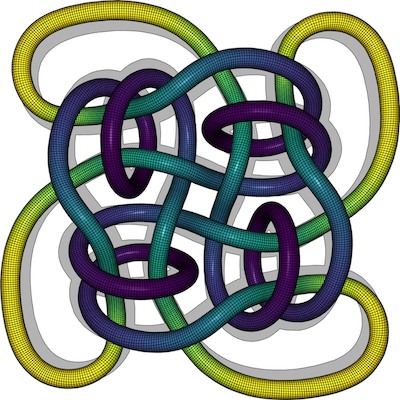}\par
	\includegraphics[width=0.45\linewidth]{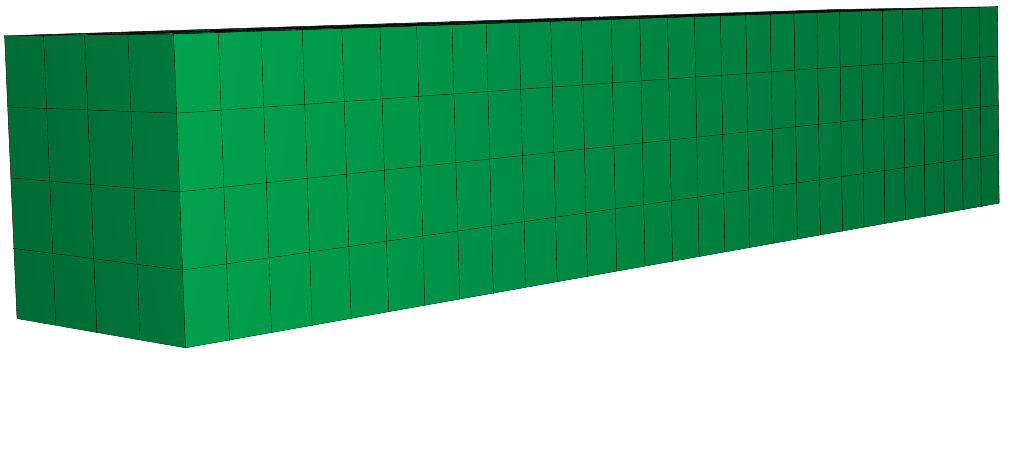}\hfill
	\includegraphics[width=0.45\linewidth]{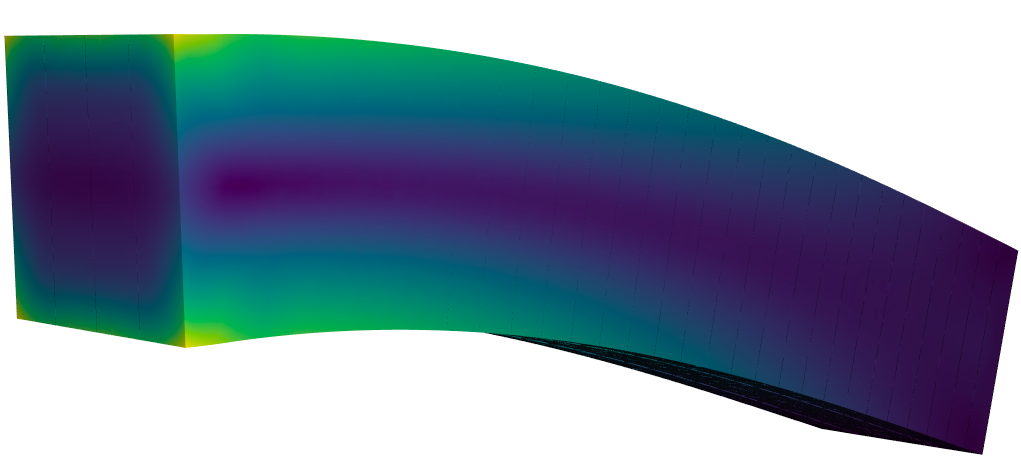}\par
	\caption{Displacements computed solving linear elasticity on a pure hexahedral 3D model, using spline bases. Top a complicated model with $\lambda=1$ and $\mu=1$, bottom a bended bar $\nu=0.35$ and large young modulus $E=2$e5.}
	\label{fig:elasticity}
\end{figure}

\paragraph{Elasticity.}
While most of our testing was done for the Poisson equation, we have performed some testing of linear elasticity problems.
\Cref{fig:elasticity} top shows the solution of a linear elasticity problem on a pure hexahedral mesh. The outer loops of the knots are pulled outside of the figure, deforming the knot. The color in the figure represents the magnitude of the displacement vectors. On the bottom we show the result for a Young's modulus of 2e5.

\begin{figure}
    \centering
	\makebox[\columnwidth][c]{
	\rotatebox{90}{\parbox{0.3\linewidth}{\centering\footnotesize Error}}\hfill
	\begin{overpic}[width=0.52\linewidth]{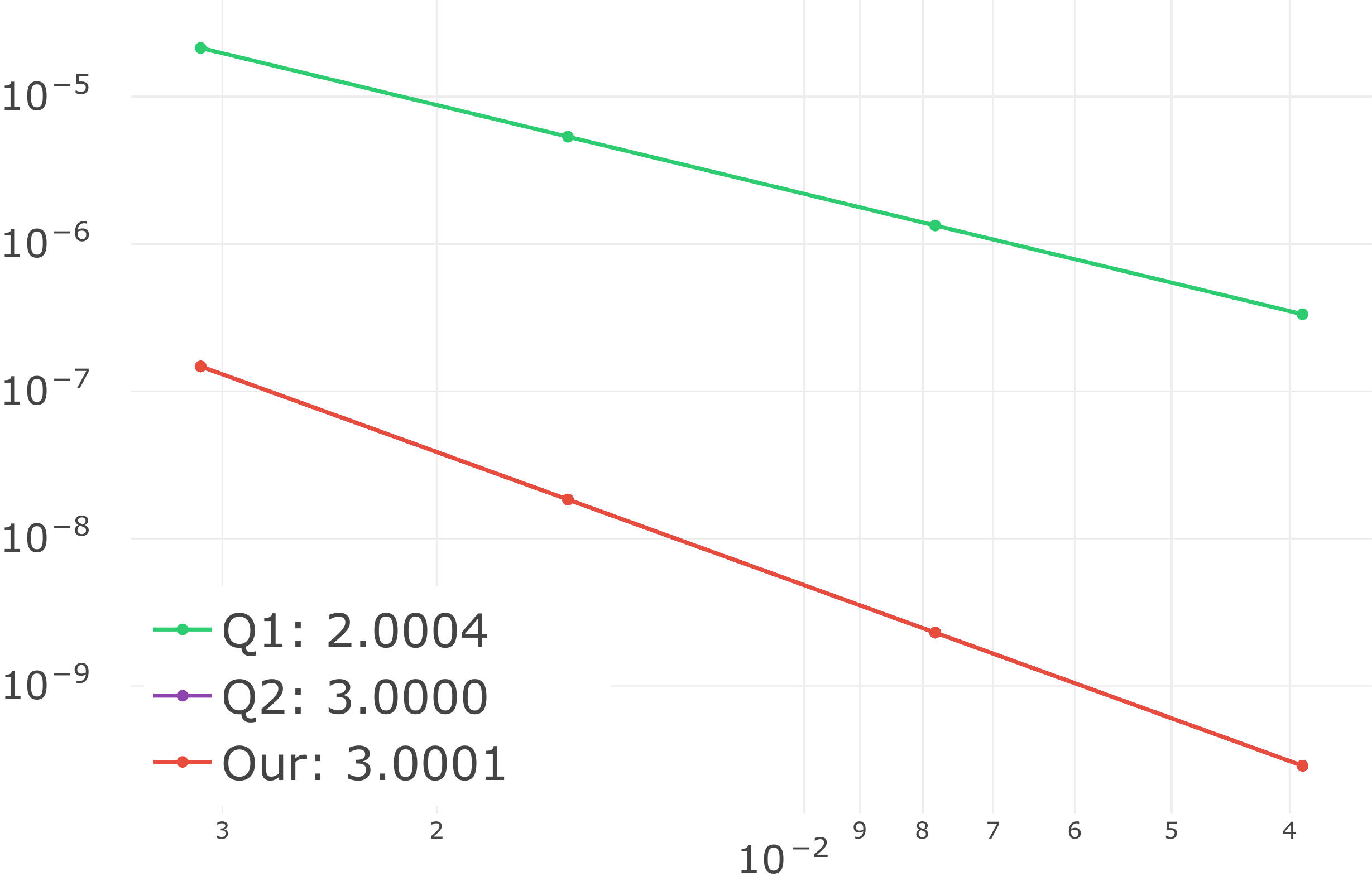}
	\put (50,60) {$L_2$}
	\end{overpic}\hfill{}
	\begin{overpic}[width=0.52\linewidth]{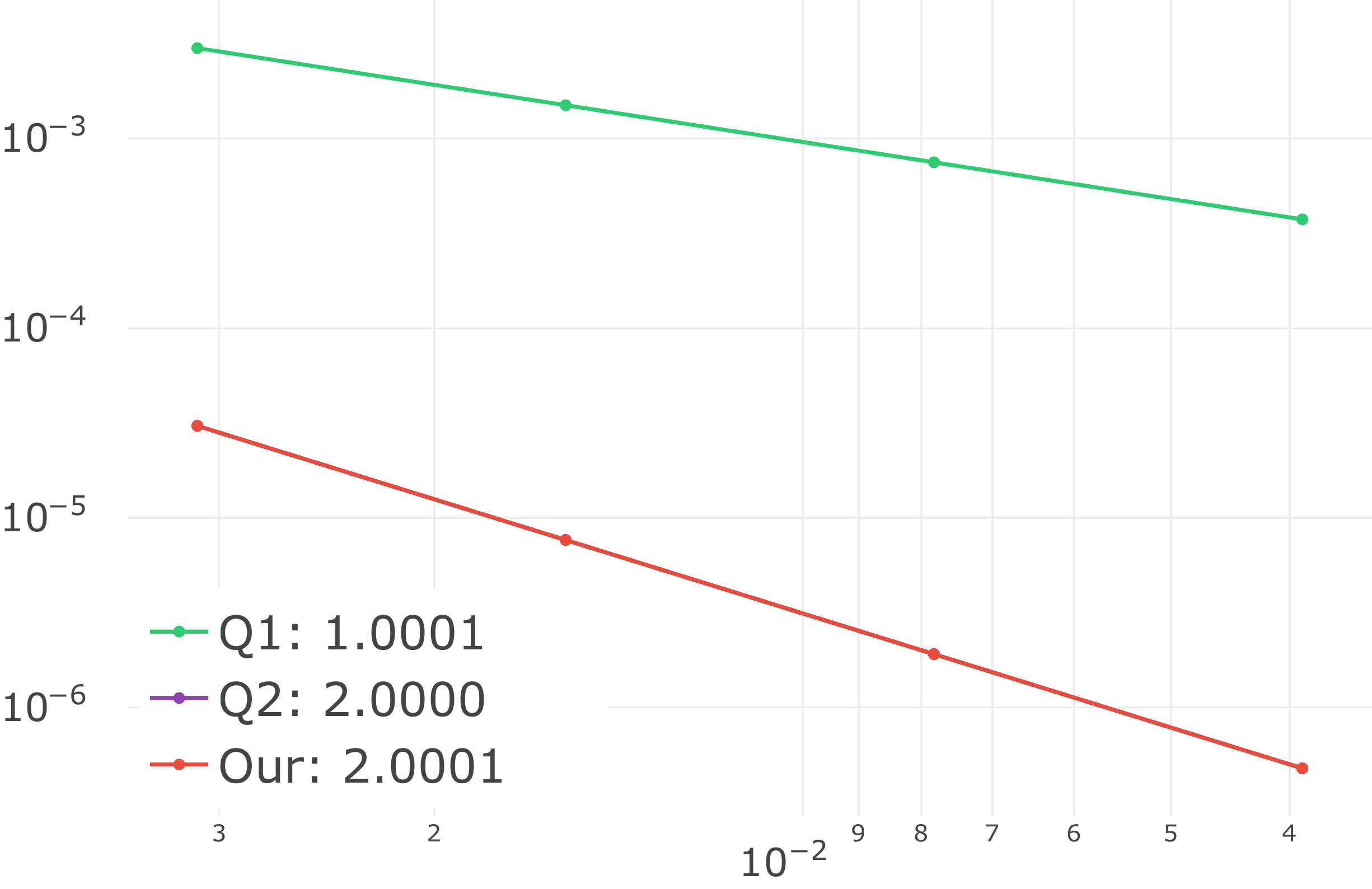}
		\put (50,60) {$H_1$}
	\end{overpic}
	}\\
	{\footnotesize Max edge length}
	\caption{
	    Linear elasticity convergence plot in $L_2$ and $H_1$ norm on a regular grid in 2D.
    }
	\label{fig:convergence_regular_elast}
\end{figure}

\begin{figure}
    \centering
	\makebox[\columnwidth][c]{
	\rotatebox{90}{\parbox{0.3\linewidth}{\centering\footnotesize Error}}\hfill
	\begin{overpic}[width=0.52\linewidth]{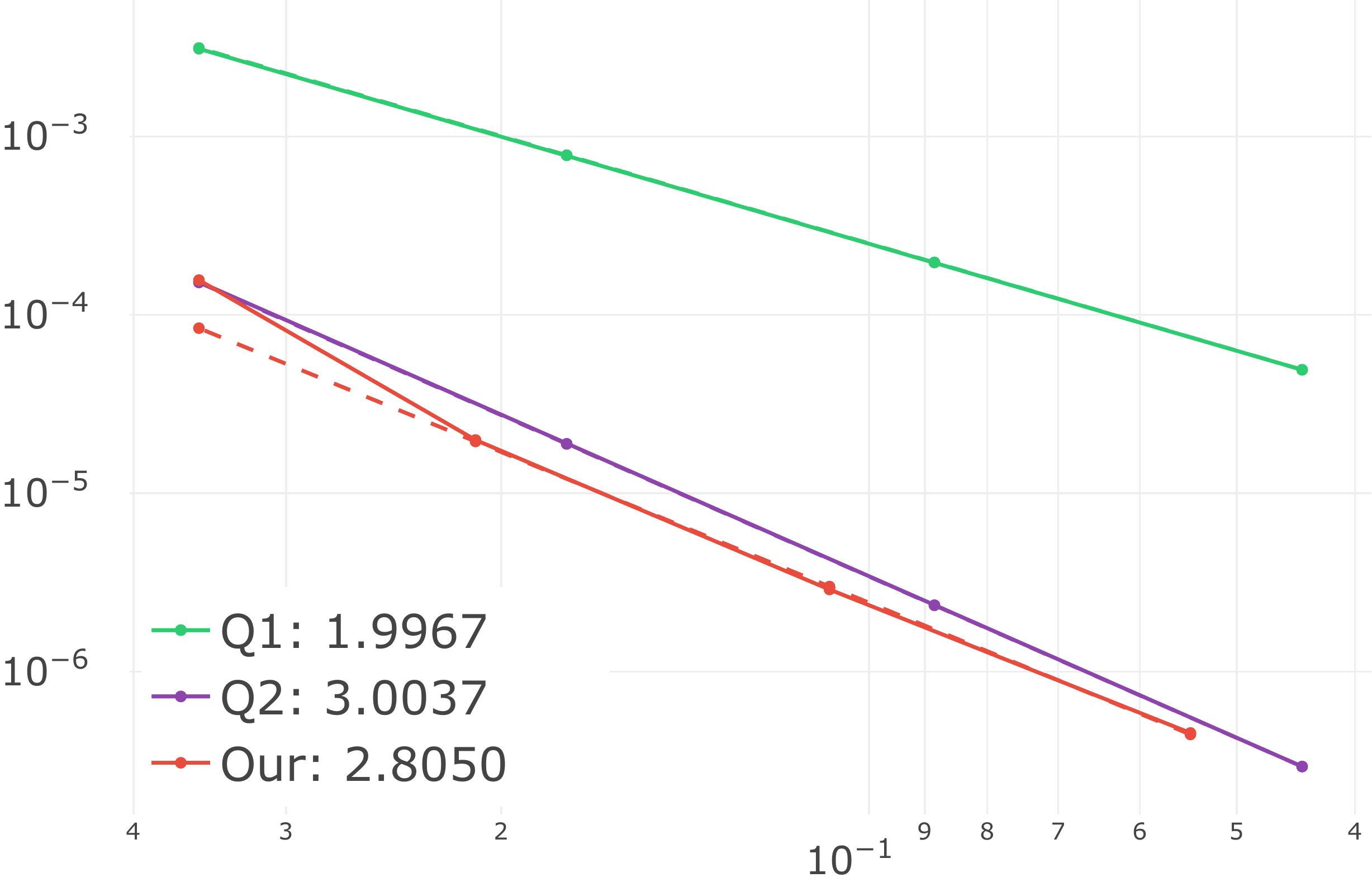}
	\put (50,60) {$L_2$}
	\end{overpic}\hfill{}
	\begin{overpic}[width=0.52\linewidth]{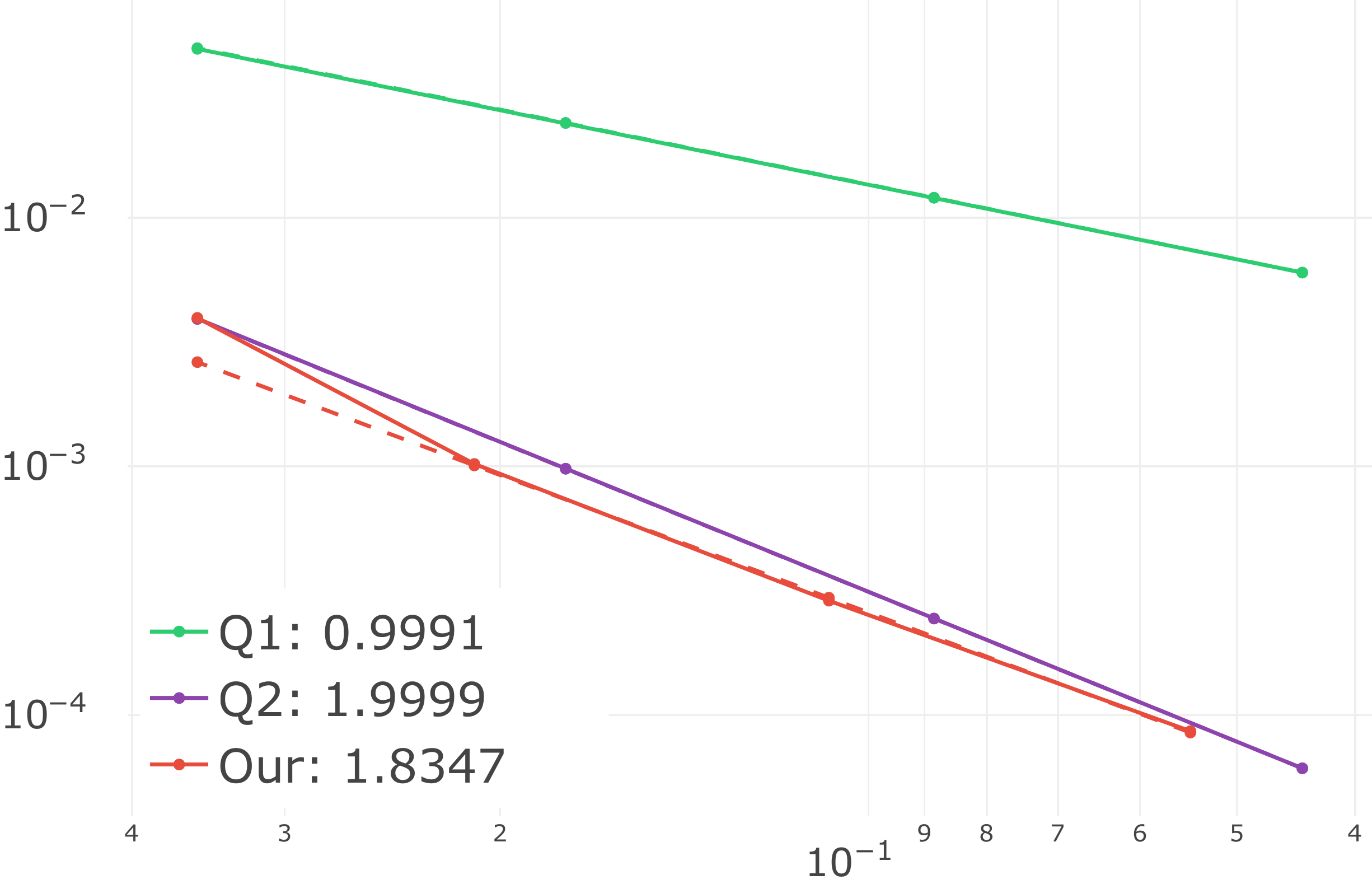}
		\put (50,60) {$H_1$}
	\end{overpic}
	}\\
	{\footnotesize Max edge length}
	\caption{
	    Linear elasticity convergence plot in $L_2$ and $H_1$ norm on a hybrid mesh in 2D.
    }
	\label{fig:convergence_hybrid_elast}
\end{figure}

Figure~\ref{fig:convergence_regular_elast} shows a plot for the linear elasticity PDE with Young's modulus 200 and Poisson's ratio 0.35 on a regular grid, and similar results are obtained an hybrid mesh, Figure~\ref{fig:convergence_hybrid_elast}.
The convergence plots for $Q_1$ and $Q_2$ are obtained by mixing regular $Q_1$/$Q_2$ bases with the polyhedral construction (Section~\ref{sec:polyhedral}).

\begin{table}
    \centering
    \makebox[\linewidth][c]{
    {\footnotesize
	\begin{tabular}{@{}c@{}l@{}rrr@{\quad}r@{\quad}rS[output-exponent-marker=\text{e}]S[output-exponent-marker=\text{e}]}
		   & & Num dofs & \multicolumn{1}{@{}c@{}}{Solver} & \multicolumn{1}{@{}c@{}}{Bases} & Assembly & \specialcell[c]{Memory \\ (\si{\mebi\byte})} & {$L_2$ Error} & {$L_\infty$ Error} \\
		\toprule

		\multirow{5}{*}{$Q_1$}
			& mean   & 174,335 & \ang{;  8;26} & \ang{; 15;40} & \ang{;0;37} & 1,132 & \num{7.60e-05} & \num{1.11e-03} \\
			& std    & 177,192 & \ang{; 20;29} & \ang{; 19;44} & \ang{;0;40} & 1,589 & \num{2.27e-04} & \num{2.69e-03} \\
			& min    &   3,035 & \ang{;  0; 0} & \ang{;  0;17} & \ang{;0;1 } &     5 & \num{5.57e-07} & \num{1.98e-05} \\
			& median & 105,451 & \ang{;  1; 2} & \ang{;  9;8 } & \ang{;0;23} &   500 & \num{3.18e-05} & \num{3.39e-04} \\
			& max    & 926,938 & \ang{;182;32} & \ang{;101;27} & \ang{;5;20} & 9,329 & \num{2.88e-03} & \num{2.09e-02} \\
		\midrule
		\multirow{5}{*}{$Q_2^\star$ }
			& mean   &   552,583 & \ang{; 63;43} & \ang{;13;11} & \ang{;0;55} &  5,716 & \num{3.62e-06} & \num{6.34e-05} \\
			& std    &   355,783 & \ang{; 72;42} & \ang{;11;17} & \ang{;0;59} &  4,382 & \num{1.80e-05} & \num{2.00e-04} \\
			& min    &    21,525 & \ang{;  0; 5} & \ang{; 0;19} & \ang{;0; 3} &     94 & \num{5.85e-09} & \num{9.58e-08} \\
			& median &   457,358 & \ang{; 34;17} & \ang{; 8;38} & \ang{;0;40} &  4,586 & \num{6.31e-07} & \num{1.06e-05} \\
			& max    & 1,709,712 & \ang{;289;19} & \ang{;52;4 } & \ang{;6;56} & 15,677 & \num{1.87e-04} & \num{1.50e-03} \\
		\midrule
		\multirow{5}{*}{\rotatebox[origin=c]{90}{our$^\star$}}
			& mean   &   239,245 & \ang{;34;30} & \ang{;14;59} & \ang{;1;33} &  3,728 & \num{1.65e-05} & \num{2.88e-03} \\
			& std    &   178,979 & \ang{;62;19} & \ang{;12;34} & \ang{;1;15} &  3,787 & \num{5.57e-05} & \num{1.82e-02} \\
			& min    &     9,987 & \ang{;0;1}   & \ang{; 0;21} & \ang{;0; 3} &     61 & \num{4.09e-08} & \num{8.04e-07} \\
			& median &   189,880 & \ang{;9;13}  & \ang{;10;13} & \ang{;1;11} &  2,391 & \num{3.46e-06} & \num{2.62e-04} \\
			& max    & 1,033,492 & \ang{;324;8} & \ang{;55;11} & \ang{;7; 9} & 15,681 & \num{5.85e-04} & \num{2.30e-01} \\
		\bottomrule
	\end{tabular}
	}}

	\caption{
	    Dataset 3D pure hexahedra + star-shaped polyhedra (188 models in total). The memory is the total peak memory (in \si{\mebi\byte}) as reported by the solver Pardiso. \emph{$~^\star$does not include the models that went out of memory}.
	    From left to right, the total number of DOFs, the time required to solve the system, the time used to build the bases, the time employed to assembly the stiffness matrix, the peak memory, the $L_2$ error, and the $L_\infty$ error.
	}
	\label{tab:statistics}
\end{table}

\section{Limitations and concluding remarks}
\label{sec:conclusion}

We introduced Poly-Spline FEM, an integrated meshing and finite element method designed to take advantage of recent developments in hexahedral-dominant meshing, opening the doors to black box analysis with an high-order basis and cubic convergence under refinement. Our approach is to use the best possible basis for each element of the mesh and is amenable to continuous improvement, as the mesh generation methods and basis constructions improve. 
For instance, in this setting, one can avoid costly manual mesh repair and improvement, 
at the expense of modest increases in solution time, by switching to more expensive, but much less shape-sensitive elements when a hexahedron is badly shaped.

While our basis construction is resilient to bad element quality, the geometric map between the splines and the $Q_2$ elements might introduce distortion (and even inversions in pathological cases), lowering convergence rate. These effects could be ameliorated by optimizing the positions of the control points of the geometric map, which is an interesting avenue for future work. 

Our current construction always requires an initial refinement step to avoid having polyhedra adjacent to other polyhedra or to the boundary. This limitation could be lifted by generalizing our basis construction, and would allow our method to process very large datasets, that cannot be refined due to memory considerations. 
Another limitation of our method is that the consistency constraints in our basis construction (\Cref{sec:poly-contraints}) are PDE-dependent, and they thus require additional efforts to be used with a user-provided PDE: a small and reusable investment compared to the cost of manually meshing with hexahedra every surface that one wishes to analyze using $Q_2$ elements. 
The code can be found at \url{https://polyfem.github.io/} and 
provides an automatic way to generate such constraints relying on both the local assembler and automatic differentiation.

Poly-Spline FEM is a practical construction in-between unstructured $Q_2$ and fully-structured pure splines: it requires a smaller number of dofs than $Q_2$ (thanks to the spline elements) while preserving cubic convergence rate. We believe that our construction will stimulate additional research in the development of heterogeneous FEM methods that exploit the regularity of spline basis and combine it with the flexibility offered by traditional FEM elements. To allow other researchers and practitioners to immediately build upon our construction, we will release our entire software framework as an open-source project.

\section*{Acknowledgements}

We are grateful to the NYU HPC staff for providing computing cluster service. This work was partially supported by the \grantsponsor{NSFC}{NSF CAREER}{} award \grantnum{NSFC}{1652515}, the NSF grant \grantnum{NSFC}{IIS-1320635}, the NSF grant \grantnum{NSFC}{DMS-1436591}, the NSF grant \grantnum{NSFC}{1835712}, the SNSF grant \grantnum{SNSF}{P2TIP2\_175859}, a gift from Adobe Research, and a gift from nTopology.

\bibliographystyle{ACM-Reference-Format}
\bibliography{jabbrv,99-bib}

\appendix

\section{Brief Finite Element Introduction}
\label{app:fem}

Many common elliptic partial differential equations have the general form
\[
\mathcal{F}(\bx, u, \nabla u, \Delta u) = f(\bx), \qquad \bx\in\Omega,
\]
subject to
\[
u(\bx)=d(\bx),~\bx\in\partial\Omega_D \qquad\text{and}\qquad 
\nabla u(\bx) \cdot \vect{N}(\bx) = n(\bx),~\bx\in\partial\Omega_N
\]
where $\vect{N}(\bx)$ is the surface normal, $\partial\Omega_D$ is the Dirichlet boundary where the function $u$ is constrained (e.g., positional constraints) and $\partial\Omega_N$ is the Neumann boundary where the gradient of the function $u$ is constrained. The most common PDE in this class is the Poisson equation $-\Delta u = f$.

\paragraph{Weak Form} The first step in a finite element analysis consists of introducing the weak form of the PDE: find $u$ such that
\[
\int_\Omega\mathcal{F}(\bx, u, \nabla u, \Delta u)\, v(\bx)\,\dd\bx = \int_\Omega f(\bx)\, v(\bx)\,\dd\bx,
\]
holds for any \emph{test function} $v$ vanishing on the boundary. This reformulation has two advantages: (1) it simplifies the problem, and (2) it weakens the requirement on the function $u$. For instance, in case of the Poisson equation, the strong form is well defined only if $u$ is twice differentiable, which is a difficult condition to enforce on a discrete tesselation. However, the weak form requires only that the second derivatives of $u$ are integrable, allowing discontinuous jumps. Using integration by parts it can be further relaxed to
\[
\int_\Omega \nabla u(\bx) \cdot \nabla v(\bx)\,\dd\bx  = \int_\Omega f(\bx)\, v(\bx)\,\dd\bx,
\]
where only the gradient of $u$ needs to be integrable, that is $u\in H^1$, and can thus be represented using piecewise-linear basis functions. 

\begin{figure}
    \centering
    \includegraphics[width=0.22\linewidth]{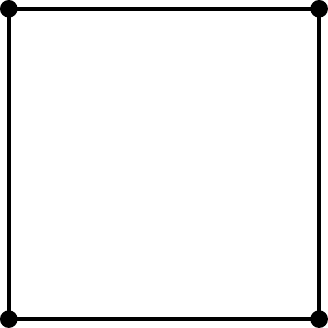}\hfill
    \includegraphics[width=0.22\linewidth]{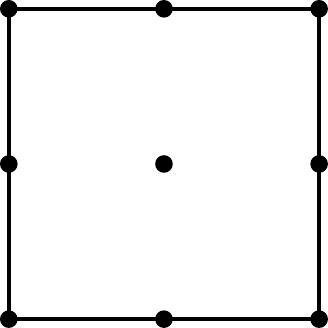}\hfill
    \includegraphics[width=0.22\linewidth]{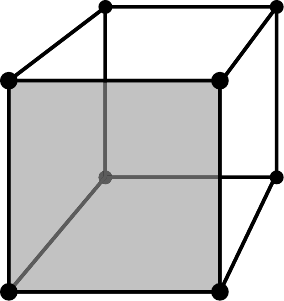}\hfill
    \includegraphics[width=0.22\linewidth]{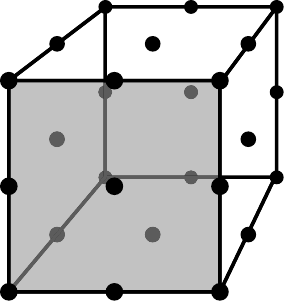}\par
    \parbox{0.22\linewidth}{\centering $Q_{1}$}\hfill
    \parbox{0.22\linewidth}{\centering $Q_{2}$}\hfill
    \parbox{0.22\linewidth}{\centering $Q_{1}$}\hfill
    \parbox{0.22\linewidth}{\centering $Q_{2}$}\par
    \caption{Node position for the linear and quadratic bases in two and three dimensions.}
    \label{fig:basis-nodes}
\end{figure}

\paragraph{Basis Functions} The key idea of a finite element discretization is to approximate the solution space via a \emph{finite} number of basis functions $\basis_i$, $i=1,\hdots, N$, which are independent from the PDE we are interested in.
The number of nodes (and basis functions) per element and their position is directly correlated to the order of the basis, see Figure~\ref{fig:basis-nodes}. We note that the nodes coincide the mesh vertices only for linear basis functions.
Instead of solving the PDE, the goal becomes finding the coefficients $u_i$, $i=1,\hdots, N$ of the discrete function $u_h(\bx)=\sum_{i=1}^N u_i\basis_i(\bx)$ that approximates the unknown function $u$. For a linear PDE this results in a linear system $\mat{K}\vect{u}=\vect{f}$, where $\mat{K}$ is the $N\times N$ stiffness matrix, $\vect{f}$ captures the boundary conditions, and $\vect{u}$ is the vector of unknown coefficients $u_i$. For instance, for the Laplace equation the entries of the stiffness matrix are
\[
K_{ij} = \int_\Omega \nabla\basis_i(\bx) \cdot \nabla\basis_j(\bx)\,\dd\bx.
\]

\paragraph{Local Support}Commonly used basis functions are locally supported. As a result, most of the pairwise intgerals are zero, leading to a sparse stiffness matrix. The pairwise integrals can be written as a sum of integrals over the elements (e.g., quads or hexes) on which both functions do not vanish. This representation enables so-called \emph{per-element assembly}: for a given element, a local stiffness matrix is assembled.

    For instance, if and element $C$ has four non-zero basis functions $\basis_i$, $\basis_j$, $\basis_k$, $\basis_l$ (this is the case for linear $Q_{1}$ quad) the local stiffness matrix $\mat{K}^L\in\reals^{4\times 4}$ for the Poisson equation is
\[
K^L_{o,p} = \int_C \nabla\basis_n(\bx) \cdot \nabla\basis_m(\bx)\,\dd\bx,
\]
where $o,p=1,\hdots,4$ and $m,n\in\{i,j,k,l\}$.
By using the mapping of local indices $(o,p)$ to global 
indices $(m,n)$, the local stiffness matrix entries are summed to yield the global stiffness matrix entries. 

\paragraph{Geometric Mapping} The final piece of a finite element discretization is the geometric mapping $\bg$. The local integrals need to be computed on every element. The element stiffness matrix entries are computed as integrals over a \emph{reference}  element $\hat{C}$ (e.g., a regular unit square/cube) through change of variables
\[
\int_C \nabla\basis_n(\bx) \cdot \nabla\basis_m(\bx) \,\dd\bx = 
\int_{\hat{C}} 
(\dD\bg^{\mathsf{-T}} \nabla\hat\basis_n(\bx)) \cdot 
(\dD\bg^{\mathsf{-T}} \nabla\hat\basis_m(\bx)) \abs{\dD\bg} \,\dd\bx,
\]
where $\dD\bg$ is the Jacobian matrix of the geometric mapping $\bg$, and $\hat\basis=\basis\circ \bg$ are the bases defined on the reference element $\hat{C}$. While usually $\bg$ is expressed by linear combination of $\basis_i$, leading to isoparametric elements, the choice of $\bg$ is independent from the basis.

\paragraph{Quadrature} All integrals are computed numerically by means of quadrature points and weights, which translates the integrals into weighted sums. Although there are many strategies to generate quadrature data (e.g., Gaussian quadrature), all of them integrate exactly polynomials up to a given degree to ensure an appropriate approximation order. For instance, if we use one quadrature point in the element's center with weight 1, we can integrate exactly constant functions.

\paragraph{Right-hand Side} The setup of the right-hand side $\vect{b}$ is done in a similar manner: its entries are $b_i = \int_\Omega \basis_i(\bx) f(\bx)\,\dd\bx$.

Dirichlet boundary conditions are treated as constrained degrees of freedom. The Neumann boundary conditions are imposed by setting
\[
b_j = \int_{\partial\Omega_N} \basis_j(\bx) \cdot n(\bx)\,\dd\bx
\]
for any node $j$ in $\partial\Omega_N$.

As for the stiffness matrix assembly the basis and node construction for the right-hand side is performed locally.

\section{Test Functions}
\label{app:franke}
In our experiments, we use the test functions proposed in \cite{Franke79}.
For 2D: 
\begin{align*}
f_{2\text{D}}(x_1,x_2) & = \frac{3}{4}{\rm e}^{-\frac{(9x_1-2)^2+(9x_2-2)^2}{4}}+\frac{3}{4}{\rm e}^{-\frac{(9x_1+1)^2}{49}-\frac{9x_2+1}{10}} \\
&+\frac{1}{2} {\rm e}^{-\frac{(9x_1-7)^2+(9x_2-3)^2}{4}}-\frac{1}{5} {\rm e}^{-(9x_1-4)^2-(9x_2-7)^2},
\end{align*}
and 3D:
\begin{align*}
f_{3\text{D}}&(x_1,x_2,x_3) = \\&\frac{3}{4}{\rm e}^{-\frac{(9x_1-2)^2+(9x_2-2)^2+(9x_3-2)^2}{4}}+\frac{3}{4} {\rm e}^{-\frac{(9x_1+1)^2}{49}-\frac{9x_2+1}{10}-\frac{9x_3+1}{10}} \\ 
&+\frac{1}{2} {\rm e}^{-\frac{(9x_1-7)^2+(9x_2-3)^2+(9x_3-5)^2}{4}}-\frac{1}{5} {\rm e}^{-(9x_1-4)^2-(9x_2-7)^2-(9x_3-5)^2}.
\end{align*}

\section{Basis functions}
\label{app:pk}

We use Lagrange tensor product function to interpolate between the nodes in quadrilateral and hexahedral elements. We provide the explicit formulation for $Q_{1}$ and $Q_{2}$ both in 2D, the 3D formulation follows. The four linear bases are constructed from the 1D linear bases
\[
\alpha_1(t) = 1-t \qquad \text{and}\qquad\alpha_2(t) =  t
\]
as the tensor products
\begin{align*}
\basis_1(u,v) = \alpha_1(u)\, \alpha_1(v),\qquad
\basis_2(u,v) = \alpha_1(u)\, \alpha_2(v),\\
\basis_3(u,v) = \alpha_2(u)\, \alpha_1(v),\qquad
\basis_4(u,v) = \alpha_2(u)\, \alpha_2(v).
\end{align*}
\noindent
Similarly the nine quadratic bases follow from the three quadratic polynomials
\[
\theta_1(t) = (1 - t) \, (1 - 2 t),\quad
\theta_2(t) = 4 t \, (1 - t), \quad
\theta_3(t) = t \, (2 t - 1)
\]
as
\begin{align*}
\basis_1(u,v) = \theta_1(u)\, \theta_1(v),\quad
\basis_2(u,v) = \theta_1(u)\, \theta_2(v),\quad
\basis_3(u,v) = \theta_1(u)\, \theta_3(v),\\
\basis_4(u,v) = \theta_2(u)\, \theta_1(v),\quad
\basis_5(u,v) = \theta_2(u)\, \theta_2(v),\quad
\basis_6(u,v) = \theta_2(u)\, \theta_3(v),\\
\basis_7(u,v) = \theta_3(u)\, \theta_1(v),\quad
\basis_8(u,v) = \theta_3(u)\, \theta_2(v),\quad
\basis_9(u,v) = \theta_3(u)\, \theta_3(v).
\end{align*}

\section{Polyhedral basis constraints}
\label{app:constraints}

We restrict the detailed explanation to 2D, the three-dimensional case follows.
Let $\vect{p}_i=(x_i, y_i)$, $i=1,\dots,s$, be the set of collocation points, that is, the points where we know the function values.
For the FEM basis $\basis_j$ that is nonzero on the polyhedral element $P$, we want to solve the least squares system $\vect{A w} = \vect{b}$, where

\begin{equation*}
    \mat{A} =
\begin{pmatrix}
    \kernel_1(\vect{p}_1)& \dots&   \kernel_k(\vect{p}_1)& 1&      x_1&    y_1&     x_1\, y_1& x_1^2&  y_1^2 \\
\vdots&      \ddots&  \vdots&      \vdots& \vdots& \vdots&  \vdots&    \vdots& \vdots&\\
\kernel_1(\vect{p}_s)& \dots&   \kernel_k(\vect{p}_s)& 1&      x_s&    y_s&     x_s\, y_s& x_s^2&   y_s^2
\end{pmatrix},
\end{equation*}
$\vect{b}$ is the evaluation of the basis $\basis_j$ from the neighbouring elements on the on the collocation points, and $\vect{w}=(w_1, \dots, w_k, a_{00}, a_{10}, a_{01}, a_{11}, a_{20}, a_{02})$.

Now to ensure consistency we need that
\[
-\int_{\bg(\param\cM)}\Delta \monom \basis_j = \int_{\bg(\param\cM)}\nabla \monom \cdot \nabla\basis_j
\]
holds for any of the 5 monomials.
We now split the previous integral over the polygon $P$ and over the known non-polygonal part $\overline P = \bg(\param\cM) \setminus P$
\[
\int_{P}\Delta \monom \basis_j + \int_{P}\nabla \monom  \cdot \nabla\basis_j =
-\int_{\overline P}\Delta \monom \basis_j - \int_{\overline P}\nabla \monom  \cdot \nabla\basis_j.
\]
We remark that the right-hand side of this equation is known since the bases on $\overline P$ are given, we call the five term $c_{ij}$ following the same indices as $a_{ij}$ (e.g., $c_{20} = 
-\int_{\overline P}\Delta x^2 \basis_j - \int_{\overline P}\nabla x^2  \cdot \nabla\basis_j
$).

We now evaluate the left-hand side for the five 2D monomials
\begin{align*}
    \int_{P}\Delta x  \basis_j +  \int_{P}\nabla x  \cdot  \nabla\basis_j   &= \int_{P}\pdiff{x}{\basis_j},\\
    \int_{P}\Delta y  \basis_j +  \int_{P}\nabla x  \cdot  \nabla\basis_j   &= \int_{P}\pdiff{y}{\basis_j},\\
    \int_{P}\Delta (xy) \basis_j +  \int_{P}\nabla (xy)  \cdot  \nabla\basis_j  &= \int_{P}y\pdiff{x}{\basis_j} + \int_{P}x\pdiff{y}{\basis_j},\\
    \int_{P}\Delta x^2 \basis_j +  \int_{P}\nabla x^2  \cdot  \nabla\basis_j  &= \int_{P} 2 \basis_j + \int_{P}2x\pdiff{x}{\basis_j},\\
    \int_{P}\Delta y^2 \basis_j +  \int_{P}\nabla y^2  \cdot  \nabla\basis_j  &= \int_{P} 2 \basis_j + \int_{P}2y\pdiff{y}{\basis_j}.
\end{align*}
By plugging the definition of $\basis_j$ over $P$ we obtain the following consistency constraints for the coefficients $a_{00}, a_{10}, a_{01}, a_{11}, a_{20}, a_{02}$:
\begin{multline*}
\sum_{i=1}^k w_i^j \int_P\pdiff{x}{\kernel} + a_{10}\abs{P} + a_{11}\int_P y + 2a_{20}\int_P x = c_{10}\\
\sum_{i=1}^k w_i^j \int_P\pdiff{y}{\kernel} + a_{01}\abs{P} +\\ a_{11}\int_P x + 2a_{02}\int_P y= c_{01},\\
\sum_{i=1}^k w_i^j \int_P(y\pdiff{x}{\kernel} + x\pdiff{y}{\kernel}) +\\ a_{10} \int_P x + a_{01}\int_P y + a_{11}\int_P x^2 + y^2 + 2(a_{20} + a_{02})\int_P xy= c_{11}\\
2\sum_{i=1}^k w_i^j \int_P(\kernel + x\pdiff{x}{\kernel})  + 2a_{00} +\\ 4a_{10}\int_P x + 2a_{01}\int_P y + 4a_{11}\int_P x y + 6a_{20}\int_P x^2 + 2a_{02}\int y^2= c_{20},\\
2\sum_{i=1}^k w_i^j \int_P(\kernel + y\pdiff{y}{\kernel}) + 2a_{00} +\\ 2a_{10}\int_P x + 4a_{01}\int_P y + 4a_{11}\int_P x y + 2a_{20}\int_P x^2 + 6a_{02}\int y^2= c_{02}.\\
\end{multline*}

\end{document}